\documentclass[11pt]{amsart}

\usepackage{amssymb,array}
\usepackage[pdftex]{graphicx,epsfig,color} 
\usepackage[all]{xy} 

\def\ov#1{{\overline{#1}}}
\def\un#1{{\underline{#1}}}
\def\wh#1{{\widehat{#1}}}
\def\wt#1{{\widetilde{#1}}}

\newcommand{\MV}{\operatorname{MV}}
\newcommand{\ord}{\operatorname{ord}}

\newcommand{\MI}{\operatorname{MI}}

\newcommand{\Conv}{\operatorname{Conv}}
\newcommand{\NP}{\operatorname{NP}}

\newcommand{\Res}{\operatorname{Res}}
\renewcommand{\div}{\operatorname{div}}
\newcommand{\Graph}{\operatorname{Graph}}
\newcommand{\Vol}{\operatorname{Vol}}

\renewcommand{\Im}{\operatorname{Im}}
\newcommand{\rank}{\operatorname{rank}}
\newcommand{\init}{\operatorname{init}}

\newcommand{\hot}{\mbox{ h.o.t. }}
\newcommand{\Faces}{\operatorname{Faces}}
\newcommand{\Slope}{\operatorname{Slopes}}
\newcommand{\Card}{\operatorname{Card}}
\newcommand{\mult}{\operatorname{mult}}

\newcommand{\geom}{{{\boldsymbol{\cA}}}}
\newcommand{\arith}{{\boldsymbol{\widehat{\cA}}}}
\newcommand{\I}{{\operatorname{I}}}
\newcommand{\Ef}{{\operatorname{F}}}

\renewcommand{\d}{{\operatorname{d}}}
\newcommand{\dist}{\operatorname{dist}}
\newcommand{\Perm}{\operatorname{Perm}}
\newcommand{\codim}{\operatorname{codim}}

\def \A{\mathbb{A}}
\def \C{\mathbb{C}}

\def \K{\mathbb{K}}

\def \N{\mathbb{N}}
\def \P{\mathbb{P}}
\def \Q{\mathbb{Q}}
\def \R{\mathbb{R}}
\def \Es{\mathbb{S}}
\def \T{\mathbb{T}}
\def \Z{\mathbb{Z}}

\def \k{\mathbb{K}}

\def\cA {{\mathcal A}}
\def\cB {{\mathcal B}}
\def\cC {{\mathcal C}}

\def\cE {{\mathcal E}}

\def\cO {{\mathcal O}}

\def\cS {{\mathcal S}}

\newcommand{\bfb}{{\boldsymbol{b}}}
\newcommand{\bfc}{{\boldsymbol{c}}}
\newcommand{\bfd}{{\boldsymbol{d}}}
\newcommand{\bff}{{\boldsymbol{f}}}
\newcommand{\bfg}{{\boldsymbol{g}}}
\newcommand{\bfell}{{\boldsymbol{\ell}}}
\newcommand{\bfs}{{\boldsymbol{s}}}
\newcommand{\bft}{{\boldsymbol{t}}}
\newcommand{\bfu}{{\boldsymbol{u}}}

\newcommand{\bfx}{{\boldsymbol{x}}}
\newcommand{\bfy}{{\boldsymbol{y}}}

\newcommand{\bfD}{\boldsymbol{D}}
\newcommand{\bfF}{\boldsymbol{F}}
\newcommand{\bfP}{\boldsymbol{P}}
\newcommand{\bfQ}{\boldsymbol{Q}}
\newcommand{\bfU}{\boldsymbol{U}}
\newcommand{\bfV}{\boldsymbol{V}}

\newcommand{\bfcalB}{{\boldsymbol{\cB}}}

\newcommand{\bfalpha}{{\boldsymbol{\alpha}}}
\newcommand{\bfbeta}{{\boldsymbol{\beta}}}
\newcommand{\bfdelta}{{\boldsymbol{\delta}}}
\newcommand{\bfvartheta}{{\boldsymbol{\vartheta}}}
\newcommand{\bfrho}{{\boldsymbol{\rho}}}

\newcommand{\bftau}{{\boldsymbol{\tau}}}

\newcommand{\bfzero}{\boldsymbol{0}}
\newcommand{\bfun}{\boldsymbol{1}}

\numberwithin{equation}{section}

\newtheorem{defn}{Definition}[section]
\newtheorem{lem}[defn]{Lemma}

\newtheorem{prop}[defn]{Proposition}
\newtheorem{thm}[defn]{Theorem}
\newtheorem{cor}[defn]{Corollary}

\newtheorem{rem}[defn]{Remark}
\newtheorem{exmpl}[defn]{Example}
\newtheorem{example}{Example}[subsection]

           {\vspace{3.3mm}
           \noindent{\bf #1}\it}%
           {\vspace{3.3mm}}

\newenvironment{Proof}[1]{\noindent{\it #1}}{\hfill\mbox{$\Box$} }

\font\ninebf=cmbxti10 at 8pt
\def\enmarge#1{\goodbreak\ifodd\count0\marginpar{\vbox{\hsize=28mm\ninebf\char62#1}}
\else\marginpar{\vbox{\hsize=33mm\leftskip5mm\ninebf\char62#1}}\fi}
      
\begin{document}

\title[A refinement of the Ku\v{s}nirenko--Bern\v{s}tein estimate]
{A refinement of the Ku\v{s}nirenko--Bern\v{s}tein estimate}

\author[Patrice Philippon]{Patrice Philippon} 

\address{Institut de Math{\'e}matiques de 
Jussieu - U.M.R. 7586 du CNRS, \'Equipe G{\'e}om{\'e}trie et Dynamique.
175 rue du Chevaleret, 75013 Paris, France.}
\email{pph@math.jussieu.fr}
\urladdr{http://www.math.jussieu.fr/$\scriptstyle\sim$pph/}


\author[Mart{\'\i}n~Sombra]{Mart{\'\i}n~Sombra}

\address{Universitat de Barcelona,
Departament d'{\`A}lgebra i Geometria.
Gran Via 585,
08007 Barcelona, Spain.}
\email{sombra@ub.edu}
\urladdr{http://atlas.mat.ub.es/personals/sombra/}
\thanks{M. Sombra was supported by the 
Ram\'on y Cajal Program of the Ministerio de Educaci\'on y Ciencia, Spain.}


\date{\today}


\begin{abstract}
\setlength{\baselineskip}{12pt}
A theorem of Ku\v snirenko and Bern\v stein shows that the number of iso\-lated roots of a system of polynomials in a torus is bounded above by the mixed volume of the Newton polytopes of the given polynomials, and this upper bound is generically  exact. 
We improve on this result by introducing refined combinatorial invariants of polynomials and a generalization of the mixed volume of convex bodies: the mixed integral of concave functions. 
The proof is based on new techniques and results from relative toric geometry.
\end{abstract}

\maketitle

\overfullrule=1mm

\vspace{-6mm} 

\setcounter{tocdepth}{1}\tableofcontents

\vspace{-8mm} 
\setcounter{section}{0}
\section{Introduction}

The main purpose of this text is to establish a new upper bound for the number of roots of a system of polynomial equations. More precisely, let $\k$ be an algebraically closed field and consider a family of $n+1$ Laurent polynomials
$f_0,\dots,f_n \in \k[s][t_1^{\pm 1},\dots,t_n^{\pm 1}]$ in the variables $\bft=(t_1,\dots,t_n)$ with coefficients polynomials in the single variable $s$. How many isolated solutions $\xi\in\k\times(\k^\times)^{n}$ are there to the system of equations 
$$
f_0(\xi)=\cdots=f_n(\xi)=0 
\enspace?
$$
Let $P_i\subset \R^{n+1}$ denote the Newton polytope of $f_i$ when regarded as a Laurent polynomial in all of the variables $s,t_1,\dots,t_n$. The classical theorem of Ku\v{s}nirenko and Bern\v{s}tein asserts that the number (counting multiplicities) of those isolated points lying in $(\k^\times)^{n+1}$ is bounded above by the mixed volume $\MV_{n+1}(P_0,\dots,P_n)$,  with equality when 
$f_0,\dots,f_n$ is generic among those systems with Newton polytopes $P_0,\dots,P_n$ \cite{Kus76,Ber75}. This result is a cornerstone of toric geometry and polynomial equation solving, {\it see} for instance~\cite{GKZ94,Stu02}.

\medskip
For the formulation of our results we introduce some combinatorial invariants associated to the given system of polynomials. Let $f=\sum_{j=0}^N\alpha_j(s)\bft^{a_j}\in\k(s)[t_1^{\pm 1},\dots,t_n^{\pm 1}]$ be a non-zero Laurent polynomial and for each $v\in \P^1$ consider the {\it $v$-adic Newton polytope of $f$} defined as the convex hull
$$
\NP_v(f):=\Conv\big( (a_0,-\ord_v(\alpha_0)), \dots,(a_N,-\ord_v(\alpha_N))\big) \subset \R^{n+1}\enspace,
$$
where $\ord_v(\alpha_j)$ denotes the {\em order of vanishing of $\alpha_j$ at $v$} viewed as a rational function on $\P^1$. This polytope sits above the Newton polytope relative to the variables $\bft$
$$\NP(f):=\Conv(a_0,\dots, a_N)\subset \R^n$$
{\it via} the natural projection $\R^{n+1}\to \R^n$ that forgets the last coordinate. 
Consider the {\em roof function of $f$ at $v$} defined as 
$$
\vartheta_v(f):\NP(f)\rightarrow\R\quad,\quad {u}\mapsto\max\{z\in\R:({u},z)\in\NP_v(f)\}
\enspace,$$
that is the concave and piecewise affine function parameterizing the upper envelope of $\NP_v(f)$ above $\NP(f)$. 
For $f=0$ we set for convenience 
$\NP(f):=\{0\}\subset \R^n$ and for $v\in \P^1$ 
we define $\vartheta_v(f):\{0\}\to \R$ to be the zero function. It is worth mentioning that this roof function appears also in tropical geometry as the {Legendre-Fenchel dual} of the {``tropical polynomial''}: $\R^n \to \R, u\mapsto \min_j(\ord_v(\alpha_{j})+\langle a_{j},u\rangle)$, associated to $f$ with respect to the  valuation $\ord_v$, {\it see} for instance~\cite{Mik04}. 

\smallskip
For concave functions $\rho:Q\rightarrow\R$ and $\sigma: R\rightarrow\R$ defined on convex sets $Q,R\subset\R^n$ respectively, we consider their {\it sup-convolution}
$$\rho\boxplus\sigma:Q+R\rightarrow\R\quad,\quad {u}\mapsto\max\{\rho({v})+\sigma({w}):{v}\in Q,{w}\in R,{v}+{w}={u}\}
\enspace,$$
which is a concave function defined on the Minkowski sum $Q+R$. This operation is dual under the Legendre-Fenchel conjugation to the pointwise sum of concave functions~\cite{Roc97} (whence the name ``convolution''\footnote{This notion comes from convex analysis, but since this theory deals mostly with convex functions rather than with concave, the corresponding operation of {inf-convolution} (usually denoted $\Box$) is more common in this context; the connection with our notation reads: $-(\rho\boxplus\sigma)=(-\rho)\Box(-\sigma)$.}) and extends the Minkowski sum to concave functions. 

\begin{defn}[\cite{PS03}] \label{defmultiintegrale}
For a family of $n+1$ concave functions $\rho_0:Q_0\rightarrow\R$, \dots, $\rho_n:Q_n\rightarrow\R$ defined on {convex bodies} of $\R^n$, the {\em mixed integral} is 
$$\MI_n(\rho_0,\dots,\rho_n):=\sum_{j=0}^n(-1)^{n-j}\sum_{0\leq i_0<\dots<i_j\leq n}\int_{Q_{i_0}+\dots+Q_{i_j}}(\rho_{i_0}\boxplus\dots\boxplus\rho_{i_j})({u}) \, 
\d u_1\cdots \d u_n\enspace.$$
\end{defn}

This is the natural extension to concave functions of the mixed volume of convex bodies and as such it satisfies analogous properties: it is symmetric in $\rho_0,\dots,\rho_n$, linear with respect to $\boxplus$ in each variable $\rho_i$, and for a concave function $\rho:Q\to \R$ we have: $\MI_n(\rho, \dots, \rho) = (n+1)! \, \int_Q \rho(u) \, \d u_1\cdots \d u_n$~\cite[Prop.~IV.5(a,b)]{PS03}. In~\S~\ref{bk_multiint} we  establish further properties of this notion, in particular its monotonicity (proposition~\ref{fantine}) and a decomposition formula (proposition~\ref{marius}) expressing the mixed integral in terms of lower dimensional mixed integrals and mixed volumes, analogous to the decomposition formula for mixed volumes, {\it see}~\cite[Thm.IV.4.10]{Ewa96}. 

\smallskip
For a system of Laurent polynomials $f_0,\dots,f_n\in\k[s][\bft^{\pm 1}]$ let $V(f_0,\dots,f_n)\subset \k\times(\k^\times)^{n}$ denote the set of solutions of $f_0=\cdots=f_n=0$ and $V(f_0,\dots,f_n)_0$ the subset of those solutions 
that are isolated. For each of those isolated points $\xi$ we denote by ${\rm mult}(\xi|f_0,\dots,f_n)$ the intersection multiplicity of $f_0,\dots, f_n$ at $\xi$, {\it see} formula~(\ref{mult_intro}). We say that $f_i$ is {\it primitive} if it has no non-constant factor in $\K[s]$. The following is our first main result. 

\begin{thm} \label{mainthm}
Let $f_0,\dots,f_n\in\K[s][\bft^{\pm1}]=\k[s][t_1^{\pm 1},\dots,t_n^{\pm 1}]$ be a family of primitive Laurent polynomials in the variables $\bft$ with coefficients in $\k[s]$. For $0\le i\le n$ and $v\in\P^1$ let $\vartheta_{i,v}:\NP(f_i)\to \R$ be the roof function of $f$ at $v$, then 
\begin{equation} \label{mainthm_display}
\sum_{\xi\in V(f_0,\dots,f_n)_0}{\rm mult}(\xi|f_0,\dots,f_n) 
\le \sum_{v\in \P^1} \MI_n(\vartheta_{0,v},\dots,\vartheta_{n,v}) \enspace.
\end{equation}
Furthermore, this is an equality for $f_0,\dots,f_n$ generic among systems with given functions $(\vartheta_{i,v}: 0\le i\le n, v\in \P^1)$.
\end{thm}

Specializing to the unmixed case we obtain the following estimate. 

\begin{cor} \label{cor_mainthm}
With notation as in the above theorem, let $Q\subset \R^n$ be a polytope containing $\NP(f_i)$ for all $i$, and for each  $v\in \P^1$ let $\vartheta_v : Q\to \R$ be a concave function such that $\vartheta_{v}\ge \vartheta_{i,v}$ for all $i$. Then
\begin{equation} \label{cor_mainthm_display}
\sum_{\xi\in V(f_0,\dots,f_n)_0}{\rm mult}(\xi|f_0,\dots,f_n) \le (n+1)!\, \sum_{v\in \P^1} \int_Q \vartheta_v(u) \, \d u_1 \cdots \d u_n \enspace.
\end{equation}
\end{cor}

\medskip
As an illustration, consider for $k\ge 1$ the polynomials 
$$
f=(s-1)^{2k}+(s-1)^k t-s t^{2} \quad,\quad g=-3(s-1)^{2k}+(s-1)^k t+s t^{2}  \quad\in\K[s,t]\enspace.
$$
The system $f=g=0$ has the  only solution $(2,1)$ in $\K\times\K^\times$. Standard and bihomogeneous B\'ezout theorem give respectively the upper bounds $\deg(f)\deg(g)=4k^2$ and $\deg_s(f)\deg_t(g)+\deg_t(f)\deg_s(g)= 8k$ for the number  of isolated roots, while the Ku\v{s}nirenko-Bern\v{s}tein theorem predicts at most 
$4k+1$ roots in $(\K^\times)^2$. On the other hand, corollary~\ref{cor_mainthm} 
gives the exact estimate $1$,
hence this system is generic with respect to our estimate but not with 
respect to Ku\v{s}nirenko-Bern\v{s}tein's one, {\it see} example~\ref{exmpl 2} for the details. 

\smallskip

In the setting of theorem~\ref{mainthm}, note that $\vartheta_{i,v}=0$ for almost all $v$ and so the number of non-zero terms in the right hand side of
estimate~(\ref{mainthm_display}) is finite. The only positive contribution in this sum comes 
from the place $v=\infty$, 
because $-\ord_{\infty}(\alpha_{i,j})=\deg(\alpha_{i,j})$ and therefore $\vartheta_{i,\infty}\ge 0$ while for $v\in \P^1\setminus \{\infty\}$ we have $-\ord_v(\alpha_{i,j}) \le 0$ and therefore $\vartheta_{i,v}\le 0$, together with
the monotonicity of the mixed integral.  

The function $\vartheta_{i,\infty}$ ({\it resp.} $-\vartheta_{i,0}$) parameterizes the upper ({\it resp.} lower) envelope of $P_i\subset\R^{n+1}$, 
the Newton polytope of $f_i$ with respect to all of the variables $s$ and $\bft$, and by proposition~\ref{bagne}
$$\MI_n(\vartheta_{0,0},\dots,\vartheta_{n,0}) + \MI_n(\vartheta_{0,\infty},\dots,\vartheta_{n,\infty}) = \MV_{n+1}( P_0, \dots,  P_n)
\enspace.$$
This shows that~(\ref{mainthm_display}) improves upon Ku\v snirenko-Bern\v{s}tein's estimate in the case of primitive $f_i$'s in $ \K[s][\bft^{\pm1}]$, besides the fact that it counts the isolated roots of the system in a set larger than $(\k^\times)^{n+1}$. 

A discrepancy between both estimates will actually occur when 
at least one of the mixed integrals in~(\ref{mainthm_display}) 
corresponding to a place $v\ne 0, \infty$ is strictly negative. 
This might happen when some of the coefficients $\alpha_{i,j}$ share common zeros, as in the example above, though of course the amount of improvement depends  on the exact configuration of the $\vartheta_{i,v}$'s. These remarks extend to general Laurent polynomials  in $\k[s^{\pm1},\bft^{\pm1}]$, 
{\it see} inequality~(\ref{comparisonwithBK}). 

\medskip
The conditions for the estimate in theorem~\ref{mainthm} to be exact
can be specified in terms of lower dimensional systems of equations. Let $f\in \K(s)[\bft^{\pm1}]$ and  $\tau\in \R^n$, for $v\in \P^1\setminus \{\infty\}$ the {\it $\tau$-initial part of $f$ at $v$} is the Laurent polynomial $\init_{v,\tau}(f)\in\K[\bft^{\pm1}]\setminus \{0\}$ such that
$$
f(s, s^{-\tau_1}t_1,\dots,s^{-\tau_n}t_n) = (s-v)^c (\init_{v,\tau}(f)(\bft)  + o(1)) 
$$ 
for some $c\in\Z$ and $o(1)$ going to $0$ as $s$ tends to $v$, while the {\em $\tau$-initial part of $f$ at $\infty$} is just defined as the $\tau$-initial part of $f(s^{-1},\bft)$ at $0$.

\begin{prop}\label{equality in mainthm}
With notation as in theorem~\ref{mainthm}, if 
for all $v\in \P^1\setminus \{\infty\}$ and $\tau \ne \bfzero$, and for 
$v=\infty$ and all $\tau\in \R^n$, 
the system of equations
\begin{equation}\label{eponine}
\init_{v,\tau}(f_0)(\xi) = \cdots= \init_{v,\tau}(f_n)(\xi)=0 
\end{equation}
has no solution 
in $(\K^\times)^n$, then the estimate~(\ref{mainthm_display}) is an equality. 
\end{prop}

Though it is not evident from the formulation above, these genericity conditions
are equivalent to a {\it finite} number of systems of equations in $\le n$ variables, {\it see}~\S~\ref{bk_equality}. A further situation where we find equality in~(\ref{mainthm_display}) is when $\sum_{v\in \P^1} \MI_n(\vartheta_{0,v},\dots,\vartheta_{n,v}) =0$, a condition that can be characterized in terms of the rank of some $\Z$-modules (proposition~\ref{anzelma}). In this case, the estimate is obviously an equality for any system with given functions $(\vartheta_{i,v}: 0\le i\le n, v\in \P^1)$. 

\bigskip

It is natural to try to extend theorem~\ref{mainthm} to an arbitrary base (of dimension $1$) instead of $\P^1$. In this direction, we consider the case of a smooth complete curve $S$ equipped with a family of line bundles $L_0,\dots, L_n$. For $0\le i\le n$ consider the $\K$-vector space of Laurent polynomials with coefficients global sections of $L_i$
$$
\Gamma(S; L_i)[\bft^{\pm1}] := \Gamma(S;L_i)\otimes_\K\K[\bft^{\pm1}] = \bigoplus_{a\in\Z^n} \Gamma(S;L_i)\otimes_\K \bft^{a} 
$$
together with an element $f_i=\sum_{i=0}^{N_i}\sigma_{i,j} \bft^{a_{i,j}} \in\Gamma(S;L_i)[\bft^{\pm1}]$. The set of zeros ({\it resp.} isolated zeros)  in $S\times(\K^\times)^n$ of the system $f_0=\dots=f_n=0$ is well-defined and as before we denote it by $V(f_0,\dots,f_n)$ or $V(\bff)$ ({\it resp.}  $V(f_0,\dots,f_n)_0$ or $V(\bff)_0$). We also extend in the natural way the notions of $v$-adic Newton polytope $\NP_v(f_i)\subset \R^{n+1}$ and corresponding roof function $\vartheta_v(f_i):\NP(f_i)\to\R$.

To take into account the possibility that the coefficients of some $f_i$ might have common zeros, we introduce for each $0\le i\le n$ and $v\in S$  an extra function $\overline{\vartheta}_{v}(f_i)$ defined as the constant function $\ord_v(f_i):=\min_j \ord_v(\sigma_{i,j})$ on the Newton polytope of the Laurent polynomial $f_i(v,\cdot)\in\K[\bft^{\pm1}]$. And for a set $Q$ and  $c\in\R$ we denote by $c|_Q$ the constant function $c$ with domain $Q$. 

\smallskip

In this general setting we have the following extension of theorem~\ref{mainthm}. 

\begin{thm} \label{genthm}
Let $S$ be a smooth complete curve equipped with line bundles $L_i$ for $0\le i\le n$, and for each $i$ let $f_i \in \Gamma(S;L_i)[\bft^{\pm 1}]\setminus\{0\}$ be a non-zero  Laurent polynomial in the $\bft$-variables with coefficients in $\Gamma(S;L_i)$.
 
Let $\delta_i:=\deg(L_i)$ and set for short $\bfdelta|_{\NP(\bff)}:=\left(\delta_0|_{\NP(f_0)},\dots,\delta_n|_{\NP(f_n)}\right)$ and $\bff:=(f_0,\dots,f_n)$. For $v\in S$ let $\vartheta_{v}(f_i):\NP(f_i)\to \R$ denote the  roof function of $f$ at $v$ and $\overline{\vartheta}_{v}(f_i):\NP(f_i(v,\cdot))\to \R$ the constant function $\ord_v(f_i)$, then
\begin{equation}\label{genthm_display}
\sum_{\xi\in V(\bff)_0}{\rm mult}(\xi|\bff) \le 
\MI_{n}\left(\bfdelta|_{\NP(\bff)} \right)+
\sum_{v\in S} \left(\MI_{n} (\vartheta_v(\bff)) + \MI_{n}(\overline{\vartheta}_{v}(\bff)) \right)\enspace.
\end{equation}
Furthermore, this is an equality for $\bff$  generic among systems with given functions
$(\vartheta_{i,v}: 0\le i\le n, v\in S)$.
\end{thm}

This result is only significant when $\deg(L_i)\ge 0$ for all $i$, 
otherwise $L_i$ does not admit 
non-zero global sections.
On the other hand, the smoothness hypothesis is not strictly necessary, 
and the result can be extended to a singular base curve (theorem~\ref{gengenthm}). 
As for theorem~\ref{mainthm}, it is possible to give explicit conditions 
for equality in the  estimate~(\ref{genthm_display}) in terms of lower dimensional systems of equations (proposition~\ref{genegalite}).

For all $v\in S$ we have $\MI_{n} (\vartheta_v(\bff)) + \MI_{n}(\overline{\vartheta}_{v}(\bff))$ $
\le 0$
because of 
$\NP(f_i(v,\cdot))\subset \NP(f_i)$, 
$\vartheta_v(f_i)\le -\ord_v(f_i) \le 0$, 
and the monotonicity of the mixed integral. 
Hence the only positive contribution in the right hand side of~(\ref{genthm_display}) comes from the first term. Besides, we show that $ \MI_{n}(\overline{\vartheta}_{v}(\bff))\ne 0$ if and only if $v$ is a base point of exactly one  of the $f_i$'s and in that case, this mixed integral can be expressed as a $n$-dimensional mixed volume, {\it see} remark~\ref{remarqueMI}. In particular, when the $f_i$'s have no base points, the functions $\ov \vartheta_{i,v}$ do not contribute to the estimate at all. This is precisely the situation in theorem~\ref{mainthm} because of the assumption that the $f_i$'s are primitive. Indeed, applying theorem~\ref{genthm} to the case $S=\P^1$ extends theorem~\ref{mainthm} to possibly non-primitive polynomials, {\it see}~inequality~(\ref{mainthm_whithbasepoints}).

The present generalization of theorem~\ref{mainthm} 
to an arbitrary base curve allows to treat systems of equations
over a semi-abelian variety $G$, extension $0\to(\K^\times)^n\to G\to E\to 0$  of an elliptic curve $E$ by a torus, {\it see} example~\ref{exempl 4} for the case of a torus of dimension $1$. 
We refer to~\S~\ref{Comparing} for another kind of situation that 
is out of reach of 
theorem~\ref{mainthm} but can be sucessfully treated with  
theorem~\ref{genthm}.

\medskip 

As a consequence of theorem~\ref{genthm}, we obtain a bound for the degree of cycles of $S\times(\K^\times)^n$ of positive dimension $d$, defined by $n+1-d$ equations. For $1\le i\le d$ let $L_i$ be a line bundle over $S$ and $g_i=\sum_{j=0}^{N}\sigma_{i,j}\bft^{a_{i,j}} \in \Gamma(S;L_i)[\bft^{\pm 1}]\setminus\{0\}$. The {\it degree with respect to  $g_1,\dots, g_d$} of a pure $d$-dimensional cycle $Z\subset S\times(\K^\times)^n$ is defined as 
$$
\deg_{\bfg}(Z) := \deg\big(Z\cdot \varphi_{g_1}^{-1}(E_1)\cdot \, \cdots \, 
\cdot \varphi_{g_d}^{-1}(E_{d})\big) 
\enspace,$$
where $\cdot$ denotes the intersection product, $\varphi_{g_i}$ denotes the map
$$
S\times(\K^\times)^n \to 
\P^{N_i} \quad ,\quad 
(s,\bft) \mapsto (\sigma_{i,0}(s)\bft^{a_{i,0}}: \cdots : \sigma_{i,N_i}(s)\bft^{a_{i,N_i}})\enspace, 
$$
and $E_i$ is a generic hyperplane of $\P^{N_i}$.

\begin{cor}\label{corgenthm}
Let $S$ be a smooth complete curve equipped with line bundles $L_i$ for $0\le i\le n$. Let  $m\leq n$, for $0\le i\le m$ and $m+1\le k\le n$ let be given  $f_i
\in\Gamma(S;L_i)[\bft^{\pm 1}]$  and $g_k\in\Gamma(S;L_k)[\bft^{\pm 1}]$ such that 
the base locus of each $g_k$ is empty. 
 
Let $\delta_i:=\deg(L_i)$ and set for short $\bff:=(f_0,\dots, f_m)$, $\bfg:=(g_{m+1},\dots, g_n)$, then $\bfdelta|_{\NP(\bff)}:=(\delta_0|_{\NP(f_0)},\dots,\delta_m|_{\NP(f_m)})$ and $\bfdelta|_{\NP(\bfg)}:=(\delta_{m+1}|_{\NP(g_{m+1})},\dots,\delta_n|_{\NP(g_n)})$. Let $Z(\bff)_{n-m}$ denote the $(n-m)$-dimensional part of the intersection cycle $\div(f_0)\cdot\, \cdots \, \cdot \div(f_m)$ in 
$S\times(\K^\times)^n$, then $\deg_{\bfg}(Z(\bff)_{n-m})$ is bounded above by
$$ \MI_n(\bfdelta|_{\NP(\bff)}, \bfdelta|_{\NP(\bfg)}) + \sum_{v\in S} \big(\MI_n(\vartheta_{v}(\bff),\vartheta_{v}(\bfg)) + \MI_n(\overline{\vartheta}_{v}(\bff),\overline{\vartheta}_{v}(\bfg)) \big)
\enspace.$$
\end{cor} 

\bigskip
The proof of these results is based on intersection theory applied to a suitable compactification  of the ambient space $S\times (\K^\times)^n$, 
{\it see} sections~\ref{bk_degree} and~\ref{bk_bezout}. The system of Laurent polynomials $f_i\in\Gamma(S;L_i)[\bft^{\pm1}]$ is naturally associated to a linear system on a multiprojective toric variety $X$ over the curve $S$. 
These varieties are related to the toric varieties over a discrete valuation ring studied by A.L.~Smirnov in a similar context~\cite{Smi97}. 
Sections~\ref{The geometry of a toric variety over a curve}
and~\ref{orbitdecomposition} are devoted to a thoughtful study of 
toric varieties over a curve, and in particular we show that such a 
variety 
is naturally endowed with a family of concave piecewise affine functions $\Theta_X:=(\vartheta_{i,v}:\NP(f_i)\to\R: 0\le i\le n, v\in S)$ which plays for this variety the r\^ole of the polytope for a
projective toric variety over a field.

We show that the geometry of $X$ can be made explicit  in terms of $\Theta_X$, in particular its dimension, fiber structure and mixed degrees. The estimate for the number of roots is deduced from the computation of a certain mixed degree, while the genericity conditions are obtained from a fine study of the structure of the fibers over $S$. This strategy is reminiscent of B.~Teissier's approach to the  Ku\v snirenko-Bern\v stein theorem that is implicit in~\cite{Tei79}, {\it see} also~\cite[chap. 5]{Ful93}.

Other approaches to the  Ku\v snirenko-Bern\v stein theorem might extend to the setting of theorems~\ref{mainthm} and~\ref{genthm}. Based on a preliminary version of this text, M.I.~Herrero has recently proposed an alternative proof of theorem~\ref{mainthm} for the case of bivariate polynomials, close in spirit to Bern\v stein's original article~\cite{Her07}.

\bigskip
In practical situations, the computation of the  estimate~(\ref{mainthm_display})
can be substantially simplified by some observations. The functions $\vartheta_{i,v}$ can be directly obtained from factorizations 
$$
\alpha_{i,j}(s)= \lambda_{i,j} \prod_{p\in P} p(s)^{e_p(i,j)} \quad \mbox{ for } 0\le i\le n \mbox{ and } 0\le j\le N_i
$$ 
for some finite set $P\subset \K[s]$ of pairwise coprime polynomials, $e_p(i,j)\in \N$ and $\lambda_{i,j}\in \K^\times$ (proposition~\ref{mainthm_calculable}) and 
such factorizations can be computed with gcd computations only (lemma~\ref{factorization}). For a system $\bff$ defined over an effective field (say $\Q$)
the $\vartheta_{i,v}$'s can thus be determined with no need to access to the roots of the $\alpha_{i,j}$'s,  not even to completely factorize them over $\Q[s]$. In addition, the relevant mixed integrals can be calculated by applying the decomposition formula in proposition~\ref{marius}, thus avoiding the costly computation of sup-convolutions.

\smallskip
We close this introduction by pointing out a recent application of theorem~\ref{mainthm}, to the determination of the Newton polygon of the equation of a rational plane curve in terms of a given parameterization~\cite{DS07}.


\bigskip 

\noindent{\bf Acknowledgments.\hspace{1.5mm}---}\hspace{1mm} 
We thank Bernard Teissier for clarifying discussions on the notion of multiplicity.

\bigskip

\section{The geometry of a toric variety over a curve}
\label{The geometry of a toric variety over a curve}

Let $S$ be a curve defined over an algebraically closed field $\K$, with field of $\K$-rational functions  $\K(S)$; we assume that $S$ is complete and smooth unless otherwise stated. Let $\T^d:=(\K^\times)^d$ be the algebraic torus of dimension $d$ with coordinates $\bft=(t_1,\dots,t_d)$ and for some positive integers $N_0,\dots,N_n$ let $\P:=\P^{N_0}\times\cdots\times\P^{N_n}$ be the corresponding multiprojective space, with coordinates $\bfx_i=(x_{i,0}:\dots:x_{i,N_i})$ for $0\le i\le n$. A variety is supposed to be defined over $\K$, reduced and irreducible.  
For a cycle $Z$ we denote by $|Z|$ its underlying algebraic set. A property depending on parameters is said {generic} if it holds for all points in a dense 
open subset of the  parameter space. We denote by $\N$ the set of  all natural integers including $0$. 

\bigskip

The present and next sections are devoted to the study of the structure of multiprojective toric varieties over $S$. The necessary background on toric varieties over a field can be found for instance, in~\cite{Ful93,GKZ94,Ewa96}, while some details on multiprojective toric varieties over a field are worked out in~\cite[\S~1]{PS04}. We simultaneously introduce notations to be used throughout the text, with the proviso that $d=n$ from section~\ref{bk_degree} on. In the following sections  we note by $s$ a point of $S$ while in the introduction and sections~\ref{bk_degree},~\ref{bk_bezout} and \ref{bk_equality}
we use the letter $v$. The reason for this will become apparent only in theorem~\ref{gengenthm}.

\medskip 

\subsection{Torus action and associated maps}\label{torusaction}
For $0\le i\le n$ consider a vector $\cA_i=(a_{i,0},\dots,a_{i,N_i})\in(\Z^d)^{N_i+1}$ of vectors of $\Z^d$ and a vector $\alpha_i=(\alpha_{i,0},\dots,\alpha_{i,N_i})\in(\K(S)^\times)^{N_i+1}$ of non-zero rational functions on $S$, then we set $\widehat{\cA}_i:=(\cA_i,\alpha_i)$. 

We also set $\geom:=(\cA_0,\dots,\cA_n)$, $\boldsymbol{\alpha}:=(\alpha_0,\dots,\alpha_n)$ and $\arith:=(\geom,\bfalpha)$. This latter data defines a map
$$
\varphi_\arith:S\times\T^d\dashrightarrow S\times \P\quad,\quad (s,\bft) \mapsto 
\big(s,\big(\alpha_{i,j}(s) \bft^{a_{i,j}}: 0\le i\le n, 0\le j\le N_i \big) \big)\enspace,
$$
rational in $s$ and monomial in $\bft$. This map extends to a regular one over $S\times \T^d$ because $S$ is smooth. We define $X_\arith \subset S\times \P$, the (multiprojective) {\it toric variety over $S$ associated to  $\arith$}, as the Zariski closure of the image of this map, equipped with the natural projection $\pi:X_{\arith}\twoheadrightarrow S$. 

Similarly, we consider the monomial map $\varphi_\geom:\T^d\rightarrow\P$, $\bft\mapsto(\bft^{a_{i,j}}:i,j)$, and we set  $X_\geom\subset\P$ 
for the standard multiprojective toric variety associated to $\cA$, defined as the Zariski closure of the image of $\varphi_\geom$.

\smallskip
The data $\geom$ induces a diagonal action of the torus $\T^d$ on $S\times \P$ 
$$
*_{\geom}: \T^d\times (S\times \P)\to S\times\P  \quad ,\quad
(\bft,(s,\bfx))\mapsto \bft*_\geom(s,\bfx)
$$
defined by  $\bft*_\geom(s,\bfx) := \big(s,\big(\bft^{a_{i,j}}x_{i,j}: 0\le i\le n, 0\le j\le N_i \big) \big)$. We identify each $\alpha_i$ with the rational map
$$
\alpha_i:S\dashrightarrow \P^{N_i}  \quad ,\quad s\mapsto \big(\alpha_{i,j}(s): 0\le j\le N_i \big) \enspace, 
$$
{\it a priori} only defined on a dense open subset of $S$, but which as before 
extends to a regular one because $S$ is smooth. The image of $\varphi_\arith$ coincides with the orbit under the action $*_{\geom}$ of the graph of the map $\bfalpha=(\alpha_0,\dots,\alpha_n):S\rightarrow\P$. In particular, both $\Im(\varphi_\arith)$ and $X_\arith$ are equivariant. Besides, if $\alpha_{i,j}=1$ for all $i,j$ 
then $X_\arith = S\times X_\geom$.

\medskip 

\subsection{Dimension}\label{dimensions}
The dimension of $X_\arith$ and related varieties can be characterized 
combinatorially as the rank of some $\Z$-modules. For $ 0\le i \le n$ consider the affine span by the vector components of $\cA_i$
$$
L_{\cA_i}:= \sum_{0\le j,k\le N_i}\Z\, (a_{i,j}-a_{i,k}) \subset\Z^d
$$
and let 
$
L_{\geom}:= \sum_{i=0}^n L_{\cA_i} \subset \Z^d$. For each $s\in S$ we also set
$$
L_{i,s}:= \sum_{j,k\ :\ \alpha_i(s)_{j}, \alpha_i(s)_{k}\ne0}\Z\, (a_{i,j}-a_{i,k}) \subset\Z^d 
\enspace, 
$$
the sum being over all $0\le j,k \le N_i$ such that the $j$-th and $k$-th coordinates of the evaluation at $s$ of the {\it regular} map $\alpha_i:S\to \P^{N_i}$ are non-zero, and let $ L_{s}:=\sum_{i=0}^nL_{i,s}$. 

\smallskip
Without loss of generality, we can assume that $L_{\geom}=\Z^d$, modulo a reparameterization of $X_\arith$: if $L_{\geom}\ne\Z^d$ let $\ell : \Z^r \hookrightarrow \Z^d$ be an injective linear map such that $\ell(\Z^r) = L_\geom$ and set $b_{i,j}:= \ell^{-1}(a_{i,j}-a_{i,0}) \in \Z^r$ then 
$\boldsymbol{\cB}:= (b_{i,j}: 0\le i\le n, 0\le j\le N_i)$. We then have $L_{\boldsymbol{\cB}}=\Z^r$ and $X_{(\boldsymbol{\cB}, \bfalpha)} = X_\arith$, {\it see}~\cite[chap.~5,~prop.~1.2, p.~167]{GKZ94}. Consequently, from now on we will assume that
$$
L_\geom=\Z^d
$$
unless otherwise explicitely stated. As a byproduct of the above discussion, in case $r=d$ we obtain that $\varphi_\arith$ factorizes through
\begin{equation}\label{diagrammereduc}
\xymatrix{
S\times\T^d \ar[d]_\psi \ar[r]^{\varphi_\arith} &X_\arith\\
S\times\T^d \ar[ru]_{\varphi_{(\cB,\bfalpha)}}  &}
\end{equation}
where $\psi$ is a monomial map corresponding 
to the linear map $\ell$ above,  \'etale 
of degree $|\det(\ell)|=
[\Z^d:L_\geom]$.

\smallskip

Set
$$
\left(X_\arith\right)_s :=\pi^{-1}(s) \subset \{s\}\times\P \quad\mbox{and} \quad 
\varphi_s := \varphi_\arith (s,\cdot) :\T^d\to \{s\}\times \P\enspace.
$$

\begin{prop}\label{oubliee}
With the above notation, we have
$$
\dim(X_\arith) 
= \dim(X_\geom)+1 = d+1 \quad ,\quad 
\dim(\Im(\varphi_s)) = \rank_\Z (L_{s})\enspace.
$$
\end{prop}

\begin{proof}
Note first  that $\varphi_s$ is the monomial map $\bft\mapsto (s,(\alpha_i(s)_{j}\,\bft^{a_{i,j}} : i,j))$ and so the Zariski closure  $\overline{\Im(\varphi_s)}$ 
is a multiprojective toric variety contained in $\left(X_\arith\right)_s$. The equality on the right is the standard formula for the dimension of such a variety. Similarly, it is well known that $\dim(X_\geom)=\rank_\Z(L_\geom)=d$. 

For the equalities on the left let $s$ be a generic point in $S$, we have $L_{s}=L_\geom$. Applying the theorem on dimension of fibers to the projection $\Im(\varphi_\arith)\to S$ we obtain
$$
\dim\left(X_\arith\right) = \dim\left(\Im\left(\varphi_\arith\right)\right) 
= \dim(\Im(\varphi_s))+1 = \rank_\Z (L_{s})+1 = d+1\enspace.
$$ 
\end{proof}

This result shows that $\dim(\Im(\varphi_s)) =  \dim\left(\left(X_\arith\right)_s\right) = d$ for a generic $s$, but it may happen that $\dim(\Im(\varphi_s)) < \dim\left(\left(X_\arith\right)_s\right)$ for a particular $s$.
 
\medskip

We now turn to the projection of $X_\arith$ into $\P$. With notation as in the previous subsection, for $0\le i\le n$ we set
$$L_{\widehat{\cA}_i}:=\sum_{\substack{s\in S\\ 0\le j,k\le N_i}} \Z\, (a_{i,j}-a_{i,k},-\ord_s(\alpha_{i,j}/\alpha_{i,k})) \subset \Z^{d+1}$$
and further $L_{\arith}:=\sum_{i=0}^n L_{\widehat{\cA}_i}$, which is also a
submodule of $\Z^{d+1}$.

\begin{lem}\label{enjolras}
Let $\varpi:S\times\P\to\P$ be the natural projection onto the second factor, then $\displaystyle \dim(\varpi(X_\arith))=\rank_\Z(L_{\arith})$.
\end{lem}

\begin{proof}
We note first that the image $\varphi_\geom(\T^n)$ is a subtorus of the torus 
$$
\P^\circ:=\P\setminus \, \bigcup_{i,j}Z(x_{i,j}) \simeq \T^{\sum_{i=0}^n N_i} 
\enspace,
$$
and setting $S^0:=\{s\in S:\alpha_i(s)_j\not=0 \mbox{ for all } i,j\}$ we have also that $\bfalpha(S^0)\subset\P^\circ$. The dimension of $\varpi(X_\arith)$ is either $d$ or $d+1$, depending on whether the curve $\bfalpha(S^0)$ lies in a translate of $\varphi_\geom(\T^d)$ or not: in case $\bfalpha(S^0)\subset \bfalpha(s_1) \cdot \varphi_\geom(\T^d)$ for some $s_1\in S$ then $\varpi(\varphi_\arith(S^0\times\T^d)) =\bfalpha(s)\cdot \varphi_\geom(\T^d)$ and so $\dim(\varpi(X_\arith))=\dim(X_\geom)=d$. Otherwise, the dimension of $\varpi(X_\arith)$ is strictly bigger and so equal to $d+1$, since this is the biggest it can be. 

The condition $\bfalpha(S^0)\subset \bfalpha(s_1)\cdot \varphi_\geom(\T^d)$ is equivalent to 
$$
\bfalpha(s)^\bfb= \bfalpha(s_1)^\bfb
$$
for all $s\in S^0$ and $\bfb\in\prod_{i=0}^n \Z^{N_i+1}$ such that $\sum_{i,j} b_{i,j}a_{i,j}=0$ and $\sum_{i,j} b_{i,j}=0$~\cite[chap.~VII, thm.~3.12, pp.~278-9]{Ewa96}. This holds if and only if
$$
\sum_{i,j} b_{i,j}\ord_{s}(\alpha_{i,j})=0 \quad ,\quad \mbox{ for all } s\in S \enspace.
$$
This means that the vectors $ a_{i,j}-a_{i,k} $ and $(a_{i,j}-a_{i,k},-\ord_s(\alpha_{i,j}/\alpha_{i,k}))$ satisfy the same set of linear relations, which in this context is equivalent to $\rank_\Z(L_{\arith})=d$. This completes the proof: by proposition~\ref{oubliee} if $\bfalpha(S^0)\subset \bfalpha(s_1)\cdot \varphi_\geom(\T^d)$ then $\dim(\varpi(X_\arith)) = d = \rank_\Z(L_{\arith})$, otherwise $\dim(\varpi(X_\arith)) = d+1 =\rank_\Z(L_{\arith})$. 
\end{proof}

\medskip

\subsection{Finiteness}\label{finiteness}
Our objective in this subsection is to 
determine the subset of $X_\arith$ where $\varphi_\arith$ is finite or an isomorphism. Recall that a map of algebraic varieties $f:X\to Y$  is {\em finite} ({\it resp.} an {\em isomorphism}) at a point $y\in Y$ whenever there are affine open sets $V\subset Y$ with $y\in V$ and $U\subset f^{-1}(V)$ such that $f:U\to V$ is finite ({\it resp.} an isomorphism).

\medskip
Always under the assumption $L_\geom=\Z^d$, consider the dense open subsets of $S$ 
\begin{equation}\label{densesubsetsofS}
\begin{array}{c}
S^\Ef:=\{s\in S: \rank_\Z(L_{s})=d\} \enspace ,\quad 
S^\I:=\{s\in S: L_{s}=\Z^d\}\enspace,\\[2mm]
S^0:= \{s\in S: \alpha_i(s)_{j}\ne 0 \mbox{ for all } i,j \}\enspace,
\end{array}
\end{equation}
and the corresponding subsets of $X_\arith$
$$\begin{array}{c}
X^\Ef_\arith:=\varphi_\arith(S^\Ef\times\T^d) \enspace,\quad 
X^\I_\arith:=\varphi_\arith(S^\I\times\T^d)\enspace,\quad
X^0_\arith:=\varphi_\arith(S^0\times \T^d)\enspace.
\end{array}$$
Notice the chains of inclusions $S^0\subset S^\I\subset S^\Ef\subset S$ and $X^0_\arith\subset X^\I_\arith\subset X^\Ef_\arith\subset\Im(\varphi_\arith)$.

\begin{lem}\label{isomorphismegenerique}
The map $\varphi_\arith:S\times \T^d \to X_\arith$ is finite ({\it resp.} an isomorphism) 
at a point $(s,\bfx)\in X_\arith$ if and only if 
$(s,\bfx)\in X^\Ef_\arith$ ({\it resp.} $(s,\bfx)\in X^\I_\arith$). 
\end{lem}

\begin{proof}
The ``only if'' direction is easy: the map $\varphi_\arith$ is finite ({\it resp.} an isomorphism) at $(s,\bfx)$ only if the fiber $\varphi_\arith^{-1}(s,\bfx)$ is 
non-empty and finite ({\it resp.} consists of only one point) which is equivalent to $(s,\bfx)\in X_\arith^\Ef$ ({\it resp.} $(s,\bfx)\in X_\arith^\I$). 

Conversely, let $(s,\bfx)\in X_\arith^\Ef$ so that $\rank_\Z(L_s)=d$, and take a basis $h_1,\dots,h_d$ of $\Z^d$ such that $\delta_1h_1,\dots,\delta_dh_d$ is a basis of $L_{s}$ for some $\delta_i\ge 1$ (invariant factors) such that $\delta_1|\delta_2|\cdots|\delta_d$. Take then $\lambda_{i,j,k}^\ell \in \Z$ for $1\le \ell\le d$, $0\le i\le n$ and $0\le j,k \le N_i$ such that $\alpha_i(s)_{j}, \alpha_i(s)_{k}\ne 0$, satisfying
$$
\sum_{i,j,k}\lambda_{i,j,k}^\ell (a_{i,j}-a_{i,k})=\delta_\ell h_\ell \in L_{s}\enspace.$$
Consider the map, well-defined in a neighborhood of $\{s\}\times\T^d$, 
$$X_\arith^\Ef\dashrightarrow S^\Ef\times \T^d \quad,\quad (s,\bfx)\mapsto \Big(s, \Big( 
\prod_{i,j,k} \Big(\frac{x_{i,j} \alpha_{i}(s)_k}{x_{i,k} \alpha_{i}(s)_j}\Big)^{\lambda_{i,j,k}^\ell} : 1\le \ell\le n\Big)\Big)\enspace. 
$$
The composition $S^\Ef\times\T^d\, \displaystyle\mathop{\longrightarrow}^{\varphi_\arith}\, X_\arith^\Ef\longrightarrow S^\Ef\times\T^d$ is  the finite map 
$$
(s,\bft)\mapsto(s,(\bft^{\delta_\ell h_\ell}: 1\le \ell\le d))\enspace, 
$$
and so $\varphi_\arith$ is finite at $(s,\bfx)$. In case $(s,\bfx)\in X_\arith^\I$, we take the $h_i$'s as the standard basis of $\Z^d$ so that $\delta_i=1$ for all $i$ and the composite map is an isomorphism, hence so is $\varphi_\arith$. 
\end{proof}

This result implies that $X_\arith$ is birational to $S\times \T^d$. The above proof gives the additional information that for $(s,\bfx)\in \Im(\varphi_\arith)\setminus X^\Ef_\arith$ the inverse image $\varphi_\arith^{-1}(s,\bfx)$ contains no isolated point, while of course for $(s,\bfx)\in X_\arith\setminus \Im(\varphi_\arith)$ it is empty. 

\medskip

\subsection{Parameterizations and initial coefficients}\label{parameterizations}
For $s\in S$, a {\it parameterization of $S$ at $s$} is defined as a local analytic isomorphism $g_s:\K\to S$ such that $g_s(0)=s$. For such a $g_s$ and a rational function $\beta\in \K(S)^\times$ we have 
$$
\beta\circ g_s(z) = \lambda z^{\ord_s(\beta)} + \mbox{ higher order terms } 
\quad (\mbox{\rm h.o.t.})
$$ 
for some $\lambda\in \K^\times$. We then set $\lambda_{g_s}(\beta) :=\lambda$ for the {\it initial coefficient of $\beta$ at $s$} relative to $g_s$. Given a second  parameterization $h_s$ at the same point $s$ we have $g_s^{-1}\circ h_s(z) = \nu z + \hot$ for some $\nu\in\K^\times $ and so 
\begin{align*}
\beta\circ h_s(z) &= (\beta\circ g_s)\circ(g_s^{-1}\circ h_s)(z)\\[1mm]
&= \lambda_{g_s}(\beta) (g_s^{-1}\circ h_s(z))^{\ord_s(\beta)} + \hot\\[1mm]
&= \lambda_{g_s}(\beta) \nu^{\ord_s(\beta)} z^{\ord_s(\beta)} + \hot \enspace, 
\end{align*}
hence $\lambda_{h_s}(\beta)=\nu^{\ord_s(\beta)}
\lambda_{g_s}(\beta) 
$. 
That is, for $\bfalpha\in\prod_{i=0}^n(\K(S)^\times)^{N_i+1}$, changing  the local parameterization $g_s$ acts on the vector
$$
\lambda_{g_s}(\bfalpha):=\big(\lambda_{g_s}(\alpha_{i,j}):0\le i\le n,0\le j\le N_i\big) \in\prod_{i=0}^n(\K^\times)^{N_i+1}
$$ 
as the 1-dimensional torus action associated to $(\ord_s(\alpha_{i,j}): i,j) \in \prod_{i=0}^n\Z^{N_i+1}$. We denote $\lambda_s(\bfalpha)$ the point $\lambda_{g_s}(\bfalpha)$ modulo this action.

\bigskip

\section{Orbit decomposition}
\label{orbitdecomposition}

The variety $X_\arith$  decomposes as the union of its fibers over $S$: with the notation in subsections~\ref{torusaction} and~\ref{dimensions} we have 
$$
X_{\arith}= \bigsqcup_{s\in S}\left(X_{\arith}\right)_s\enspace.
$$
Our next objective is to study the geometry of these fibers, which turns to be governed by the upper envelope of the corresponding family of $s$-adic polytopes. More precisely, for $s\in S$ and $0\le i\le n$ consider the {\it $s$-adic polytope associated to $(\cA_i,\alpha_i)$}
$$
Q_{i,s}:=\Conv\big((a_{i,j},-\ord_s(\alpha_{i,j})) : 0\le j\le N_i\big) \subset \R^{d+1}
\enspace,$$
which sits above the polytope  $Q_i:=\Conv(a_{i,0},\dots,a_{i,N_i})\subset\R^n$ {\it via} the natural projection $\R^{d+1}\to \R^d$; then set $\boldsymbol{Q}_s:=(Q_{0,s},\dots,Q_{n,s})$  and $\boldsymbol{Q}:=(Q_{0},\dots,Q_{n})$ for the families of those polytopes. 
For $\tau\in \R^{d}$ we define the {\it slope   of $Q_{i,s}$ in the direction $(\tau,1)$} as
$$
Q_{i,s}^{(\tau,1)}:= \{u \in Q_{i,s} : \langle u,(\tau,1)\rangle \ge \langle w,
(\tau,1)\rangle \mbox{ for all } w \in Q_{i,s}\} 
\enspace, $$
which is a face of the upper envelope of $Q_{i,s}$. Put further
$$ 
\bfQ_s^{(\tau,1)} :=\big(Q_{0,s}^{(\tau,1)},\dots, Q_{n,s}^{(\tau,1)} \big)
$$
and let $\Slope(\bfQ_s)$ denote the set of  families of slopes  $\bfQ_s^{(\tau,1)}$ obtained as $\tau$  varies. 

\medskip

Fix both $s\in S$ and a parameterization $g_s$ of $S$ at $s$. For each family of faces $\bfF=(F_0,\dots,F_n)\in\Slope(\bfQ_s)$ consider the point $\bfx_{g_s,\bfF} \in \P$ defined by $(\bfx_{g_s,\bfF})_{i,j}= \lambda_{g_s}(\alpha_{i,j})$ if $(a_{i,j},-\ord_s(\alpha_{i,j}))\in F_i$ and $(\bfx_{g_s,\bfF})_{i,j}= 0$ if not.

\begin{lem}\label{welldefined}
The orbit $\T^d*_\geom (s, \bfx_{g_s,\bfF})$ does not depend on the choice of $g_s$. 
\end{lem}

\begin{proof}
Take $\tau$  such that $\bfF=\bfQ_s^{(\tau,1)}$, then for 
all $i,j$ such that $(a_{i,j},-\ord_s(\alpha_{i,j}))\in F_i$ 
$$
\langle (a_{i,j},-\ord_s(\alpha_{i,j})), (\tau,1)\rangle = c_i
$$
for some $c_i\in \Z$ independent of $j$, or equivalently $\ord_s(\alpha_{i,j})= \langle a_{i,j},\tau\rangle - c_i $. 
Now let $h_s$ be a second parameterization of $S$ at $s$. By the results in subsection~\ref{parameterizations}, there exists $u\in \K$ such that $ \lambda_{h_s}(\alpha_{i,j}) = \lambda_{g_s}(\alpha_{i,j}) u^{\ord_s(\alpha_{i,j})}$. 
Together with the above, this implies $(s,\bfx_{h_s,\bfF}) = (u^{\tau_1},\dots, u^{\tau_n}) *_\geom (s,\bfx_{g_s,\bfF})$ and so $\T^d*_\geom (s, \bfx_{h_s,\bfF})=\T^d*_\geom (s, \bfx_{g_s,\bfF})$, which concludes the proof. 
\end{proof}

Consequently we set 
\begin{equation}\label{defXSF}
X_{s,\bfF}:= \T^d*_\geom (s, \bfx_{g_s,\bfF}) \subset \{s\}\times\P
\end{equation}
for the orbit of this point under the action  $*_\geom$. Note that this orbit is a translate of a torus embedded in some coordinate subspace of $\{s\}\times\P$. 

\medskip
The following proposition gives the orbit decomposition of $\left(X_\arith\right)_s$ and {\it a fortiori} that of $X_\arith$. The unmixed case ($n=0$) was established by A.L.~Smirnov: localizing at $s$, $X_\arith$ can be viewed as a toric variety over the discrete valuation ring $\cO_{S,s}$ and $\left(X_\arith\right)_s$ as its special fiber; the result can be found in these terms in~\cite[\S~2.4]{Smi97}. 

\begin{prop}\label{X_v}
With the notation introduced, we have $\displaystyle \left(X_\arith\right)_s=\bigsqcup_{\bfF\in \Slope(\bfQ_s)}X_{s,\bfF}$.
\end{prop}

\begin{proof}
For $\xi\in \left(X_\arith\right)_s$ let $\cC_\xi\subset X_\arith$ be the germ of an analytic curve containing $\xi$ such that $\cC_\xi \cap X^\I_\arith \ne \emptyset$ and $\cC_\xi\not\subset \left(X_\arith\right)_s$. Let $\eta:\K\to \cC_\xi$ be a parameterization of $\cC_\xi$ at $\xi$. By lemma~\ref{isomorphismegenerique}, the restriction of $\varphi_\arith$ to $S^\I\times\T^d$ is an isomorphism, and so restricting to the Zariski dense subset $U:=\eta^{-1}(\cC_\xi\cap X^\I_\arith)\subset \K$ we can factorize $\eta$ through $S^\I\times\T^d$: 
\begin{displaymath}
\xymatrix{
U \ar[r]^\eta \ar[dr]_\theta& X^\I_\arith\\
 &S^\I\times \T^d \ar[u]_{\varphi_\arith} 
}
\end{displaymath}
where $\theta(z) =(C(z), \bfD(z))$ for some $C:U\to S$ and $\bfD=(D_1, \dots,D_d):U\to \T^d$. The hypothesis $\cC_\xi\not\subset X_s$ implies that $C(z)$ is not constant; besides this analytic function extends to a regular one in a neighborhood of $s$, such that $C(0)=s$. 

\smallskip
Let $g_s: \K\to S$ be a parameterization of $S$ at $s$. The previous considerations imply 
$$
g_s^{-1}\circ C(z)= \gamma z^c + \hot \quad \mbox{ for some }\gamma\in\K^\times \mbox{ and } c\in\N^\times\enspace.
$$
Modifying the parameterization $g_s$ if necessary, we can even assume $\gamma=1$. We also have for $1\le i\le d$ 
$$
D_i(z)= \delta_i z^{-c\tau_i}+ \hot \enspace, 
$$ 
for some $\tau_i\in \R$ and $\delta_i\in \K^\times$. 
Thus
$$
\alpha_{i,j} (C(z))= (\alpha_{i,j} \circ g_s) (g_s^{-1}\circ C(z)) = \lambda_{g_s}(\alpha_{i,j}) z^{c\cdot \ord_s(\alpha_{i,j})} + \hot
$$
and, putting $\tau:=(\tau_1,\dots,\tau_d)$, 
$\bfD:=(D_1,\dots, D_d)$ and  
$\bfdelta:=(\delta_1,\dots,\delta_d)$,   
\begin{align*}
\eta(z) &= \varphi_\arith\circ\theta(z) \\[0mm]
&= \Big(C(z), \Big( \alpha_{i,j}(C(z)) \bfD(z)^{a_{i,j}} : i,j\Big)\Big) \\[1mm]
&=\Big(C(z), \Big( \lambda_{g_s}(\alpha_{i,j}) \bfdelta^{a_{i,j}} z^{\langle -c(\tau,1),(a_{i,j},-\ord_s(\alpha_{i,j}))\rangle} + \hot : i,j \Big)\Big)\enspace.
\end{align*}
For $z$ going to $0$ and $i$ fixed, only survive the initial parts of the 
$j$-th coordinates such that the scalar product $\langle (\tau,1),(a_{i,j},-\ord_s(\alpha_{i,j}))\rangle$ is maximal for $0\le j\le N_i$. Therefore, $\xi= \eta(0)$ satisfies $\xi_{i,j}= \lambda_{g_s}(\alpha_{i,j}) \bfdelta^{a_{i,j}}$ if $(a_{i,j},-\ord_s(\alpha_{i,j}))\in Q_{i,s}^{(\tau,1)}$ and $\xi_{i,j}=0$ if not, which 
implies that $\xi = \bfdelta *_\geom (s,\bfx_{s,\bfQ_s^{(\tau,1)}})$ belongs to the orbit $X_{s,\bfF}$ for $\bfF=\bfQ_s^{(\tau,1)}$.

\smallskip

Conversely, let $\bfF\in\Slope(\bfQ_s)$ and take any $\xi\in X_{s,\bfF}$. By definition, there is some $\bft\in\T^d$ such that $\xi_{i,j}=\lambda_{g_s}(\alpha_{i,j}) \bft^{a_{i,j}}$ if $(a_{i,j},-\ord_s(\alpha_{i,j}))\in F_i$ and  $\xi_{i,j}=0$ if not. With the above notations, take  $\tau=(\tau_1,\dots,\tau_d)$ such that $\bfF=\bfQ_s^{(\tau,1)}$ and consider an analytic function $\eta:\K\to X$ of the form
$$
\eta(z)=\varphi_\arith (g_s(z), z^{-\tau_1}t_1,\dots,z^{-\tau_d}t_d) \enspace.
$$
We readily verify that $\eta(0)=\lim_{z\to 0} \eta(z)=\xi$ which implies $\xi\in \left(X_\arith\right)_s$, as desired. 
\end{proof}

This shows that our toric variety over $S$ decomposes as the infinite union of orbits 
\begin{equation} \label{orbitdecompositionformula}
X_\arith=\bigsqcup_{s\in S}\ \bigsqcup_{\bfF\in \Slope(\bfQ_s)}X_{s,\bfF} \enspace.
\end{equation}
In the sequel we show that the orbits over the dense open subset $S^0= \{s\in S: \bfalpha_i(s)_{j}\ne 0 \mbox{ for all } i,j \}$ of $S$ can be glued together into a finite number of families, {\it see} identity~(\ref{X^0decomposition}) below. This remark will be of importance for the analysis of the equality conditions at the end of section~\ref{bk_equality}. 

For $0\le i \le n$ and $\sigma\in \R^d$ let $Q_i^\sigma \subset \R^d$ be the face 
made of the points in $Q_i$ maximizing the functional $u\mapsto\langle \sigma, u\rangle$ over the polytope, then set
$$
\Faces(\bfQ) :=\{\bfQ^\sigma=(Q_{0}^\sigma, \dots,Q_{n}^\sigma): \sigma\in \R^d\} 
$$
for the families of faces obtained in this way. For $s\in S^0$, the set $\Slope(\bfQ_s)$ is in bijection with the set $\Faces(\bfQ)$: indeed, any element in $\Slope(\bfQ_s)$ is of the form $(F_0'\times\{0\},\dots,F_n'\times\{0\})$ for some $\bfF'=(F'_0,\dots, F_n') \in \Faces(\bfQ)$. In particular, $\Slope(\bfQ_s)$ is independent of $s\in S^0$.  

For such a  $\bfF'\in \Faces(\bfQ)$ we set $X_{S^0,\bfF'}$ for the image of the map $S^0\times \T^d \to S^0\times \P, (s,\bft)\mapsto (s,\bfx)$ where $x_{i,j}= \alpha_{i,j}(s)\bft^{a_{i,j}}$ if $a_{i,j} \in F_i$ and $x_{i,j}=0$ otherwise. We have $X_{S^0,\bfF'} = \bigsqcup_{s\in S^0}X_{s,\bfF}$ with $\bfF:=(F_0'\times\{0\},\dots,F_n'\times\{0\})$ and so it follows from proposition~\ref{X_v}
\begin{equation}\label{X^0decomposition}
\bigsqcup_{s\in S^0}\left(X_\arith\right)_s=\bigsqcup_{\bfF'\in \Faces(\bfQ)} X_{S^0,\bfF'} \enspace. 
\end{equation}

Besides, note that for $s\in S^0$ each fiber $\left(X_\arith\right)_s$ is linearly isomorphic to $X_\geom$ and thus we recover the orbit decomposition of a multiprojective toric variety over $\K$, {\it see}~\cite[\S~1]{PS04}. 

\medskip

With notation as in subsection~\ref{finiteness}, consider the chain of inclusions
$$
X^0_\arith\subset X^\I_\arith\subset X^\Ef_\arith \subset \Im(\varphi_\arith) 
\enspace.$$
These are all 
equivariant subsets of $X_\arith$ and, in particular, an orbit is contained in one of these subsets if and only if
it contains a point in it. 
The following proposition shows that for each $s\in S$, these subsets contain at most one specific orbit.

\begin{prop}\label{combeferre} 
Let $s\in S$ and $\bfF\in \Slope(\bfQ_s)$, then 
\begin{enumerate}
\item\label{Im(varphi)} $X_{s,\bfF}\subset \Im(\varphi_\arith)$ if and only if $\bfF=\bfQ_s^{(\bfzero,1)}$;\smallskip
\item\label{X^Ef} $X_{s,\bfF}\subset X^\Ef_\arith$ if and only if $\bfF=\bfQ_s^{(\bfzero,1)}$ and $s\in S^\Ef$;\smallskip
\item\label{X^I} $X_{s,\bfF}\subset X^\I_\arith$ if and only if $\bfF=\bfQ_s^{(\bfzero,1)}$ and $s\in S^\I$;\smallskip
\item\label{X^0} $X_{s,\bfF}\subset X^0_\arith$ if and only if $\bfF=\bfQ_s^{(\bfzero,1)}$ and $s\in S^0$.
\end{enumerate}

In particular, $X_\arith^0 = \bigsqcup_{s\in S^0}X_{s,\bfQ_s^{(\bfzero,1)}}$, $X_\arith^\I = \bigsqcup_{s\in S^\I}X_{s,\bfQ_s^{(\bfzero,1)}}$, $X_\arith^\Ef = \bigsqcup_{s\in S^\Ef}X_{s,\bfQ_s^{(\bfzero,1)}}$ and $\Im(\varphi_\arith) = \bigsqcup_{s\in S}X_{s,\bfQ_s^{(\bfzero,1)}}$.
\end{prop}

\begin{proof}
The restriction of the image of  $\varphi_\arith$ to the fiber $\left(X_\arith\right)_s$ coincides with the image of the monomial map $\bft\mapsto (s,(\alpha_i(s)_{j}\bft^{a_{i,j}}: i,j))$ and hence with the orbit of the point $(s,\bfalpha(s))$. For any parameterization $g_s$ of $S$ at $s$, we have $\alpha_i(s)_{j}=\lambda_{g_s}(\alpha_{i,j}) $ for $i,j$ such that $-\ord_s(\alpha_{i,j})$ is maximal for $0\le j\le N_i$  and $\alpha_i(s)_{j}=0$ otherwise, and so $\bfx_{s,\bfF}=\bfalpha(s)$ for $\bfF=\bfQ_s^{(\bfzero,1)}$, which proves (1). The statements concerning the subsets $X_\arith^\Ef$, $X_\arith^\I$ and $X_\arith^0$ are direct consequences of the previous, together with their definition. And, in view of (1)-(4), the orbit decompositions for these subsets and $\Im(\varphi_\arith)$ result directly from the orbit decomposition for $X_\arith$ in~(\ref{orbitdecompositionformula}).
\end{proof}

\bigskip

\section{Mixed degrees and mixed integrals}\label{bk_degree}

In this section we obtain a combinatorial formula for a certain mixed degree of a toric variety over a curve (proposition~\ref{fauchelevent}). This is the function field analog of the formula for the normalized height of a toric variety defined over $\overline{\Q}$ in~\cite{PS04} and it constitutes the core of the proof of our main results. We also obtain a characterization for the vanishing of this mixed degree in terms of the rank of some $\Z$-modules (proposition~\ref{anzelma}).

\bigskip
Let $\cS$ be an arbitrary variety and consider again the multiprojective space $\P =\P^{N_0}\times \dots\times\P^{N_n}$. Let $\pi_i:\cS\times \P \to \P^{N_i}$  be the natural projection to the $i$th factor. Given an equidimensional cycle $Z$ of $\cS\times \P$ of dimension $d$ and a  multi-index $c\in\N^{n+1}$ of length $c_0+\dots+c_n=d$,  the corresponding {\it mixed degree} (or {\it multi-degree}) is
$$\deg_{c}(Z) := \deg\big(Z\cdot \pi_0^{-1}(E_0)\cdot \, \cdots \, \cdot \pi_{n}^{-1}(E_{n})\big) \ge 0
\enspace,$$
where $\cdot$ denotes the intersection product and $E_i\subset \P^{N_i}$ a generic linear subvariety of codimension $c_i$. 

\smallskip
It is  useful to know that mixed degrees can be interpreted in terms of resultants,
whenever $\cS$ is projective and $Z$ is given as a subscheme of $\cS\times\P$. In the sequel we explain this in the case which concerns us: $\cS=S$ a complete smooth curve, 
$Z\subset S\times \P$ a variety of dimension $n+1$ and the multi-index $c=(1,\dots,1)\in \N^{n+1}$.

Choose an embedding $S\hookrightarrow \P^{N_{-1}}$ so that $Z$ becomes a subvariety 
of $\P^{N_{-1}}\times \P$. For $-1\le i\le n$ we introduce a group of $N_i+1$ variables $\bfU_i=\{U_{i,0},\dots, U_{i,N_i}\}$ corresponding to the coefficients of the general linear form $L_i=\sum_{j=0}^{N_i} U_{i,j}x_{i,j}$ in the variables $\bfx_i= \{ x_{i , 0}, \dots, x_{i, N_i}\} $. We then consider the corresponding resultant of $Z$
$$
\Res_{Z} \in \K[\bfU_{-1},\bfU_0, \dots, \bfU_n] 
$$ 
as defined and studied in~\cite[chap.~5]{Rem01b}, {\it see} also~\cite[\S~I.2]{PS04}. In the terminology of these references, this is the resultant of the multihomogeneous ideal $I(Z)\subset \K[\bfx_{-1},\bfx_0, \dots, \bfx_n]$ with 
respect to the vectors  $e_{-1},e_0, \dots, e_n \in \Z^{n+2}$ of the standard basis of $\R^{n+2}$. We prompt the reader to the above references for the exact definition and fundamental properties of the resultant, we only note here that, whenever $\Res_Z\ne1$, then $\Res_Z(\bfu_{-1},\dots,\bfu_n)=0$ for given  $\bfu_i\in \K^{N_i+1}$ if and only if $Z\cap Z(L_{-1}(\bfu_{-1}),\dots, L_n(\bfu_n))\ne \emptyset$. 

The relevant mixed degree of $Z$ is given by the degree of $\Res_Z$ in the $\bfU_{-1}$-variables:
\begin{equation}\label{mixeddegree}
\deg_{(1,\dots,1)}(Z) = \deg_{\bfU_{-1}}\left(\Res_{Z}\right)\enspace,
\end{equation}
and this equality does not depend on the choice of the projective embedding of $S$ \cite[chap.~5, prop.~3.4 and 2.11]{Rem01b}. 

\medskip

From now on, we set $d=n$ and we reconsider the data introduced in \S~\ref{torusaction}:
\begin{equation}\label{datatorique}
\geom \in\prod_{i=0}^n (\Z^n)^{N_i+1}
\quad, \quad 
\bfalpha\in \prod_{i=0}^n (\K(S)^\times)^{N_i+1} \quad ,\quad 
\arith=(\geom,\bfalpha)\enspace.
\end{equation}
As before, we will assume that $L_\geom=\Z^n$, which in particular implies that $X_\arith $ has dimension $n+1$. 
The following is the main result of this section. 

\begin{prop}\label{fauchelevent}
For $v\in S$  and $0\le i\le n$ 
let  $\vartheta_{i,v}$ denote the roof function of the 
$v$-adic polytope $Q_{i,v}=\Conv((a_{i,j},-\ord_v(\alpha_{i,j}): j) $
above $Q_i=\Conv(a_{i,j}:j)$, then 
\begin{equation}\label{fauchelevent_display}
\deg_{(1,\dots,1)}(X_\arith) = \sum_{v\in S} \MI_n(\vartheta_{0,v},\dots,\vartheta_{n,v})
\enspace.\end{equation}
\end{prop}

The proof is done in two steps: first we compute this mixed degree for a $1$-dimensional deformation of an arbitrary variety in terms of mixed Chow weights, then we show that in the toric case, mixed Chow weights can be expressed  in terms of mixed integrals. 

\smallskip
\begin{defn}
Let $W\subset \P$ be a variety and $\bfalpha$ as in~(\ref{datatorique}). The {\em $\bfalpha$-deformation} $W_\bfalpha \subset S\times\P$ is defined as the Zariski closure of the set of points of the form
$$
\left\{\left({s},( \alpha_{i,j}(s) x_{i,j}: i,j)\right)   : s\in S^0, \ \bfx\in W \right\} 
\enspace.
$$
\end{defn}

Note that $W_\bfalpha$ is of dimension $\dim(W)+1$. Also note that  $X_\arith$ is the $\bfalpha$-deformation of the toric variety $X_\geom\subset \P$, so that with the above notation $X_\arith=(X_\geom)_{\bfalpha}$.

\medskip
Let $W\subset \P$ be a $n$-dimensional variety and set $\Res_{W}\in \K[\bfU_0,\dots,\bfU_n]$ for its resultant with respect to the standard basis of $\R^{n+1}$. Given a vector $\bftau \in \prod_{i=0}^n\R^{N_i+1}$, the {\it mixed Chow weight} (or {\it Chow multi-weight}) $e_{\bftau}(W)$ is defined as the weight in the $\bftau$-direction of the Newton polytope of $\Res_{W}$. Introducing an additional variable $T$,
$$
e_{\bftau}(W)= \deg_T\big(\Res_{W} (T^{\tau_{i,j}}\, U_{i,j}: 0\le i\le n, 0\le j\le N_i)\big)
\enspace.
$$
Chow weights of projective varieties were introduced and studied by D.~Mumford 
in the studying stability in geometric invariant theory~\cite{Mum77}. The extension to the multiprojective setting was done in our previous text~\cite[\S~IV.2]{PS03}. 

\begin{prop} \label{madeleine}
Let $W\subset \P$ be a variety of dimension $n$. For $v\in S$ set $\tau_{\bfalpha,v}:= (-\ord_v(\alpha_{i,j}): i,j) \in \prod_{i=0}^n\Z^{N_i+1}$, then 
\begin{equation} \label{madeleine_display}
\deg_{(1,\dots,1)}(W_\bfalpha) = \sum_{v\in S} e_{\tau_{\bfalpha,v}}(W)\enspace.
\end{equation}
\end{prop}

Choose  $S\hookrightarrow \P^{M}$ an embedding such that the linear projection $\ell:\P^{M}\to \P^1$, $\bfy\mapsto (y_0:y_1)$  induces  a finite map from  $S$ onto $\P^1$.
In particular $\K(\P^1)\hookrightarrow \K(S)$ is a finite extension. Besides, $W_\bfalpha$ becomes a subvariety of 
$\P^{M}\times \P$ and
introducing a group of variables $\bfV=\{V_0,\dots,V_M\}$ 
 we can consider the resultants
$$
\Res_{W_\bfalpha }\in 
\K[\bfV, 
\bfU_0, \dots, \bfU_n] \quad ,\quad 
\Res_{W}\in \K[ \bfU_0, \dots, \bfU_n]\enspace.
$$
The proof of the proposition above depends on the following Poisson-type formula.

\begin{lem}\label{poisson}
With notation as above,  for $r\in \K$ set $u(r):=(-r, 1, 0,\dots,0)$, then there 
exists a rational function $q\in \K(\P^1)^\times$ such that
$$
\Res_{W_\bfalpha}\Big(u\Big(\frac{\ell(s)_1}{\ell(s)_0}\Big)
 , \bfU\Big) = q(\ell(s)) 
\prod_{\sigma}\Res_{W}(\sigma(\alpha_{i,j})\big({s})U_{i,j}: i,j\big)
$$
for generic $s\in S$, 
where $\sigma :\K(S)\hookrightarrow \ov{\K(\P^1)}$ runs over all $\K(\P^1)$-embeddings
of $\K(S)$.
\end{lem}

\begin{proof}
We assume throughout the proof that  $s\in S$ is generic. Set $r:=\frac{\ell(s)_1}{\ell(s)_0}\in \K$ and $W_\bfalpha\cdot H_{r}$ for the intersection cycle of $W_\bfalpha$ with the hypersurface  $H_r:=Z( ry_0-y_1)  \subset \P^{M}\times \P$ defined by $u(r)$. By~\cite[chap.~5, prop.~3.6]{Rem01b} we have 
$$
\Res_{W_\bfalpha}(u(r),  \bfU) = q(r) \, \Res_{W_\bfalpha\cdot H_{r}; e_0,\dots, e_n }(\bfU)
$$
for some rational function $q\in \K(\P^1)^\times$, where $\Res_{W_\bfalpha\cdot H_{r}; e_0,\dots, e_n }$ denotes the resultant with respect to the $n+1$ last vectors in  the standard basis of $\R^{n+2}$, {\it see}~\cite[chap.~5, \S~3]{Rem01b} for precisions.

Since $s$ is generic, Bertini's theorem implies that $W_\bfalpha\cdot H_r$ is reduced, and so it coincides with the set theoretical intersection $ W_\bfalpha\cap H_r$. Besides,  $\ell^{-1}(r) \subset S^0$ and so $W_\bfalpha\cap H_r = \cup_{p\in \ell^{-1}(r)} \{p\}\times (\bfalpha(p) W)$, where $\bfalpha(p)W:= \{(\alpha_{i,j}(p) x_{i,j}: i,j) : \bfx\in W\} \subset \P$ and so 
\begin{align*}
\Res_{W_\bfalpha\cdot H_{r}; e_0,\dots, e_n } (\bfU)& = \prod_{p\in\ell^{-1}(r)}
\Res_{\{p\}\times \bfalpha(p) W; e_0,\dots, e_n } \big(U_{i,j}:  i,j\big)\\[0mm]
&=  \prod_{p\in\ell^{-1}(r)}\Res_{\{p\}\times W;  e_0,\dots, e_n }
\big(\alpha_{i,j}(p)U_{i,j}:  i,j)\\[0mm]
&=  \prod_{p\in\ell^{-1}(r)}\Res_{W} \big(\alpha_{i,j}(p)U_{i,j}:  i,j) 
\enspace, 
\end{align*}
the last equality comes from the fact that $\Res_{\{p\}\times W;  e_0,\dots, e_n } = \Res_{W }$, since the quotient rings $\K[\bfx]/I(W)$ and   $\K[S][\bfx]/I(\{p\}\times W)$ coincide. We finally observe that the $\K(\P^1)$-embeddings of $\K(S)$ into  $\ov{\K(\P^1)}$ act on the fiber $\ell^{-1}(r)$ (which contains $s$) by permutation of the points, this concludes the proof. 
\end{proof}

\begin{Proof}{Proof of proposition~\ref{madeleine}.}
Let notation be as in lemma~\ref{poisson}. By the very definition of the mixed Chow weight we have
\begin{align}\label{dspreuvepropmadeleine}\nonumber
e_{\tau_{\bfalpha,v}}(W) 
&= \deg_T\big( \Res_{W}(T^{-\ord_v(\alpha_{i,j})}\, U_{i,j}: i,j)\big)\\[0mm]\nonumber
&= \max_\bfb\Big(-\sum_{i=0}^n\sum_{j=0}^{N_i}b_{i,j}\ord_v(\alpha_{i,j}) \Big)\\[1mm]
&= -\ord_v\big(\Res_{W}(\alpha_{i,j}({s})U_{i,j}: i,j)\big)\enspace, 
\end{align}
where the maximum runs over the set of exponents ${\bfb}=(b_{i,j}: 0\le i\le n, 0\le j\le N_i)$ of the monomials occuring in $R$.

On the other hand,  $\Res_{W_\bfalpha }(u(r), \bfU) \in \K[r,\bfU]$ has no non-trivial factor in $\K[r]$, since otherwise this would imply that the projection of  $W_\bfalpha$ through $S\times \P\to \P^1$ is a point. Hence, this is a primitive polynomial of degree, with respect to the variable $r$, equal to the degree of $\Res_{W_\bfalpha }$, with respect to the group $\bfV$, and so 
\begin{equation*}\label{dspreuvepropmadeleinei}
\deg_{\bfV}(\Res_{W_\bfalpha }) 
= \deg_r\big(\Res_{W_\bfalpha }(u(r), \bfU)\big)
= -\sum_{w\in\P^1}\ord_w\big(
\Res_{W_\bfalpha }(u(r), \bfU) \big)
\enspace.\end{equation*}
Identity~(\ref{mixeddegree}) and lemma~\ref{poisson} together with the product formula $(\sum_{w\in\P^1}\ord_w(q)=0)$ then imply
\begin{align*}
\deg_{(1,\dots,1)}(W_\bfalpha) &=\deg_{\bfV}(\Res_{W_\bfalpha }) \\[1mm]
&= -\displaystyle\frac{1}{[\k(S):\k(\P^1)]}\sum_{w\in \P^1}\sum_{v\in\ell^{-1}(w)} \ord_v \left( \Res_{W_\bfalpha }\left(u\left(\ell(s)_1/\ell(s)_0\right),\bfU\right)\right) \\[1mm]
&= -\displaystyle\frac{1}{[\k(S):\k(\P^1)]}\sum_{w\in \P^1}\sum_{v\in\ell^{-1}(w)}\sum_\sigma\ord_v \big( \Res_W(\sigma(\alpha_{i,j})({s})U_{i,j}:i,j)\big) \\[1mm]
&= -\displaystyle\sum_{w\in \P^1}\sum_{v\in\ell^{-1}(w)}\ord_v \big( \Res_W(\alpha_{i,j}({s})U_{i,j}:i,j)\big) \\[1mm]
&= \displaystyle\sum_{v\in S}e_{\tau_{\bfalpha,v}}(W)
\enspace,
\end{align*}
as desired. The fourth equality comes from the fact that for each $w\in\P^1$ the sum 
$$\sum_{v\in\ell^{-1}(w)}\ord_v\left(\Res_W(\sigma(\alpha_{i,j})({s})U_{i,j}:i,j)\right)
$$
is independent of $\sigma$ and that the number of $\sigma$'s is equal to $[\k(S):\k(\P^1)]$.
The last equality comes from~(\ref{dspreuvepropmadeleine}).  
\end{Proof}

\medskip
Proposition~\ref{fauchelevent} follows directly from proposition~\ref{madeleine} and the following one. This latter is an extension of~\cite[prop.~IV.6]{PS03}, which supposes $L_{\cA_i}=\Z^n$ for all $i$ instead of our weaker assumption 
$L_\geom=\Z^n$.

\begin{prop}\label{propChowMI}
Let $\bftau\in\prod_{i=0}^n\R^{N_i+1}$ and for $0\le i\le n$ let $\vartheta_i$ 
denote the
roof function of $\Conv((a_{i,j},\tau_{i,j}):i,j)$ above $Q_i=\Conv(a_j:j)$, then 
$$e_{\bftau}(X_\geom)=\MI_n(\vartheta_{0},\dots,\vartheta_{n})\enspace.$$ 
\end{prop}

\begin{proof}
We can reduce without loss of generality to the case when the $\tau_{i,j}$'s are integers, since both sides of the identity are continuous in $\bftau$ and 
homogeneous of degree $1$ with respect to homotheties $\bftau\mapsto \lambda\bftau$ ($\lambda\ge 0$).

For each $0\le i\le n$ let $\mu_i\in \Z$ be such that $\mu_i\le \tau_{i,j}$ for all $j$ and consider
the vector
$$
\cB_i := \big( (a_{i,j},\tau_{i,j}-\mu_i), (a_{i,j},0): 0\le j\le N_i\big) \in (\Z^{n+1})^{2N_i+2} 
$$
and then the polytope $\wt{Q}_i:= \Conv(\cB_i) \subset \R^{n+1}$. From the identity~(\ref{masure_gorbeau})
we obtain
\begin{align} \label{dspreuveChowMIi}
\MI_n(\vartheta_{0},\dots,\vartheta_{n}) =&\MV_{n+1}(\widetilde{Q}_{0},\dots,\widetilde{Q}_{n}) \\[1mm]
&+ \sum_{i=0}^n\mu_{i}\MV_n(Q_0,\dots,Q_{i-1},Q_{i+1},\dots,Q_n) \enspace. \nonumber
\end{align}
We can interpret the mixed volumes in this identity as certain mixed degrees of some toric varieties. 
To this end, set $\wt\P:=\P^{2N_0+1}\times\dots\times\P^{2N_n+1}$ and consider the diagonal embedding 
$$\begin{array}{rccc}
\iota: &\P &\longrightarrow &\wt\P\\
&(x_0,\dots,x_n) &\longmapsto &((x_0:x_0),\dots,(x_n:x_n))
\end{array}$$
For $0\le i\le n$ consider the vector $\beta_i\in(\K(S)^\times)^{2N_i+2}$ defined by 
$$\beta_{i,j}(s):=\left\{ 
\begin{array}{ll}
s^{\tau_{i,j}-\mu_i} &\quad \mbox{ for } 0\le j\le N_i\enspace ,\\[2mm]
1 &\quad \mbox{ for } N_i+1\le j\le 2N_i+1\enspace,  
\end{array}\right.$$
then set $\bfbeta:=(\beta_0,\dots,\beta_n)$. The $\bfbeta$-deformation $\iota(X_\geom)_\bfbeta\subset\P^1\times\wt\P$ can then be identified with the toric variety $X_\bfcalB\subset\P^1\times\wt\P$ associated to the data $\bfcalB:=(\cB_{-1},\cB_0,\dots,\cB_n)$ for $\cB_0,\dots, \cB_n$ as before and  $\cB_{-1}:=\left((\bfzero,0),(\bfzero,1)\right)\in(\Z^{n+1})^2$.

\smallskip
Set $c:=(1,\dots,1) \in \N^{n+1}$ and $c_i=c-e_{i+1} \in \N^{n+1}$, where $e_{i+1}$ denotes the $(i+1)$th vector in the standard basis of $\R^{n+1}$. These are multi-indexes  of length $n+1$ and $n$ respectively, and we have~\cite[prop.~I.2]{PS04} 
\begin{align} \label{degs=MV}
\deg_{c}(\iota(X_\geom)_\bfbeta)&= \MV_{n+1}(\wt{Q}_0,\dots, \wt{Q}_n)
\enspace,\\[2mm]
\deg_{c_i}(X_\geom)&= \MV_n(Q_0,\dots, {Q}_{i-1},{Q}_{i+1},\dots, {Q}_n)
\enspace, \nonumber
\end{align}
because $L_\bfcalB= \Z^{n+1}$ and $L_\geom= \Z^n$.

\smallskip
Proposition~\ref{madeleine} applied to $\iota(X_\geom)_\bfbeta$ finishes the proof: from the construction of $\bfbeta$ we verify that $\tau_{\bfbeta,v} =\bfzero$ for all $v\in \P^1\setminus \{0,\infty\}$ and $\tau_{\bfbeta,\infty} = (\tau_{\beta_0,\infty},\dots,\tau_{\beta_n,\infty}) = -\tau_{\bfbeta,0}$, with 
$$\tau_{\bfbeta_i,\infty}=(\tau_{i,0}-\mu_i,\dots,\tau_{i,N_i}-\mu_i,0,\dots,0)\in\Z^{2N_i+2}
\enspace.
$$
Therefore, the only non zero Chow weight in formula~(\ref{madeleine_display}) corresponds to the place $v=\infty$ and we can write  
\begin{equation}\label{deg_c}
\deg_{c}\big(\iota(X_\geom)_{\bfbeta}\big) = e_{\tau_{\bfbeta,\infty}}(X_\cA) = e_{\bftau}(X_\cA) - \sum_{i=0}^n\mu_{i}\deg_{c_i}(X_\geom)
\enspace,
\end{equation}
by the definition of the mixed Chow weight.
The statement derives now from the identities~(\ref{dspreuveChowMIi}),~(\ref{degs=MV}) and~(\ref{deg_c}) above. 
\end{proof}

Proposition~\ref{fauchelevent} is now fully established. The following proposition characterizes the vanishing of the relevant mixed degree of $X_\arith$. 

\begin{prop}\label{anzelma}
With notation as in proposition~\ref{fauchelevent}, the following conditions are equivalent: 
\begin{enumerate}
\item \label{cond1} $\deg_{(1,\dots,1)}(X_\arith) >0$.
\smallskip
\item \label{cond2} For each $0\le i\le n$ there exists  $j_i\in \{1,\dots, N_i\}$ 
such that the submodule of $\Z^{n+1}$ generated by the vectors 
$$
(a_{i,j_i}-a_{i,0},-\ord_v(\alpha_{i,j_i}/\alpha_{i,0})) \quad \mbox{ for } 0\le i\le n \mbox{ and } v\in S 
$$
is of rank $n+1$.
\smallskip
\item \label{cond3} $\rank_\Z\big(\sum_{i\in I} L_{\wh{\cA}_{i}}\big) \ge  \Card(I)$ for every subset $I \subset \{0,\dots,n\}$.
\end{enumerate} 
\end{prop}

\begin{proof}
(\ref{cond2})$\Rightarrow$(\ref{cond1}): consider the projection  
$$
\rho:S\times\P
\dashrightarrow (\P^1)^{n+1}\quad ,\quad  (s,\bfx)\mapsto \big( (1:x_{0,j_0}), \dots, (1:x_{n,j_n})\big) 
\enspace.
$$
We have $\ov{\rho(X_\arith)}=\varpi(Y)$ where 
$\varpi$ denotes the projection $S\times (\P^1)^{n+1}\to (\P^1)^{n+1}$ 
and $Y$ is the toric variety over $S$ corresponding to the data $(a_{i,0},a_{i,j_i}: 0\le i\le n)$ and $(\alpha_{i,0},\alpha_{i,j_i}: 0\le i\le n) $. By condition~(\ref{cond2}) and lemma~\ref{enjolras} applied to $Y$, it comes that 
$\dim(\varpi(Y))=n+1$ and so $\varpi(Y)=(\P^1)^{n+1}$. 
Since by proposition~\ref{oubliee} we know that $\dim(X_\arith)=n+1$, 
this implies 
that the restriction of $\rho$ to $X_\arith$ is generically finite
and so 
$$
\deg_{(1,\dots,1)}(X_\arith)=\deg(\rho|_{X_\arith})\deg_{(1,\dots,1)}\big((\P^1)^{n+1}\big)
=\deg(\rho|_{X_\arith}) \ge 1 \enspace.
$$

\medskip
(\ref{cond1})$\Rightarrow$(\ref{cond3}): 
suppose that 
$\rank_\Z\big(\sum_{i\in I} L_{\wh\cA_i}\big) < \Card(I)$ for some 
$I \subset \{0,\dots,n\}$. 
Lemma~\ref{enjolras} implies that the projection of $X_\arith$ to $\prod_{i\in I}\P^{N_i}$ is of dimension 
 $ < \Card(I)$ and so 
$$
X_\arith\ \cap \ \bigcap_{i\in I}\pi_i^{-1}(E_i) =\emptyset
$$
for generic hyperplanes $E_i\subset \P^{N_i}$ for $i \in I$, 
which implies that $\deg_{(1,\dots,1)}(X_\arith) =0$.

\medskip
(\ref{cond3})$\Rightarrow$(\ref{cond2}): 
it is known that condition~(\ref{cond3}) in this setting 
implies that there exists a family of linearly independent vectors
$w_i\in L_{\wh{\cA}_i}$ $(0\le i\le n)$~\cite[lem.~5.1.8, p.~278]{Schn93}. 
We can choose the $w_i$'s among the given generators of $L_{\wh{\cA}_i}$, namely $w_i=(a_{i,j_i}-a_{i,0},-\ord_{v_i}(\alpha_{i,j_i}/\alpha_{i,0}))
$ for some $1\le j_i\le N_i$ and $v_i\in S$. These indexes $j_i$ satisfy condition~(\ref{cond2}), this finishes the proof.
\end{proof}

\bigskip

\section{Intersection cycles and the B\'ezout theorem}
\label{bk_bezout}

In this section we recall the necessary background from  multiplicities and multiprojective intersection theory, which is developed in detail in~\cite{Serre75,Bou75,FOV99,Rem01b}. These tools together with the mixed degree computation in proposition~\ref{fauchelevent}, allow to derive the upper bound in our main results and to set the path for the study of the case of equality, to be treated in the next section. 

\bigskip
Let $M$ be a smooth ambient variety, $Z$ a cycle on $M$ and $k\in\N$, we will denote by $|Z|$ the underlying algebraic set and with a subscript the (possibly empty) 
$k$-dimensional part $Z_k$ of $Z$. 

Let now $W,Z$ be pure dimensional cycles on $M$, we define the {\it intersection cycle of $W$ and $Z$} by the formula 
\begin{equation}\label{intersection}
W\cdot Z=\sum_{Y} \mult(Y| W;Z) \, Y
\end{equation}
where $Y$ runs over the irreducible components of 
$|W|\cap |Z|$ of codimension $\codim(W)+\codim(Z)$ 
and $\mult(Y|V;W) \ge 1$ denotes the intersection multiplicity of $W$ and $Z$ along~$Y$. This intersection multiplicity can be alternatively defined as
some Samuel multiplicity or through Serre's Tor-formula, and both definitions coincide~\cite[\S~V.C.1-2]{Serre75}. Note that $W\cdot Z$, as defined in~(\ref{intersection}), never has embedded component, that is a component stricly contained in another component of $W\cdot Z$. However, when the intersection is {\it proper}, namely such that $\codim(|W|\cap |Z|)= \codim(W)+\codim(Z)$, the product $W\cdot Z$ defined in~(\ref{intersection}) coincides with the usual intersection product. In particular, it is always commutative and associative as soon as all the involved intersections are proper~\cite[\S~V.C.3]{Serre75} but not in general. 

\smallskip
The case of interest for us is when $Z$ is defined in $M$ by a complete intersection of regular functions $q_1,\dots, q_r\in \cO(M)$. In this setting, 
we denote $Z(q_1,\dots,q_r)$ the cycle $Z(q_1)\cdot\cdots\cdot Z(q_r)$. 
We also denote $\mult(Y|W;q_1,\dots, q_r)$ the intersection multiplicity of $W$ and $Z$ along a component $Y$ of the proper part of the intersection. 
When $W$ is the ambient space $M$ we simply write $\mult(Y|q_1,\dots,q_r)$ instead of $\mult(Y|q_1,\dots,q_r)$.

For $W$ a variety, this intersection multiplicity is equal to the Samuel multiplicity $e_{\mathfrak q}^{\cO_{W,Y}}(\cO_{W,Y})$ of the local ring $\cO_{W,Y}$ of $W$ along $Y$, relative to the ideal ${\mathfrak q}:=(q_1,\dots,q_r)$. If $\dim(W)=r$ 
and $\xi$ is a point in the proper part of the intersection of $W$ with the zero set of $q_1,\dots, q_r$,  we have 
$$
\mult(\xi|W;q_1,\dots,q_r)\le \dim_\K(\cO_{W,\xi}/{\mathfrak q} )
$$
with equality if and only if $\cO_{W,\xi}$ is a Cohen-Macaulay ring~\cite[cor.~1.2.13, p.~18]{FOV99}. This is the case when $W$ is smooth at $\xi$, as in the setting of theorems~\ref{mainthm} and~\ref{genthm} where $W=M=S\times\T^n$. For instance, for a system of Laurent polynomials $f_0,\dots,f_n\in\k[s][\bft^{\pm 1}]$ and $\xi \in Z(f_0,\dots,f_n)_0$ we have
\begin{equation}\label{mult_intro}
{\rm mult}(\xi|f_0,\dots,f_n)=
\dim_\K\big(\K[s][\bft^{\pm1}]/(f_0,\dots,f_n)\big)_{{\frak m}(\xi)}
\end{equation}
where ${\frak m}(\xi)$ is the ideal of definition of the point $\xi$. 

\medskip
We recall that the degree of a $0$-dimensional cycle is defined as the sum of its multiplicities. The following is a version of B\'ezout theorem adapted to our purposes. 

\begin{lem}\label{Bezoutthm}
Let $W\subset S\times\P$ be a cycle of dimension $n+1$ and for $0\le i\le n$ let
$\ell_i\in  \K[\bfx_i]$ be a linear form, then setting $\bfell=(\ell_0,\dots,\ell_n)$ we have 
$$
\deg\left(W\cdot Z(\bfell)\right) \le \deg_{(1,\dots,1)} (W) \enspace, 
$$
with equality when $|W|\cap|Z(\bfell)|$ is of dimension $0$.
\end{lem}

\begin{proof}
Write  $H_i:=Z(\ell_i)$ for $0\le i\le n$. Set $W^{n+1}:=W$ and for $k=n,\dots, 0$ 
we define inductively cycles $B^{k+1}$, $Z^{k+1}$ and $ W^k$ of pure dimension $k+1$, $k+1$ and $k$ respectively as follow: $B^{k+1}$ ({\it resp.} $Z^{k+1}$) is the part of $W^{k+1}$ properly intersected by ({\it resp.} contained in) $H_k$, so that 
$$
W^{k+1}=B^{k+1}+Z^{k+1} \enspace, 
$$
while $ W^k:= B^{k+1}\cdot H_k$. The $Z^{k+1}$'s form the locus of improperness of the successive intersections of $W$ with the linear forms $\ell_n,\ell_{n-1},\dots ,\ell_0$. Since points $\xi\in Z^{k+1}$ cannot appear as (isolated) components of $|W|\cap|Z(\bfell)|$ of dimension $0$, it results that $|W^0|$ is actually the zero dimensional part of $|W|\cap|Z(\bfell)|$. Restricting to $M=(S\times\P)\setminus\cup_{k=0}^nZ^{k+1}$, all the intersections $B^{k+1}|_M\cdot H_k|_M$ are proper and the associativity of the intersection product together with the invariance of the multiplicity with respect to localization give 
\begin{equation}\label{equalcycles}
W^0 = W|_M\cdot Z(\bfell)|_M = W\cdot Z(\bfell)\enspace. 
\end{equation}

Consider the indices 
$$
c_{k}:=(\underbrace{1, \dots, 1}_{\scriptstyle k \ {\rm times}},
\underbrace{0,\dots,0}_{\scriptstyle n+1-k \ {\rm times}})\in \N^{n+1}
 \quad \mbox{ for } 0\le k\le n+1 \enspace.
$$
The multihomogeneous B\'ezout theorem~\cite[chap.~7, thm.~3.4]{Rem01b} implies 
\begin{equation}\label{bezoutak}
\deg_{c_k}(W^k) = \deg_{c_{k+1}}(B^{k+1}) \le \deg_{c_{k+1}}(W^{k+1})\enspace,
\end{equation}
and so $\deg(W\cdot Z(\bfell))= \deg(W^0)\le \deg_{(1,\dots, 1)}(W)$ as stated. We also note that in order to have equality in~(\ref{bezoutak}) it suffices that $\deg_{c_{k+1}}(Z^{k+1})=0$.

\medskip
Now assume that $\dim(|W|\cap|Z(\bfell)|)=0$, we have to check $\deg(W^0) = \deg_{(1,\dots,1)}(W)$ that is, with the above remark, $\deg_{c_k}(Z^k)=0$ for all $k= 1,\dots,n+1$.  

But, for $0\le k\le n$ the improper intermediate component $Z^{k+1}$ is  contained in $|W|\cap|Z(\ell_{k}, \dots, \ell_n)|$, and so 
$$
|Z^{k+1}|\cap|Z(\ell_0,\dots, \ell_k)| =\emptyset\enspace,
$$
otherwise  $|W|\cap|Z(\bfell)|$ would contain $|Z^{k+1}|\cap|Z(\ell_0,\dots,\ell_{k-1})|$ which is of positive dimension, contrary to the hypothesis. We affirm that this implies that $\deg_{c_{k+1}}(Z^{k+1})=0$: suppose this is not the case and construct inductively for $i=k+1,\dots, 0$ a variety  $Y^i \subset |Z^{k+1}|\cap|Z(\ell_i, \dots, \ell_{k})|$ such that $\deg_{c_i}(Y^i)>0$. First we take $Y^{k+1}$ to be any component of $Z^{k+1}$ of positive  $c_{k+1}$-degree, then we let $0\le i\le k$ and suppose that $Y^{i+1}$ is already constructed. We take a linear form $\ell'_i\in\K[\bfx_i]$ that cuts $Y^{i+1}$ properly, the multihomogeneous B\'ezout theorem implies 
$$
\deg_{c_i}(Y^{i+1}\cdot Z(\ell'_i)) = \deg_{c_{i+1}}(Y^{i+1}) >0
$$
and we take $Y^i$ to be any component of $ Y^{i+1}\cdot Z(\ell_i')$ of positive $c_i$-degree. We end up with a $0$-dimensional variety $Y^0\subset |Z^{k+1}|\cap|Z(\ell_0,\dots, \ell_k)|$, which is certainly not empty. It is a contradiction which establishes  $\deg_{c_{k+1}}(Z^{k+1})=0$ for all~$k$, and this concludes the proof. 
\end{proof}

\begin{rem}
We may still have equality in lemma~\ref{Bezoutthm} even when the intersection
$|W|\cap|Z(\ell_0,\dots,\ell_n)|$ has positive dimension. With notation as in the above proof, this can happen for instance if some $Z^{k+1}$ such that $\deg_{c_{k+1}}(Z^{k+1})=0$ survives after intersection with $\ell_0,\dots,\ell_{k}$. 
\end{rem}

\smallskip

In the sequel we set the notation~\label{notation} for the rest of this section and the following one. As in the setting of theorem~\ref{genthm}, we assume that for each  $0\le i\le n$ we are given a line bundle $L_i$ together with a Laurent polynomial $f_i \in \Gamma(S;L_i)[\bft^{\pm 1}]\setminus \{0\}$. As remarked in the introduction, this presupposes that $\Gamma(S;L_i)\ne 0$ or equivalently $\deg(L_i)\ge 0$. We write 
$$
f_i =\sum_{j=0}^{N_i} \sigma_{i,j} \bft^{a_{i,j}}
$$
for some $ \sigma_{i,j}\in \Gamma(S; L_i)\setminus \{0\}$ and $a_{i,j}\in \Z^n$. 
We fix a non-zero global section $\rho_i$ of $\Gamma(S;L_i)$, for instance $\rho_i:=\sigma_{i,0}$, then we set $\alpha_{i,j}:=\rho_i^{-1}\sigma_{i,j}$  which is a rational function on $S$. Put
$$
\cA_i:=(a_{i,j}: 0\le j\le N_i) \in (\Z^n)^{N_i+1} 
\quad , \quad 
\alpha_i:=(\alpha_{i,j}: 0\le j\le N_i)\in (\K(S)^\times)^{N_i+1} 
$$
then $\geom:=(\cA_0,\dots, \cA_n) $, $\bfalpha:=(\alpha_0,\dots,\alpha_n)$ and 
$\arith:=(\geom,\bfalpha)$. In this section we will {\it not} assume $L_\geom=\Z^n$, unless otherwise explicitly stated. 

Consider the map $\varphi_\arith:S\times\T^n \to S\times\P$ and the variety 
$X_\arith$ associated to the data $\arith$ as explained in subsection~\ref{torusaction}; both are independent of the choice of the $\rho_i$'s. Now for each $v\in S$ take a further section 
$$
\rho_{i,v}\in \Gamma(S;L_i)
$$ 
such that $\ord_v(\rho_{i,v}) = \ord_v(f_i) = \min(\ord_v(\sigma_{i,j}): 0\le j\le N_i)$, for instance $\rho_{i,v}=\sigma_{i,j(i,v)}$ for some index $j(i,v)$ realizing this minimum. In a neighborhood of a  point $(v,\bft)\in S\times\T^n$, the regular map $\varphi_\arith$ can then be written
$$
(s,\bft)\mapsto\varphi_\arith(s,\bft)=\big(s, ((\rho_{i,v}^{-1}\sigma_{i,j}(s)\, \bft^{a_{i,j}}:i,j)\big) \in S\times \P
\enspace.$$

\smallskip
Next consider the linear form $\ell_i:=\sum_{j=0}^{N_i}x_{i,j}$, which can be interpreted as a global section of the line bundle $\pi_i^*(O(1))$, pull-back of the universal line bundle of $\P^{N_i}$ {\it via} the projection $\pi_i:S\times \P\to \P^{N_i}$. Then $\varphi_\arith^* ( \pi_i^*(O(1)))$ trivializes over $S\times \T^n$ and so its sections are functions of $S\times\T^n$, in particular $\varphi_\arith^*(\ell_i)\in \cO_{S\times\T^n}$. Setting
$$
f_{i,v} := \rho_{i,v}^{-1}f_i\in\cO_{S\times\T^n,\{v\}\times\T^n}\enspace, 
$$ 
we verify 
\begin{equation}\label{formesection}
f_{i,v}^{-1}\varphi^*_\arith(\ell_{i})\in \cO_{S\times\T^n, \{v\}\times\T^n}^\times\enspace,\end{equation}
as this is a regular function which does not vanish on $\{v\}\times\T^n$, and so 
$f_{i,v}$ is an equation for the divisor of $\varphi_\arith^*(\ell_i)$ in a neighborhood of the fiber $\{v\}\times\T^n$. 

\medskip

The following result allows us to treat the intersection multiplicities in 
$S\times \T^n$ by passing to the variety $X_\arith$. 

\begin{lem}\label{multiplicites locales}
Let notation be as above and assume $L_\geom=\Z^n$. 
Set $\bfell:=(\ell_0,\dots, \ell_n)$ and let ${\bfx}=(v,\bfy)$ be an isolated point of $X_\arith^\Ef\cap Z(\bfell)$, then
$${\rm mult}({\bfx}|X_\arith; Z(\bfell)) = \sum_{\xi\in\varphi_\arith^{-1}({\bfx})}
{\rm mult}(\xi| f_{0,v},\dots,   f_{n,v})
\enspace.$$
\end{lem}
\begin{proof}
We consider the local rings
$$
A:=\cO_{X_\arith,\bfx} \quad, \quad B:=\cO_{S\times\T^n,\varphi_\arith^{-1}(\bfx)}
$$ 
of $X_\arith$ at $\bfx$ and of $S\times\T^n$ at $\varphi_\arith^{-1}(\bfx)$ respectively, together with the ring homomorphism $\varphi_\arith^*:A\hookrightarrow B$. We have that $A$ is a reduced local ring with maximal ideal ${\mathfrak m}$ corresponding to the point ${\bfx}$, whereas $B$ is a semi-local finite extension of $A$, according to lemma~\ref{isomorphismegenerique}. This lemma also implies that $S\times \T^n$ and $X_\arith$ are birationally equivalent, hence the field of fractions of $A$ and $B$ coincide. On the other hand, for each maximal ideal ${\mathfrak n}$ of $B$ we have $A\cap{\mathfrak n}={\mathfrak m}$ and $B/{\mathfrak n}\simeq A/{\mathfrak m}\simeq \K$, and therefore  the residual extension is of degree $1$.

Consider now the  ${\mathfrak m}$-primary ideal ${\mathfrak q}=(\ell_{0},\dots, \ell_{n}) \subset A$. By the previous considerations, we are in the hypothesis of~\cite[chap.~VIII, \S~7.3, prop.~6, pp.~75-76]{Bou75} from which results the equality of Samuel multiplicities $e_{\mathfrak q}^A(A) = e_{{\mathfrak q}B}^B(B)$. Besides, \cite[chap.~VIII, \S~7.1, cor., p.~73]{Bou75} implies $e_{{\mathfrak q}B}^B(B) = \sum_{{\mathfrak n}} e_{{\mathfrak q}B_{{\mathfrak n}}}^{B_{{\mathfrak n}}}(B_{{\mathfrak n}})$, the sum running over all the maximal ideals of $B$, and so 
$$ 
e_{{\mathfrak q}}^A(A) = \sum_{{\mathfrak n}} e_{{\mathfrak q}B_{{\mathfrak n}}}^{B_{{\mathfrak n}}}(B_{{\mathfrak n}}) \enspace. 
$$

Since ${\bfx}$ is an isolated point of $X_\arith^\Ef\cap Z(\bfell)$, the intersection multiplicity of $X_\arith^\Ef$ and $Z(\bfell)$ at $\bfx$ is by definition $e_{{\mathfrak q}}^A(A)$. On the other hand, (\ref{formesection}) implies  
$$
{\mathfrak q}B= \varphi_\arith^*({\mathfrak q})= (f_{0,v}, \dots, f_{n,v}) \subset B \enspace, 
$$
and so any $\xi\in\varphi_\arith^{-1}({\bfx})$ is an 
isolated point of $Z(f_{0,v},\dots,f_{n,v})$, defined by some maximal ideal ${\mathfrak n}$ of $B$, and the corresponding multiplicity is $e_{{\mathfrak q}B_{{\mathfrak n}}}^{B_{{\mathfrak n}}}(B_{{\mathfrak n}})$, this concludes the proof.\end{proof}

\begin{cor} \label{cormultiplicites locales} 
With the notation introduced, 
$$
\deg\big(Z\big(\varphi^*_\arith(\bfell)\big)_0\big) = \sum_{v\in S} \kern5pt \sum_{(v,\bft) \in|Z( f_{0,v},\dots, f_{n,v})_0|} \kern-15pt {\rm mult}((v,\bft) | f_{0,v},\dots, f_{n,v}) = \deg(X^\Ef_\arith\cdot Z(\bfell)) 
\enspace.$$
\end{cor}

\begin{proof}
The first equality results directly from the fact that $f_{i,v}$ is an equation of the divisor cut by $\varphi^*_\arith(\ell_i)$ on $S\times\T^n$ in a neighborhood of
$\{v\}\times\T^n$, for any $v\in S$. And the second equality is a consequence of the previous lemma, since by lemma~\ref{isomorphismegenerique} the isolated points of $Z(\varphi^*_\arith(\bfell))$ are exactly the inverse image by $\varphi^{-1}_\arith$ of the isolated  points of $X^\Ef_\arith\cdot Z(\bfell)$.
\end{proof}

Recall that for $0\le i\le n$  and $v\in S$ we denote by 
$\overline{\vartheta}_{v}(f_i)$ the constant function  $\ord_v(f_i)$ on the polytope 
$\NP(f_i(v,\cdot))$.
On the other hand, set 
$$
B(f_i):=\{v\in S: \ord_v(f_i)>0\}\subset S
$$
for the {\it base locus} of $f_i$. For $v\in B(f_i)$, 
the evaluation $f(v,\cdot)$ is zero and so $\NP(f(v,\cdot))=\{0\}$ by convention. The following lemma will allow us to control the contribution of the base points to the intersection of the $f_i$'s.  

\begin{lem}\label{lemmultpointsbase}
With the notation introduced, for  $v\in S$ we have 
$$
\sum_{\substack{\bft\in \T^n\\(v,\bft)\in|Z(\bff)|_0}} \big({\rm mult}((v,\bft)|\bff) - {\rm mult}((v,\bft)| f_{0,v},\dots, f_{n,v})\big) \leq \MI_n(\overline{\vartheta}_{v}(\bff))
\enspace,$$
with equality if and only if for all $0\le i\le n$ and $v\in S$ 
$$
Z(\bff)_0 = Z(f_{0,v},\dots, f_{n,v})_0 + Z(f_0, \dots, f_{i-1}, \rho_{i,v}, f_{i+1},\dots, f_n)_0 
$$
in a neighborhood of $\{v\}\times\T^n$ and
$$
\sum_{\bft\in \T^n} \mult(\bft| f_k(v,\cdot) : k\ne i) = \MV_n(\NP(f_k(v,\cdot)): k\ne i) \enspace.
$$
\end{lem}

\begin{proof}
In case $v\in B(f_i)\cap B(f_k)$ for some $i\ne k$, then $Z(\bff)\cap (\{v\}\times \T^n) = Z(\bff(v,\cdot))$ is a subset of $\{v\}\times \T^n$ defined by $ \le  n-1$ equations, therefore $Z(\bff)$ has no isolated components above $v$ and the left-hand side of the inequality is zero. Besides $f_i(v,\cdot)=f_k(v,\cdot)=0$ and so $\NP(f_i(v,\cdot))=\NP(f_k(v,\cdot))=\{0\}$, which implies that the mixed integral is zero because of  formula~(\ref{masure_gorbeau}) and the basic properties of the mixed volume. Hence 
the  inequality reduces to $0=0$.

On the other hand, if $v$ is not a base point of any of the $f_i$'s, then the $f_i$'s and the $f_{i,v}$'s define the same cycles in a neighborhood of $\{v\}\times \T^n$, and so the left-hand side is zero. For the right-hand side, we have $\ov\vartheta_{v}(f_i)\equiv 0$ for all $i$ and so the corresponding mixed 
integral is zero, also by formula~(\ref{masure_gorbeau}). The inequality reduces again to $0=0$.

Hence, the only interesting  case is when $v\in B(f_i)$ for exactly one $i$. We will assume without loss of generality $i=0$, up to a reordering of the indices. With this assumption,  $\rho_{k,v}(v)\ne 0$ for all $1\le k\le n$ and so $Z(\bff) = Z(f_{0,v},\dots, f_{n,v}) + Z(\rho_{0,v}, f_1,\dots, f_n)$ which implies that in a neighborhood of $\{v\}\times\T^n$ 
$$
Z(\bff)_0 \subset Z(f_{0,v},\dots, f_{n,v})_0 + Z(\rho_{0,v}, f_1,\dots, f_n)_0 \enspace.  
$$
This shows that the sum of ${\rm mult}((v,\bft)|\bff) - {\rm mult}((v,\bft)| f_{0,v},\dots, f_{n,v}) $ over $\bft\in \T^n$ such $(v,\bft)\in Z(\bff)_0 $ is bounded above by 
\begin{align*}
\ord_v(\rho_{0,v}) \sum_{\bft\in \T^n} {\rm mult}( \bft| f_k(v,\cdot): 1\le k\le n) 
&\le \ord_v(f_0)\MV_n\big(\NP(f_k(v,\cdot)): 1\le k\le n\big) \\[-1mm]
& \le \MI_n(\overline{\vartheta}_{v}(\bff)) \enspace. 
\end{align*}
The first estimate is the Ku\v snirenko-Bern\-\v{s}tein's theorem, while the last one is just formula~(\ref{masure_gorbeau}) again.
\end{proof}

In what follows we establish the upper bound in theorem~\ref{genthm}. The study of the conditions for  this estimate to be  exact is postponed to section~\ref{bk_equality}. 

\medskip
\begin{Proof}{\it Proof of theorem~\ref{genthm}.\/}
By summing the estimate in lemma~\ref{lemmultpointsbase} over $v\in S$ together with corollary~\ref{cormultiplicites locales}, we obtain 
\begin{equation}\label{sumoverS} 
\sum_{\xi\in|Z(\bff)|_0}{\rm mult}(\xi|\bff) \leq \deg\big(Z(\varphi_\arith^*(\bfell))_0\big)  + \sum_{v\in S} \MI_n(\overline{\vartheta}_{v}(\bff)) \enspace.
\end{equation}

Suppose for the moment $L_\arith=\Z^n$. Applying successively corollary~\ref{cormultiplicites locales}, lemma~\ref{Bezoutthm} and proposition~\ref{fauchelevent} we get
\begin{align}\label{dspreuvemainthm}
\deg(Z(\varphi_\arith^*(\bfell))_0)
&= \deg(X_\arith^\Ef\cdot Z(\bfell)) \\[1.5mm]\nonumber
&\leq \deg_{(1,\dots,1)}(X_\arith)\\[1.5mm]\nonumber
&\leq \sum_{v\in S} \MI_n(\vartheta_{v}(\rho_0^{-1}f_{0}), \dots, \vartheta_{v}(\rho_n^{-1}f_{n}))
\enspace,
\end{align}
with equality if $X_\arith\cap|Z(\bfell)|$ is of dimension $0$ and entirely contained in $X_\arith^\Ef\cap|Z(\bfell)|$. By definition of $\vartheta_v$ we have $\vartheta_{v}(\rho_i^{-1}f_i)= \vartheta_{v}(f_i) + \ord_v(\rho_i)$ and so by formula~(\ref{multiintegralefonctionstranslates}) it comes that $\MI_n(\vartheta_{v}(\rho_0^{-1}f_0),\dots,\vartheta_{v}(\rho_n^{-1}f_n))$ equals 
$$
\MI_n(\vartheta_{v}(\bff)) + \sum_{i=0}^n \ord_v(\rho_i)\MV_n(\NP(f_k): 0\le k\le n, k\ne i) 
\enspace.
$$
By summing over $S$ and applying formula~(\ref{masure_gorbeau})
we find
\begin{align*}
\sum_{v\in S}
\sum_{i=0}^n \ord_v(\rho_i)\MV_n(\NP(f_k): k\ne i) 
&= 
\sum_{i=0}^n \deg(L_i)\MV_n(\NP(f_k): k\ne i) \\[1mm]
&= 
\MI_n(\bfdelta|_{\NP(\bff)}) 
\end{align*}
because $\sum_{v} \ord_v(\rho_i)=\deg(L_i)$. Therefore,
\begin{equation}\label{pullback}
\deg(Z(\varphi_\arith^*(\bfell))_0) \leq \MI_n(\bfdelta|_{\NP(\bff)}) + \sum_{v\in S} \MI_n(\vartheta_v(\bff)) \enspace,  
\end{equation}
which together with~(\ref{sumoverS}) 
proves the estimate~(\ref{genthm_display}) for $L_\geom=\Z^n$. 

\smallskip 
In case $\rank_\Z(L_\geom)=n$ but $L_\geom\not=\Z^n$, we can reparameterize the toric variety as explained 
in the diagram~(\ref{diagrammereduc}). 
With the notation therein, for $\zeta=(v,\bft)\in S\times \T^n$ 
$$
\sum_{\xi\in \psi^{-1}(\zeta)} \mult(\xi| f_{0,v},\dots, f_{n,v}) = 
[\Z^n:L_\geom] \, \mult\big(\zeta| (\psi^*)^{-1}(f_{0,v}),\dots, (\psi^*)^{-1}(f_{n,v})\big)
$$
since $\psi$ is finite flat map of degree $[\Z^n:L_\geom]$, {\it see} for 
instance~\cite[chap.~VIII, \S~7.2, prop.~4, p.~73]{Bou75}. But, $\psi^*\circ\varphi^*_{(\bfcalB,\bfalpha)}(\ell_i)=\varphi^*_\arith(\ell_i)$ and $(\psi^*)^{-1}(f_{i,v})$ is an equation of the divisor of $\varphi_{(\bfcalB,\bfalpha)}^*(\ell_i)$ on $\{v\}\times\T^n$, whence $\deg(Z(\varphi_\arith^*(\bfell))_0) = [\Z^n:L_\geom] \deg(Z(\varphi_{(\bfcalB,\bfalpha)}^*(\bfell))_0) $. 

On the other hand, the functions $\bfdelta|_{\NP(\cdot)}$ and $\bfvartheta_v(\cdot)$ relative to the data ${(\bfcalB,\bfalpha)}$ are just the ones corresponding to the data $\arith$ composed with the linear transformation $\ell$ associated to map $\psi$. Proposition~\ref{transfolineaire} then shows that their mixed integrals
relative to $\arith$ are $|\det(\ell)|=[\Z^n:L_\geom]$ times those corresponding to 
$(\bfcalB,\bfalpha)$. These observations show that if inequality~(\ref{pullback})
is valid for the data~$(\bfcalB,\bfalpha)$ then it is also valid for the data $\arith$. Together with~(\ref{sumoverS}), this proves the estimate~(\ref{genthm_display}) for $\rank_\Z(L_\geom)=n$. 

\smallskip 
Finally, in case $\rank_\Z(L_\geom)<n$, by the results in subsection~\ref{dimensions} we have 
$$
\dim(X_\arith)=\rank_\Z(L_\geom)+1<n+1=\dim(S\times\T^n)\enspace.
$$ 
From the theorem of dimension of fibers, it results that the fibers of $\varphi_\arith$ are either empty or positive dimensional, and in either case $Z(\varphi_\arith^*(\bfell))$ has no isolated points and so the first term in the estimate~(\ref{sumoverS}) is zero. On the other hand, the mixed integrals in the second term 
in this estimate are also zero, because the domains of the functions are contained in translates of a single proper linear subspace of $\R^n$. Thus~(\ref{sumoverS}) reduces to zero in this case, which implies that $Z(\bff)$ has no isolated points.
The same arguments show that the mixed integrals occurring in the estimate~(\ref{genthm_display}) are zero, hence this estimate also reduces to $0=0$, which completes the proof. 
\end{Proof}

\begin{rem}\label{remarqueMI}
Setting $B_i:=B(f_i)\setminus \cup_{k\ne i} B(f_k) $, 
it results from the proof of lemma~\ref{lemmultpointsbase} that
we can express the contribution of the functions $\ov{\vartheta}_v(\bff)$ to the estimate in theorem~\ref{genthm} in terms of mixed volumes as 
$$
\sum_{v\in S} \MI_{n}(\ov{\vartheta}_v(\bff)) = \sum_{i=0}^n\sum_{v\in B_i}\ord_v(f_i)\MV_n\big(\NP(f_k(v,\cdot)):0\le k\le n,k\not=i\big) 
\enspace.$$
\end{rem}

\medskip

\begin{Proof}{\it Proof of theorem~\ref{mainthm}.\/}
Let $f_0,\dots,f_n\in\k[s][\bft^{\pm 1}]$ be a family of  Laurent polynomials. For $0\le i\le n$ set $\delta_i$ for the partial degree of $f_i$ in the variable $s$, then 
$$
F_i:= \sigma^{\otimes \delta_i} f_i
\in \Gamma(\P^1;O(\delta_i))[\bft^{\pm1}] 
$$ 
where $\sigma \in \Gamma(\P^1;O(1))$ denotes the section  corresponding to the point at infinity. Theorem~\ref{genthm} with $S=\P^1$ implies 
\begin{equation}\label{genthmoverP1}
\sum_{\xi\in|Z(\bfF)|_0}{\rm mult}(\xi|\bfF) 
\le \MI_n\left(\bfdelta|_{\NP(\bfF)}\right)+ \sum_{v\in\P^1} 
\left(\MI_n(\vartheta_v(\bfF))
+\MI_n(\ov{\vartheta}_v(\bfF))\right)
\enspace.
\end{equation}
The cycle $Z(\bff)$ is the restriction of $Z(\bfF)$ to 
$\A^1\times\T^n$, and we have 
$
\vartheta_\infty(F_i)= \vartheta_\infty(f_i)-\delta_i$ while 
$\vartheta_v(F_i)= \vartheta_v(f_i) $ for $v\ne \infty$.
By construction, the point at infinity is not a base point 
of any of the $F_i$'s and so  
$\ov{\vartheta}_\infty(F_i)=0 $ while 
$\ov{\vartheta}_v(F_i) =\ov{\vartheta}_v(f_i)$ for $ v\ne \infty$. 
We thus obtain $\MI_n(\vartheta_v(\bff)) = \MI_n(\vartheta_v(\bfF))$ and $\MI_n(\ov{\vartheta}_v(\bff)) = \MI_n(\ov{\vartheta}_v(\bfF))$ for $v\not=\infty$, $\MI_n(\ov{\vartheta}_\infty(\bfF))=0$, $\MI_n(\vartheta_\infty(\bff)) = \MI_n\left(\bfdelta|_{\NP(\bfF)}\right) + \MI_n(\vartheta_\infty(\bfF))$ and the estimate 
\begin{equation}\label{mainthm_whithbasepoints}
\sum_{\xi\in |Z(\bff)|_0}{\rm mult}(\xi|\bff) 
\le  
\MI_{n}(\vartheta_\infty(\bff)) + 
\sum_{v\in\A^1} 
\left(\MI_{n}(\bfvartheta_v(\bff))
+\MI_{n}(\ov{\vartheta}_v(\bff))\right) \enspace, 
\end{equation}
with equality when~(\ref{genthmoverP1}) is an equality and moreover 
$Z(\bfF)$ has no isolated points above $\infty$.

In the setting of theorem~\ref{mainthm}, the hypothesis that the $f_i$'s are primitive is equivalent to $B(f_i)=\emptyset $ for all $i$. This implies that all of the mixed integrals $\MI_n(\ov{\vartheta}_v(\bff))$ are zero and so~(\ref{mainthm_whithbasepoints}) reduces to the estimate~(\ref{mainthm_display}). 
\end{Proof}

\medskip
Let $f_0,\dots,f_n \in\k[s^{\pm 1},\bft^{\pm 1}]$ and set $P_i\subset\R^{n+1}$  for 
the Newton polytope of $f_i$ with respect to {\it all} of the variables $s$ and $\bft$. Set $d_i\in \Z$ for the minimal exponent such that $s^{d_i}f_i\in\K[s][\bft^{\pm1}]$. In this situation, neither $0$ nor $\infty$ is a base point of any of the $s^{d_i}f_i$'s and so $ \ov{\vartheta}_0(s^{d_i}f_i)=\ov{\vartheta}_\infty(s^{d_i}f_i)=0$. Setting $s^{\bfd}\bff:= (s^{d_0}f_0,\dots, s^{d_n}f_n)$, inequality~(\ref{mainthm_whithbasepoints}) and  proposition~\ref{bagne} imply 
\begin{align}\label{comparisonwithBK}
\nonumber \sum_{\xi\in|Z(s^{\bfd}\bff)|_0}{\rm mult}(\xi|s^{\bfd}\bff) 
\le &\MI_n(\vartheta_0(\bff)) + \MI_n(\vartheta_\infty(\bff))  \\[-3mm] 
\nonumber &+\sum_{v\in\T^1} 
\big(\MI_n(\vartheta_v(s^{\bfd}\bff))
+\MI_n(\ov{\vartheta}_v(s^{\bfd}\bff))\big)\\[1mm]
\le &\MV_{n+1}(\bfP)+\sum_{v\in\T^1} \left(\MI_n(\vartheta_v(\bff))+\MI_n(\ov{\vartheta}_v(\bff))\right). 
\end{align}
The set of common zeros in $\T^{n+1}$ of the $f_i$'s coincides with that of the $s^{d_i}f_i$'s and furthermore $\MI_n(\vartheta_v(\bff)) + \MI_n(\ov{\vartheta}_v(\bff)) \le 0$ for all $v \in \T^1$. This shows that~(\ref{mainthm_display}) improves upon Ku\v snirenko-Bern\v{s}tein's estimate, besides the fact that it counts the isolated roots of the (modified) system in a set larger than $\T^{n+1}$. 

\medskip

\begin{Proof}{Proof of corollary~\ref{corgenthm}.}
For $m+1\le k\le n$ let $E_k$ be a generic hyperplane of $\P^{N_k}$. Write
$E_k=Z\big(\sum_{j=0}^{N_k}\ell_{k,j}x_{k,j}\big)$ for generic $\ell_{k,j}\in \K$ and write also
$g_k=\sum_{j=0}^{N_k} \sigma_{k,j}\bft^{a_{k,j}}$ for some $\sigma_{k,j}\in \Gamma(S;L_k)\setminus \{0\}$ and $a_{k,j}\in \Z^n$. 
Then 
$$
\varphi_{g_k}^{-1}(E_k)= Z(g_k')
$$
for $g_k':=
\sum_{j=0}^{N_k} \ell_{k,j}\sigma_{k,j}\bft^{a_{k,j}}$
and so 
$ \deg_{\bfg}(Z(\bff)_{n-m})= \deg( Z(\bff)_{n-m}\cdot Z(\bfg'))$ is the quantity we want to estimate. 

The hypothesis that the $g_k$'s have no base point implies that $Z(\bfg')$ cuts properly any (fixed) set of $S\times\T^n$. In particular, $Z(\bfg')$ cuts properly $Z(\bff)_{n-m}$ and avoids the locus of improperness of the intersection of the $f_i$'s in $Z(\bff)_{n-m}$ and so $Z(\bff)_{n-m}\cdot  Z(\bfg') = Z(\bff, \bfg')_0$. Theorem~\ref{genthm} gives then  the result, since the Newton polytope and $v$-adic Newton polytopes of each $g'_k$ coincide with that of $g_k$. 
\end{Proof}

\medskip
Finally, we extend theorem~\ref{genthm} to the singular case. Thus  we now suppose that $S$ is a complete but possibly singular curve. In this more general setting, the points of~$S$ have to be replaced by its places. In the sequel we quickly review the definitions and basic facts about these places of a curve, the details can be found in~\cite[\S~IV.2]{Wal50}.  

\smallskip 
A {\it parameterization of $S$} is a non constant analytic map $g:\K\to S$ and the point $g(0)\in S$ is called the {\it center} of the parameterization. Two parameterizations $g,h$ are {\it equivalent} if there exists a local isomorphism $\zeta: \K\to\K $ at~$0$ such that $g=h\circ \zeta$. A parameterization $g$ is said {\it irreducible} if it is injective in a neighborhood of $0$. By definition, a {\it place of $S$} is an  equivalence class of irreducible parameterizations and we denote by $V_S$ the set of all places of $S$. For $v\in V_S$ we note $g_v$  some parameterization corresponding to  $v$; the map $V_S\to S, v\mapsto g_v(0)$ is then well-defined and surjective. 

\smallskip
Given a line bundle  $L$ of $S$, the {\it order of vanishing} $\ord_v(\sigma)$ of a section $\sigma\in \Gamma(S;L)$ at a given $v\in V_S$ is defined as the order at $0$ of the analytic map $\sigma\circ g_v: \K\to L$. This definition does not depend on the choice of $g_v$. Thus for a Laurent polynomial $f \in\Gamma(S;L)[\bft^{\pm1}]$ and a place $v\in V_S$ we can  extend in the natural way the notions of $v$-adic Newton polytope $\NP_v(f)\subset \R^{n+1}$ and corresponding functions $\vartheta_v(f)$ and $\overline{\vartheta}_{v}(f)$, respectively defined on the polytopes 
$\NP(f)$ and $\NP(f(g_v(0),\cdot))$.

\begin{thm} \label{gengenthm}
Let $S$ be a complete curve and for $0\le i\le n$ let $L_i$ be a line bundle on $S$ of degree $\delta_i$ and $f_i \in \Gamma(S;L_i)[\bft^{\pm 1}]\setminus\{0\}$, then 
\begin{equation}\label{gengenthm_display}
\sum_{\xi\in|Z(\bff)|_0}{\rm mult}(\xi|\bff) \le 
\MI_{n}\left(\bfdelta|_{\NP(\bff)} \right)+
\sum_{v\in V_S} \left(\MI_{n} (\vartheta_v(\bff)) + \MI_{n}(\overline{\vartheta}_{v}(\bff)) \right)\enspace.
\end{equation}
Furthermore, this is an equality for $\bff$  generic among systems with given functions
$(\vartheta_{i,v}: 0\le i\le n, v\in V_S)$.
\end{thm}

\begin{proof}
Let $\nu:\wt S\to S$ be the normalization morphism of $S$ and for each $i$ we consider the pull-back $\wt L_i:=\nu^*(L_i)$  of $L_i$ to a line bundle of $\wt S$ and $\wt f_i:=\nu^*(f_i)$ the pull-back of $f_i$ to a Laurent polynomial in $\Gamma(\wt S;\wt L_i)[\bft^{\pm1}]$. We have  $\deg(\wt L_i)=\deg(L_i)=\delta_i$ and $\NP(\wt f_i )=\NP(f_i)$, and so applying theorem~\ref{genthm} to the system $\wt\bff=0$ we obtain
$$\sum_{\xi\in|Z(\wt\bff)|_0} {\rm mult}(\xi|\wt\bff) \le \MI_n\big(\bfdelta|_{\NP(\bff)}\big) + \sum_{s\in \wt S} \big(\MI_n(\vartheta_{s}(\wt\bff)) + \MI_n(\ov{\vartheta}_s(\wt\bff))\big)
$$
We are in a situation similar to that of lemma~\ref{multiplicites locales} and, as in the proof of this result, \cite[chap.~VIII, \S~7.3, prop.~6, pp.~75-76]{Bou75} implies that the left-hand side of the above inequality coincides with that of~(\ref{gengenthm_display}).

Now the places of $S$ are in 1-to-1 correspondence with the points of $\wt S$: the bijection is given by $V_S\to \wt S, v\mapsto \wt g_v(0)$, where for a place $v\in V_S$, we denote by $\wt g_v:\K\to \wt S$ the lifting of the parameterization $g_v$.
For a place $v\in S$ and $s(v):=\wt g_v(0)\in \wt S$ and $\sigma\in \Gamma(S;L_i)$ 
we have $\ord_{s(v)}(\nu^*(\sigma))=\ord_v(\sigma)$ and in particular $\vartheta_{s(v)}(\wt\bff)=\vartheta_{v}(\bff)$ and $\ov{\vartheta}_{s(v)}(\wt\bff)=\ov{\vartheta}_{v}(\bff)$, the result follows.
\end{proof}

\bigskip
\section{Equality conditions}\label{bk_equality}

In this section we determine sufficient conditions for the estimates in theorems~\ref{mainthm} and~\ref{genthm} to be exact, in terms of the solvability of some initial systems associated to the input system $\bff$~(proposition~\ref{genegalite} below). 

\bigskip
We place ourselves again in the setting of theorem~\ref{genthm} and we continue to use the notation from the previous section, set up in page~\pageref{notation}. In particular, $S$ is a smooth complete curve equipped with line bundles $L_i$  and  we are given non zero Laurent polynomials $f_i=\sum_{j=0}^{N_i} \sigma_{i,j}\bft^{a_{i,j}}$ with coefficients in $\Gamma(S;L_i)$.

For each  $v\in S$ we fix a parameterization  $g_v$ of $S$ at $v$ and for a rational function $\beta\in \K(S)$ we recall that $\lambda_{g_v}(\beta)\in \K^\times$ denotes its initial coefficient at $v$, as explained in subsection~\ref{parameterizations}. 
Recall that $\rho_i$ is any non-zero global section of $L_i$ and $\alpha_{i,j}$ is the rational function $\rho_i^{-1}\sigma_{i,j}$, so that $\rho_i^{-1}f_i=\sum_{j=0}^{N_i} \alpha_{i,j}\bft^{a_{i,j}}$.

\begin{lem}\label{initial}
Let $0\le i\le n$,  $v\in S$, $\tau\in \R^n$ and set 
$$
\Lambda_i(\bft):=  \sum_{j} \lambda_{g_v}(\alpha_{i,j})\bft^{a_{i,j}} \in\K[\bft^{\pm1}] \enspace,
$$ 
where the sum runs over the $0\le j\le N$ such that $(a_{i,j},-\ord_v(\sigma_{i,j})) \in \NP_v(f_i)^{(\tau,1)}$, then there exists some $c\in \Z$ such that 
$$
\big(\rho_i^{-1} f_i\big)(g_v(z), z^{-\tau_1} t_1, \dots, z^{-\tau_n} t_n) 
= z^c( \Lambda_i(\bft)+ O(z)) \enspace. 
$$
\end{lem}

\begin{proof}
We have 
\begin{align*}
\big(\rho_i^{-1} f_i\big) (g_v(z), z^{-\tau} \bft)
& = \sum_{j=0}^{N_i} \alpha_{i,j}(g_v(z)) z^{- \langle a_{i,j},\tau\rangle}\bft^{a_{i,j}} \\[0mm] 
& 
= \sum_{j=0}^{N_i} \lambda_{g_v}(\alpha_{i,j})  z^{\ord_v(\alpha_{i,j})- \langle a_{i,j},\tau\rangle}\bft^{a_{i,j}}  + \hot \\[1mm]
& = z^c
 \Lambda_i(\bft)+ \hot 
\end{align*}
for $c=-\max\{\langle (a_{i,j},-\ord_v(\alpha_{i,j})), (\tau,1)\rangle: 0\le j\le N\}$, because the exponent 
$$
\ord_v(\alpha_{i,j})- \langle a_{i,j},\tau\rangle= 
-\ord_v(\rho_i)-
\langle (a_{i,j},-\ord_v(\sigma_{i,j})), (\tau,1)\rangle
$$ 
is minimal if and only if $(a_{i,j},-\ord_v(\sigma_{i,j})) \in \NP_v(f_i)^{(\tau,1)}$.
\end{proof}

\medskip
Let $v\in S$ and $\tau \in \R^{n}$, we define the {\em initial part of $f_i$ at $v$ with respect to $\tau$} as the Laurent polynomial $\init_{v,\tau}(f_i):=\Lambda_i(\bft) \in \K[\bft^{\pm1}]$ in the lemma above. Also we set $\init_{v,\tau}(\bff):= (\Lambda_0,\dots, \Lambda_n)$. This initial system can be interpreted as the subsystem of~$\bff$ associated to the slopes determined by $\tau$ in the family of $v$-adic Newton polytopes, {\it see} the example~\ref{exmpl 1} for illustration.

As observed in subsection~\ref{parameterizations}, changing the parameterization $g_v$ acts on $(\lambda_{g_v}(\alpha_{i,j}): i,j)\in \prod_{i=0}^n(\K^\times)^{N_i+1}$ as the 1-parameter action associated to the integer vector $(\ord_v(\sigma_{i,j}):i,j)\in \prod_{i=0}^n\Z^{N_i+1}$. Besides, changing the $\rho_i$'s multiplies each  $\Lambda_i$ by a non-zero scalar factor. Hence the initial system $\init_{v,\tau}(\bff)$ is only well-defined as a point in a multiprojective space $\P=\prod_{i=0}^n\P^{N_i}$ modulo this 1-parameter action. 

\medskip
The main property of these initial systems is that they allow to detect when the linear system $\bfell=0$ intersects a certain orbit of $X_\arith$. Recall that for a family of slopes  $\bfF$ of $\NP_v(\bff)$   we denote by $X_{v,\bfF}$ the corresponding orbit of $X_\arith$, as defined in~(\ref{defXSF}).

\begin{lem}\label{pointsinorbit}
Let notation be as in lemma~\ref{initial} and set $\bfF:= \NP_v(\bff)^{(\tau,1)}$, then $|Z(\bfell)|\cap X_{v,\bfF} = \emptyset$ if and only if the system $\init_{v,\tau}(\bff) = 0$ has no solution in~$\T^n$. 
\end{lem}

\begin{proof}
Write $\bfF=(F_0,\dots,F_n)$, then 
all $\xi\in X_{v,\bfF}$ are of the form $(v,\bfx_{g_v,\bfF})*_\geom \bft$ for some $\bft\in \T^n$. Hence
$\xi_{i,j}= \lambda_{g_v}(\alpha_{i,j}) \bft^{a_{i,j}}$ if $(a_{i,j},-\ord_v(\sigma_{i,j}))\in F_i$ and $\xi_{i,j}=0$ otherwise. 
Thus $\ell_i(\xi)= \Lambda_i(\bft)$ for all $i$ and the result follows. 
\end{proof}

By lemma~\ref{welldefined}, the orbit $X_{v,\bfF}$ does not depend on the choice of the parameterization. 
Hence the previous lemma implies that the solvability of $\init_{v,\tau}(\bff)=0$ on $\T^n$ does not depend on the choice of the parameterization $g_v$ or the global section  $\rho_i$.

\medskip
Recall that  $B(f_i)\subset S$ denotes the base locus of~$f_i$. 

\begin{prop} \label{genegalite}
Let $S$ be a smooth complete curve equipped with line bundles $L_i$  for 
$0\le i\le n$, and for each $i$ let $f_i \in \Gamma(S;L_i)[\bft^{\pm 1}]\setminus\{0\}$. 
Suppose the following hold:
\smallskip 
\begin{enumerate}
\item for all $v\in S$ and $\tau\in\R^n \setminus\{\bfzero\}$, the system of equations $\init_{v,\tau}(\bff)=0$ has no solution in $\T^n$;  

\smallskip
\item in case $S=S^\Ef$, the system of equations
$\bff(\cdot,\bft)\equiv \bfzero$ has no solution  $\bft\in\T^n$;
\smallskip
\item for all $0\leq i<k\leq n$ and 
$v\in S^\Ef\cap B(f_i)\cap B(f_k)$, the system of equations 
$\init_{v,\bfzero}(\bff)=0$ has no solution in $\T^n$. 
\end{enumerate} 
Then $f_0,\dots,f_n$ intersect properly in $S\times \T^n$ and the estimate in theorem~\ref{genthm} is exact. These conditions are satisfied for $\bff$ generic among systems with given functions $(\vartheta_{i,v}:0\le i\le n , v\in S)$. 
\end{prop}

These equality conditions are analogous to those for 
the Ku\v snirenko-Bernstein theorem, 
which can be stated as follows~\cite{Ber75}: 
for 
$h\in \K[\bft^{\pm1}]$ and $\tau\in \R^n$, the {\it initial part of $h$ with respect to $\tau$} is  the Laurent polynomial $\init_\tau(h)\in \K[\bft^{\pm1}]$ such that
$$
h(z^{-\tau_1}t_1,\dots, z^{-\tau_n}t_n) = z^c (\init_\tau(h)(\bft)+ O(z))
$$
for an additional variable $z$ and some $c\in \Z$.  
Let $h_1,\dots, h_n\in \K[\bft^{\pm 1}]$ be a family of Laurent polynomials such that for all $\tau \ne \bfzero$ the initial system
$$
\init_\tau(h_1) = \cdots =\init_\tau(h_n)=0
$$
has no solution in $\T^n$, then the number of solutions 
(counting multiplicities) in $\T^n$ of the system of equations 
$h_1=\cdots=h_n=0$  equals $\MV_n(\NP(h_1),\dots, \NP(h_n))$. 

\medskip 
The following corollary to lemma~\ref{initial}
 allows to detect when the linear system $\bfell=0$ has solutions lying  
outside of $X^\Ef_\arith$.

\begin{cor}\label{propegal}
Suppose $L_\geom=\Z^n$, then  $X_\arith\cap|Z(\bfell)| \subset  X^\Ef_\arith$ if and only if  for all $v\in S$ and $\tau \in\R^n\setminus\{\bfzero\}$,  the system of equations $\init_{v,\tau}(\bff)=0$ has no solution in $\T^n$.
\end{cor}

\begin{proof}
Considering the orbit decomposition~(\ref{orbitdecompositionformula}) the condition $X_\arith\cap|Z(\bfell)| \subset X_\arith^\Ef$ is equivalent to the fact that $X_{v,\bfF}\cap|Z(\bfell)|=\emptyset$ for all $v\in S$ and $\bfF\in\Slope(\NP_v(\bff))$ such that $X_{v,\bfF}\not\subset X_\arith^\Ef$.

On one hand, for $\bfF=\NP_v(\bff)^{(\tau,1)}$ the condition $X_{v,\bfF}\cap|Z(\bff)| = \emptyset$ is equivalent, by lemma~\ref{pointsinorbit}, to the fact that $\init_{v,\tau}(\bff) = 0$ has no solution in $\T^n$.
 
On the other hand, according to proposition~\ref{combeferre}, the orbit $X_{v,\bfF}$ 
corresponding to some  $\bfF\in\Slope(\NP_v(\bff))$ lies in the equivariant set $X_\arith^\Ef$ if and only if both $v\in S^\Ef$ and $\bfF=\NP_v(\bff)^{(\bfzero,1)}$. Therefore, the condition $X_\arith\cap|Z(\bfell)| \subset X_\arith^\Ef$ is equivalent to the fact that the system $\init_{v,\tau}(\bff) = 0$ has no solution in $\T^n$, for any $v\in S\setminus S^\Ef$ and $\tau\in\R^n$ and any $v\in S^\Ef$ and $\tau\in\R^n\setminus\{0\}$.

In the case $s\in S\setminus S^\Ef$, all orbits have to be considered. But then, all components of $\NP_v(\bff)_v^{(\bfzero,1)}$ are in translates of a same subspace of dimension $<n$ and so all families of slopes including this one, are realized by some $\tau \ne \bfzero$. 
\end{proof}

Unfortunately, the condition $X_\arith\cap|Z(\bfell)| \subset  X^\Ef_\arith$ does not warrant that this intersection is of dimension~$0$. Nevertheless, the following lemma shows that the only possible higher dimensional components  are curves of a very special shape. 

\begin{lem}\label{courbeshorizontales}
Suppose $L_\geom=\Z^n$ and let $C\subset S\times\P$ be complete subvariety of positive dimension contained in $X_\arith^\Ef$, then $S=S^\Ef$ and
$C=\varphi_\arith(S\times\{\bft\})$ for some $\bft\in\T^n$. 
\end{lem}

\begin{proof}
Since the map $\varphi_\arith:S^\Ef\times\T^n\to X_\arith^\Ef$ is  
finite and $C$ is complete, 
 $\varphi_\arith^{-1}(C) \subset S^\Ef\times \T^n$ is also 
a complete variety of positive dimension.  
This implies that its projection into $\T^n$ is complete, hence a point. 
This forces $\varphi^{-1}_\arith(C)=S\times\{\bft\}$ for some $\bft\in \T^n$, and so 
$S^\Ef=S$ and 
$C=\varphi_\arith(S\times\{\bft\})$, because $\varphi_\arith$ is a  birational map.
\end{proof}

\medskip

\begin{Proof}{Proof of proposition~\ref{genegalite}.}
As in the proof of theorem~\ref{genthm} in section~\ref{bk_bezout}, we can reduce to the case $L_\geom=\Z^n$. In order to warrant equality in the estimate of this theorem, we have to keep equality in formulas~(\ref{sumoverS}) and~(\ref{dspreuvemainthm}) in  its proof. By lemmas~\ref{Bezoutthm} and~\ref{multiplicites locales} we have equality in~(\ref{dspreuvemainthm}) when $X_\arith\cap|Z(\bfell)|$ is of dimension~$0$ and contained in $X_\arith^\Ef$. But corollary~\ref{propegal} and lemma~\ref{courbeshorizontales} above show that this is equivalent to the joint conditions~(1) and~(2), as condition~(2) excludes exactly all possible components of type $\varphi_\arith(S\times\{\bft\})$ for some $\bft\in \T^n$.

By lemma~\ref{lemmultpointsbase}, to assure equality in~(\ref{sumoverS}) we must insure that for all $0\le i\le n$ and $v\in S$
\begin{equation}\label{dsgenegalitei}
Z(\bff)_0 = Z(f_{0,v},\dots, f_{n,v})_0 + Z(f_0, \dots, f_{i-1}, 
\rho_{i,v}, f_{i+1},\dots, f_n)_0
\end{equation}
in a neighborhood of $\{v\}\times \T^n$ and 
\begin{equation}\label{dsgenegaliteii}
\sum_{\bft\in \T^n} \mult(\bft| f_k(v,\cdot) : k\ne i) = \MV_n(\NP(f_k(v,\cdot)): k\ne i) \enspace. 
\end{equation}

Note that $X_\arith\cap|Z(\bfell)|$ is already of dimension~$0$ and contained in 
$X_\arith^\Ef$ because of conditions~(1) and~(2). Hence for complying with~(\ref{dsgenegalitei}) it suffices to ensure that    $Z(f_0, \dots, f_{i-1}, 
\rho_{i,v}, f_{i+1},\dots, f_n)$ is of dimension~$0$ for all 
 $0\le i\le n$ and $v\in S^\Ef$. 

\medskip 
For $v\in B(f_i)\setminus \bigcup_{k\ne i}B(f_k) $, the equality conditions for the Ku\v snirenko-Bern\v stein estimate imply both~(\ref{dsgenegalitei}) and~(\ref{dsgenegaliteii}) whenever the initial system $\init_\tau(\bff(s,\cdot))=0$ has no solution in $\T^n$ for each $\tau\ne \bfzero$. But this condition is already contained in~(1) since $\init_\tau(\bff(v,\cdot))=\init_{v,\lambda\tau}(\bff)$ for $\lambda>0$ large enough.

Finally, condition~(3) excludes the possibility that there exists a point of $X_\arith^\Ef\cap|Z(\bfell)|$ lying above some $v\in S^\Ef\cap B(f_i)\cap B(f_k)$ for $i\ne k$, and thus~(\ref{dsgenegalitei}) and~(\ref{dsgenegaliteii}) are satisfied for those $v$ as well.

\medskip
Conditions~(1) and~(2) together are equivalent to the fact that  $X_\geom\cap|Z(\bfell)|$ is of dimension $0$ and contained in $X_\arith^\Ef$. Since $\bfell=0$ is a system of $n+1$ linear forms on a variety of dimension $n+1$, this property  holds in the generic case since $X_\arith^\Ef$ contains a dense open set. On the other hand, condition~(3) involves only a finite number of systems of $n+1$ Laurent polynomials in $n$ variables, which are not solvable in the generic case. Hence conditions (1), (2) and (3) are satified for generic $\bff$. 
\end{Proof}

\medskip

\begin{Proof}{Proof of proposition~\ref{equality in mainthm}.}
Theorem~\ref{mainthm} deals with zeros of polynomials in $\A^1\times \T^n$ whereas in theorem~\ref{genthm} the base curve is complete. To achieve equality in theorem~\ref{mainthm} it therefore suffices to write the relevant conditions on the completion of the affine line and then exclude possible zeros above the point at infinity. Note that this last condition will also exclude the possibility of having 
horizontal components of the form $\varphi_\arith(S\times \{\bft\})$. Besides, the assumption that the polynomials are primitive makes vacuous the condition on the zeros of the base loci. It thus only remains condition~(1) from proposition~\ref{genegalite} plus the condition that $\bfell=0$ has no points above $\infty$. 
Lemma~\ref{pointsinorbit} shows that this latter is equivalent to the non-solvability of the system~$\init_{\infty,\bfzero}(\bff)=0$.
\end{Proof}

\bigskip

Although conditions~(1), (2) and~(3) involve infinitely 
many parameters $v$ and $\tau$, 
they can be expressed through a {\it finite} number of systems of  equations
of $\ge n+1$ equations in $n$ variables. We explain how this is done for each condition:

\smallskip
{\leftskip=\parindent\parindent=-\parindent
\un{Condition~(1):} 
above the open subset $S^0$, the orbits of $X_\arith$ can be glued into a finite number of families, as explained in~(\ref{X^0decomposition}). Similarly, the initial systems considered for $v\in S^0$ glue together into a finite number of systems. For each $\bfF'=(F_0,\dots,F_n)\in\Faces(\NP(\bff))$ we have  a system of equations over $S^0\times\T^{n-1}$   
\begin{equation}
\sum_{j:a_{i,j}\in F_i} \sigma_{i,j}(s)\bft^{a_{i,j}}=0 \enspace,\quad \mbox{ for } 0 \le i\le 
n\enspace. 
\end{equation}
The complement $S\setminus S^0$ is finite, and so is the number of initial systems corresponding to those points. 

\smallskip
\un{Condition~(2):} 
for each $0\le i\le n$ choose a basis $\beta_{i,1},\dots,\beta_{i,m_i}$ of the vector space generated by the coefficients $\sigma_{i,j}$, $j=0,\dots,N_i$, and write $\sigma_{i,j}=\sum_{k=1}^{m_i}c_{i,j,k}\beta_{i,k}$. Then the adequate system of equations is
\begin{equation}\label{systemegenerique}
\sum_{j=0}^{N_i}c_{i,j,k}\bft^{a_{i,j}} =0\enspace,\quad
\mbox{ for } 0\le i\le n \mbox{ and } 1\le k \le m_i
\enspace.
\end{equation}

\smallskip
\un{Condition~(3):} 
this is clear,  since the sets $B(f_i)$ are finite. 
\par}

\bigskip

\section{Examples and practical considerations}\label{examples}

In this section we work out a number of examples illustrating different aspects of the presented results. We also include some considerations of practical nature, in particular a procedure to compute our estimate. 

\medskip
\subsection{Dissection of the epitome of all examples}\label{exmpl 1}

Let
\begin{equation}\label{a system}
f= (s-1)+ \, (s-1)^2\, t- 3\, s\, t^2 \quad , \quad g= -7\, (s-1)+(s-1)^2\, t+ 3\, s\, t^2  \quad \in \C[s][t^{\pm1}]\enspace; 
\end{equation}
we verify that the solutions of $f=g=0$ in $\C\times\C^\times$ are the simple root $(4,1)$ and the double one $(-\frac12, -2)$ (adding $f+g$ gives $t$ in terms of $s$). This is an unmixed system. The Newton polytope of both $f$ and $g$ is the interval $[0,2]$ and the following figures show the associated $v$-adic polytopes and their roofs, for each place $v\in \P^1$: 

\begin{figure}[htbp] 
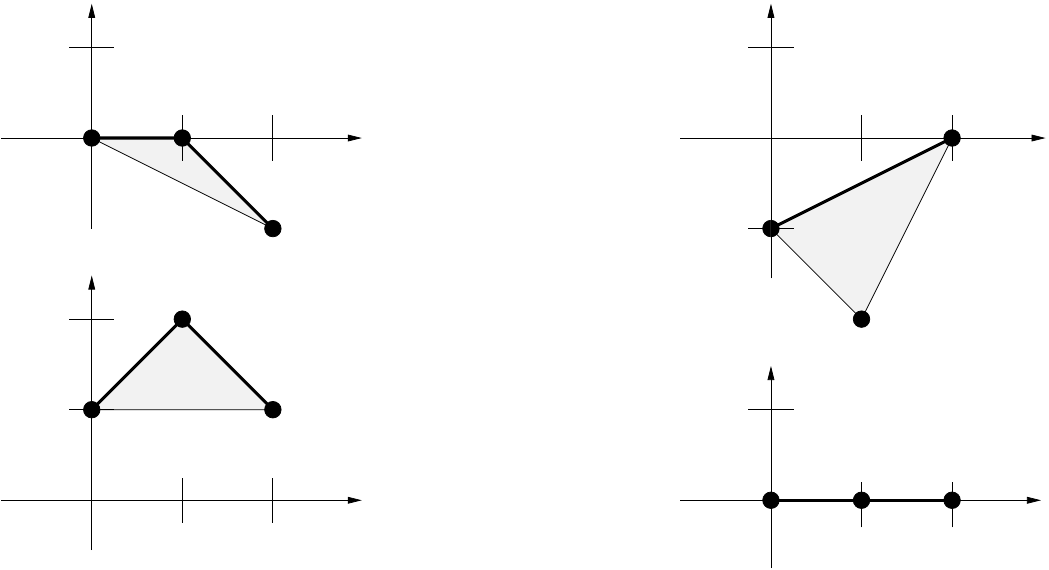
\vspace{-5mm}\caption{}\label{Figure 1}
\end{figure}

\begin{figure}[htbp]
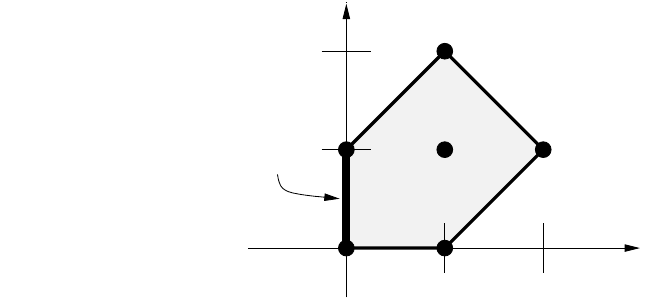
\vspace{-4mm}\caption{}\label{Figure 2}
\end{figure}

\noindent Corollary~\ref{cor_mainthm} gives the estimate (exact in this case)
$$
\sum_{\xi \in|Z(f,g)_0|} \hspace{-2mm}\mult(\xi|f,g) \,  \le \, 2!\, \Big( \int_0^2 \vartheta_0\, \d u +\int_0^2 \vartheta_1\, \d u  + \int_0^2 \vartheta_\infty\, \d u \Big) = 2\,\Big(-\frac12-1+3\Big)=3 \enspace.
$$
On the other hand, the Newton polytope of $f$ and $g$ when regarded as Laurent polynomials in the variables $s$ and $t$ is the pentagon in figure~\ref{Figure 2} below
and so the Ku\v snirenko-Bern\v stein theorem predicts at most $2! \, \Vol_2(P) =5$ solutions;
note that this corresponds to the sum of the contributions of the places $0$ and $\infty$ in the adelic formula in corollary~\ref{cor_mainthm}. Hence this system of equations is generic with respect to our estimate but not with respect to  Ku\v snirenko-Bern\v stein's one.

\medskip
The system~(\ref{a system}) corresponds to a toric surface over $\P^1$ that we denote $X$, which is the Zariski closure of the image of the map
$$
\varphi:\A^1\times\T^1 \to \P^1\times\P^2\quad,\quad (s,t)\mapsto\Big( (1:s), (s-1:(s-1)^2t:st^2)\Big)\enspace.
$$
This map extends to a regular one on $\P^1\times \T^1$, also denoted by $\varphi$. 

This is a hypersurface of $\P^1\times\P^2$ with defining bihomogeneous equation 
$$
(s_1-s_0)^3x_0x_2-s_0^2s_1x_1^2\in \C[s_0,s_1][x_0,x_1,x_2]\enspace.
$$ 
Figure~\ref{Figure 3} below shows this surface in the affine charts $\{s_0\ne 0, x_0\ne 0\}\simeq \A^1\times \A^2 $ (right in red) and  $\{s_1\ne 0, x_0\ne 0\}\simeq \A^1\times \A^2 $ (left in blue), centered at the origin and at the point at infinity of $\P^1$, respectively. The green line represents the projective line (first factor) and the dot the origin  $0=(1:0)$ (blue dot on the right) or the point at infinity  $\infty=(0:1)$  (red dot on the left). The vertical black lines are the axis defined by $x_1=0$ in the affine chart $x_0\ne0$ of $\P^2$ above the point $(1:1)$ and the horizontal ones the axis $x_2=0$ above $0$ (on the right) or $\infty$ (on the left). The two pictures glue as follows: the right part of the red surface corresponds to the part of the blue one comprised between $1$ and $\infty$, the middle part of the red surface between $0$ and $1$ corresponds to the left part of the blue surface and the part of the red surface left to $0$ glues with the part of the blue surface right to $\infty$.

\begin{figure}[htbp]
\hbox to\hsize{
\hfil\pdfximage height 5cm {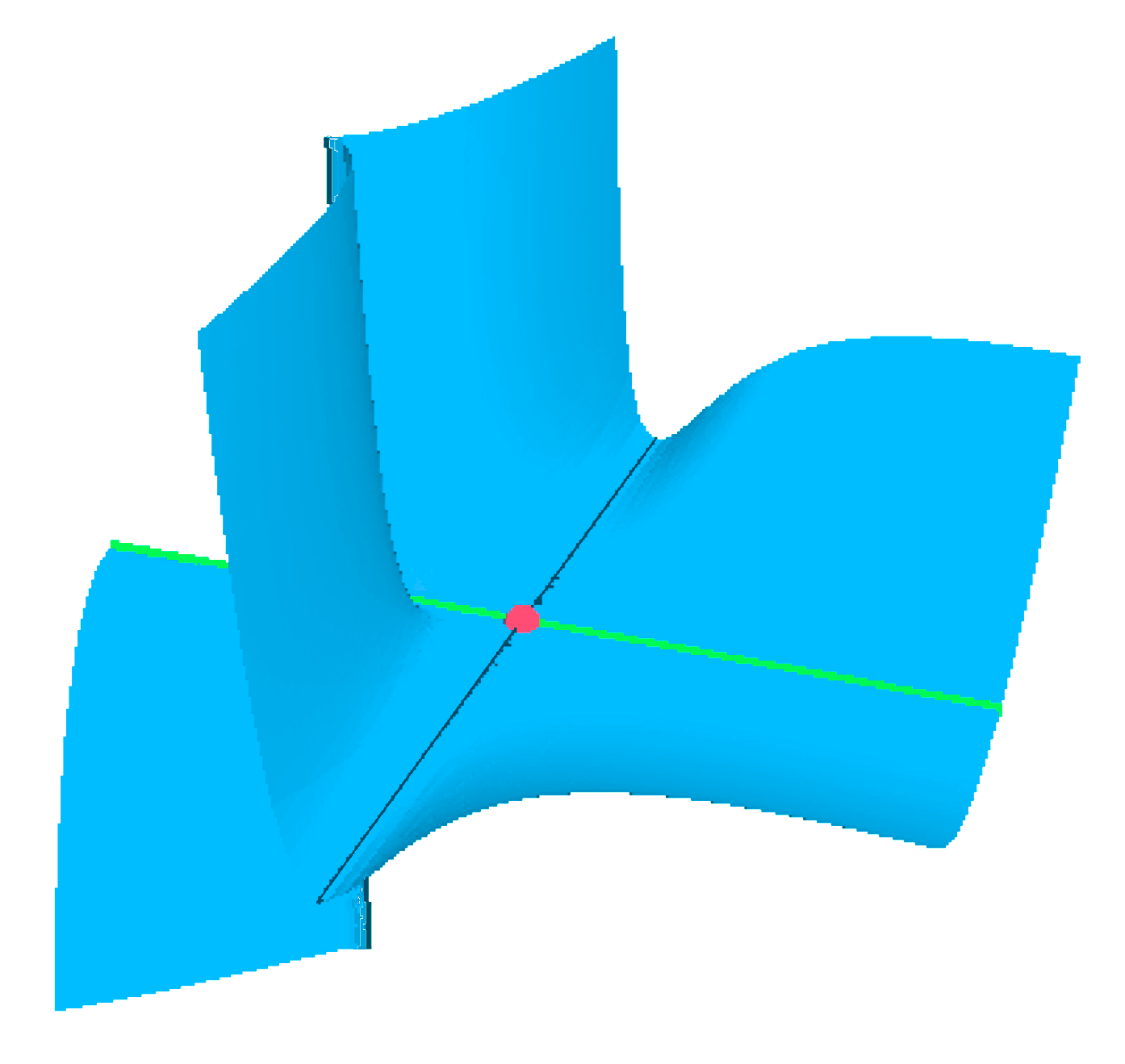}\pdfrefximage\pdflastximage 
\hfil
\pdfximage height 5cm {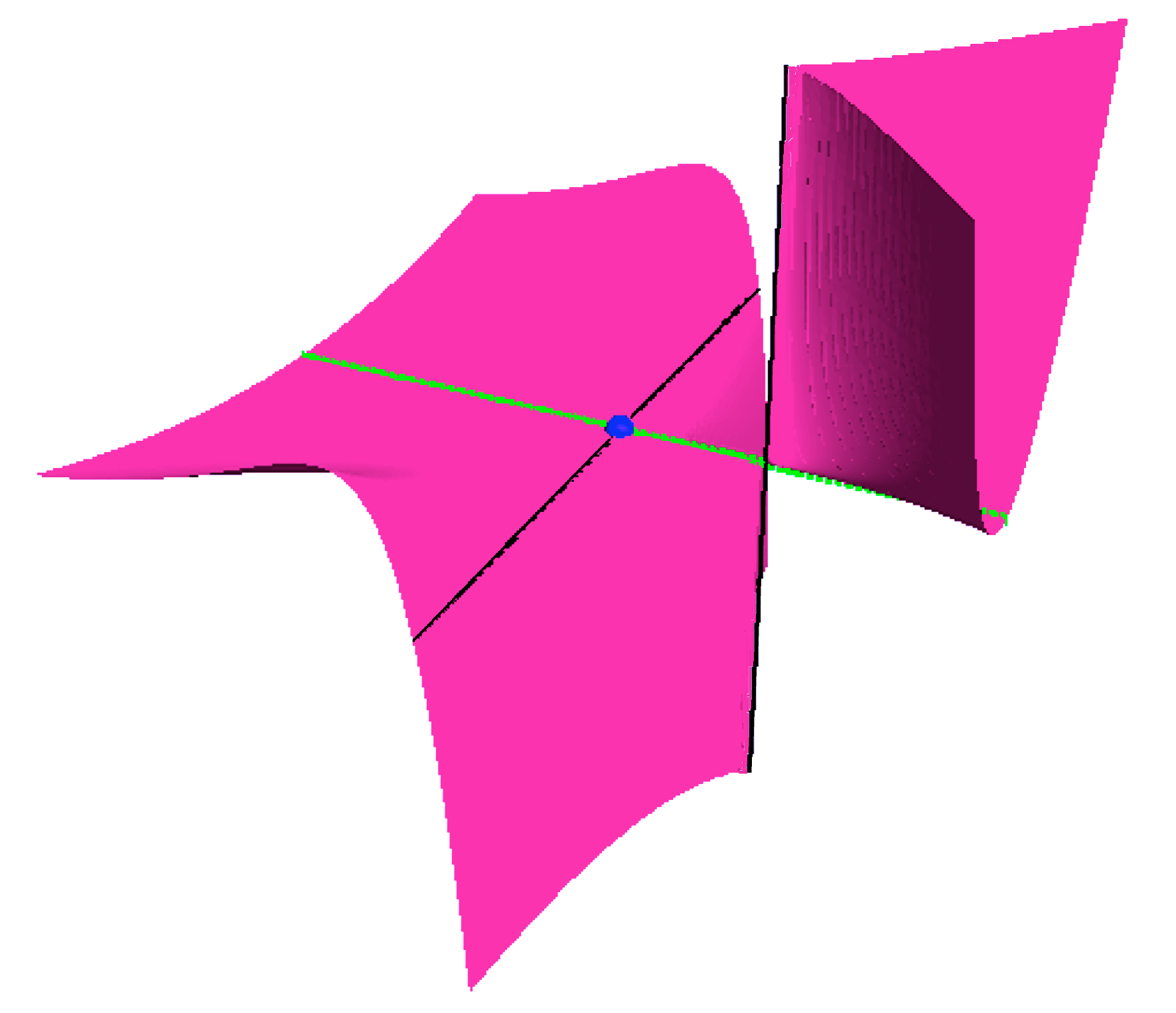}\pdfrefximage\pdflastximage\hfil} 
\vspace{-4mm}\caption{
$X$ in the affine charts  $\{s_1\ne 0, x_0\ne 0\}$
and $\{s_0\ne 0, x_0\ne 0\}$
}\label{Figure 3}
\end{figure}

\setlength{\extrarowheight}{2mm}

\begin{table}[!ht]
\begin{tabular}{|c|c|c|cr|} \hline
&\raisebox{-0.4ex}[0ex][0ex]{orbit} &&& \\
\raisebox{1.5ex}[0cm][0ex]{$v\in\P^1$} & \raisebox{1.5ex}[0cm][0ex]{dim.} & \raisebox{1.5ex}[0cm][0ex]{$F$} & \raisebox{1.5ex}[0cm][0ex]{$X_{v,F}$} & \\ \hline \hline 
& & $(0,0)$ &$ (1:0:0)$ &\\
& 0&$ (1,0)$ &$ (0:1:0)$ &\\
\raisebox{0ex}[0cm][1.5ex]{0} & & $(2,-1)$ &$ (0:0:1)$ &\\ \cline{2-5}
& & $\Conv((0,0), (1,0))$ & $\{(1:t:0):t\in\T^1\} $ & $\subset X^\I$ \\
&\raisebox{2ex}[0cm][1.5ex]{1}&$\Conv((1,0), (2,-1))$ & $\{(0:1:t):t\in\T^1\}$ &\\ \hline
& &$ (0,-1)$ & $(1:0:0)$ &\\ 
1&\raisebox{2ex}[0cm][1.5ex]{0} & $(2,0)$ & $(0:0:1)$ & \llap{$\in\Im(\varphi)$} \\ \cline{2-5}
&\raisebox{0ex}[0cm][1.5ex]{1} & $\Conv((0,-1), (2,0))$ &$ \{(1:0:t^2):t\in\T^1\}$ &\\ \hline
& & $(0,1)$ &$ (1:0:0)$ &\\
& 0& $(1,2)$ & $(0:1:0)$ &\llap{$\in\Im(\varphi)$}\\
\raisebox{0ex}[0cm][1.5ex]{$\infty$}& &$ (2,1)$ & $(0:0:1)$ &\\ \cline{2-5}
& & $\Conv((0,1), (1,2))$ & $\{(1:t:0):t\in\T^1\}$ &\\
&\raisebox{2ex}[0cm][1.5ex]{1}&$\Conv((1,2), (2,1))$ & $\{(0:1:t):t\in\T^1\}$ &\\ \hline
& & $(0,0)$ & $(1:0:0)$ &\\
{$v\ne0,1,\infty$}&\raisebox{1.5ex}[0cm][1.5ex]{0}  & $(2,0)$ & $(0:0:1)$ &\\ \cline{2-5}
&\raisebox{0ex}[0cm][2ex]{1} & $\Conv((0,0), (2,0))$ & $\{(v-1:(v-1)^2t:vt^2):t\in\T^1\}$ & 
\hspace*{-10mm}$ \subset X^0$ \\ \hline
\end{tabular}

\vspace{3mm}\caption{}\label{Table 1}
\bigskip\bigskip

\begin{tabular}{|c|ccc|c|}  \hline
&&&& \raisebox{0ex}[0ex][0ex]{initial systems giving} \\
\raisebox{2ex}[0cm][0ex]
{$v\in\P^1$} & \raisebox{2ex}[0cm][0ex]{$\lambda_{g_v}(s-1)$} &  \raisebox{2ex}[0cm][0ex]{$\lambda_{g_v}(s-1)^2$} &  \raisebox{2ex}[0cm][0ex]{$\lambda_{g_v}(s)$} & 
\raisebox{0.75ex}[0cm][0ex]{equality conditions} \\ \hline \hline 
\raisebox{0ex}[0cm][1.5ex]{0} & $-1$ & 1& 1 & $\cE_0 \cup \cE_1 \cup\cE_2\cup\cE_{1,2}$\\ 
\hline
\raisebox{0ex}[0cm][1.5ex]{1} & 1 & 1& 1 & $\cE_0\cup\cE_2\cup\cE_{0,2}$ \\ \hline
\raisebox{0ex}[0cm][1.5ex]{$\infty $} & 1 & 1 & 1 &  \hspace{2mm} 
$\cE_0\cup\cE_1\cup\cE_2\cup\cE_{0,1}\cup\cE_{1,2}$ \hspace{2mm} \\
\hline
\raisebox{0ex}[0cm][1.5ex]{$v\ne 0,1,\infty $} & $v-1$ & $(v-1)^2$& $v$ &  $\cE_0\cup\cE_2$  \\
\hline
\end{tabular}
\vspace{3mm}\caption{}\label{Table 2}
\vspace{-8mm}
\end{table}
\setlength{\extrarowheight}{0mm}

\medskip
The torus action on $X$ is 
$$
*:\T^1\times X \to X\quad,\quad (t,(\bfs,\bfx))
\mapsto\big( \bfs, (x_0:tx_1:t^2x_2)\big)\enspace.
$$
According to proposition~\ref{X_v}, the decomposition 
of each fiber  $X_v$ 
into orbits is in bijection with the slopes
of the $v$-adic polytope $Q_v\subset \R^2$ shown in figure~\ref{Figure 1}. 
This correspondence is described as follows: 
set 
$$
\alpha_0(s):=s-1 \quad ,\quad  \alpha_1(s):=(s-1)^2 \quad ,\quad  \alpha_2(s):=s 
$$ 
and 
for each $v\in \P^1$ consider a local parameterization $g_v(z)=z+v$ if $v\ne\infty$ and $g_v(z)=z^{-1}$ if $v=\infty$. 
For  $F\in\Slope(Q_v)$ consider then the point $\alpha_{v,F}\in\P^2$ defined 
by 
 $(\alpha_{v,F})_j=\lambda_{g_v}(\alpha_j)$ if $(j,-\ord_v(\alpha_j))\in F$ and 0 otherwise ({\it see} \S~\ref{parameterizations} for the notation); the initial coefficients $\lambda_{g_v}(\alpha_j)$ are explicitly given in the second column of table~\ref{Table 2}.  
The correspondence between slopes of $Q_v$ and orbits of $X_v$ is 
$$
F\mapsto  X_{v,F}:=\T^1*(v,\alpha_{v,F}) \enspace.
$$

Table~\ref{Table 1} describes the orbit decomposition of $X$, which follows readily by considering the $v$-adic polytopes in figure~\ref{Figure 1}. As in this figure, the first column distinguishes the different types of points. The second column gives for each type the possible orbit dimension, which coincides with the dimension of the different slopes $F$ listed in the third column as the convex hull of 1 or 2 points of $\R^2$ (one slope {\it per} line). The fourth column lists the orbits $X_{v,F}$ or more precisely, their projections to the second factor $\P^2$ of $\P^1\times\P^2$. 

According to proposition~\ref{combeferre}, for each $v$ there is exactly  one orbit
contained in $\Im(\varphi)$, which is the orbit  corresponding to the only ``horizontal'' slope $F$ of $Q_v$; note that such a slope might be 0-dimensional. In the right side of the fourth column of table~\ref{Table 1} we indicate the smallest
equivariant subset from the chain  $X^0\subset X^\I\subset X^\Ef\subset \Im(\varphi)$ containing this particular orbit; all other orbits lie in the complement $X\setminus\Im(\varphi)$.

In this example, the fiber $X_v$ is a parabola for  $v\ne0,1,\infty$,  while for $v=0$ and $v=\infty$ it consists in a couple of lines and for $v= 1$ it is a double line, that can be identified in figure~\ref{Figure 3}. It is interesting to note that  the image of $\varphi$ reduces over $0$ to one of these lines, while over $1$ and $\infty$ it collapses into the points~$(0:0:1)$ and $(0:1:0)$ respectively.  

\medskip
Consider the general  polynomials associated to the system~(\ref{a system}), 
or equivalently to the roof functions in figure~\ref{Figure 1}:
\begin{equation}\label{a general system}
f_i=f_{i,0} (s-1) + f_{i,1}(s-1)^2 t + f_{i,2}st^2 \quad \mbox{ for } i=0,1 \enspace.
\end{equation}
For $j=0,1,2$ consider the system of equations 
$$
\cE_j: f_{0,j}=f_{1,j}=0 
$$
and similarly for $j,k=0,1,2$ set
$$
\cE_{j,k}: f_{0,j}t^j+f_{0,k}t^k = f_{1,j}t^j+f_{1,k}t^k = 0 \quad \mbox{ for } t\in \T^1 \enspace.
$$
With this notation, table~\ref{Table 2} lists in its third column the sufficient conditions in order to have equality in our estimate, that follow from proposition~\ref{equality in mainthm}. The initial systems listed correspond to the orbits in $X\setminus X^\Ef$, and the equality condition is expressed as the non solvability of all of these  initial systems. 

We remark that the non solvability of the  three systems $\cE_{j,k}$, $\cE_j$ and $\cE_k$ is equivalent to the single inequation $\det \Big({f_{0,j} \, f_{0,k} \atop f_{1,j} \, f_{1,k}}\Big)\ne 0$. From this, the whole of the equality conditions from table~\ref{Table 2} reduces to the non vanishing of the three determinants as above for $(j,k)=(1,2)$, $(0,2)$ and $(0,1)$, that is: if
$$
\det\begin{pmatrix}
f_{0,0} & f_{0,1}\\
f_{1,0} & f_{1,1} 
\end{pmatrix} 
\det\begin{pmatrix}
f_{0,1} & f_{0,2}\\
f_{1,1} & f_{1,2} 
\end{pmatrix}
\det\begin{pmatrix}
f_{0,0} & f_{0,2}\\
f_{1,0} & f_{1,2} 
\end{pmatrix} \ne 0 
$$
the system~(\ref{a general system}) has exactly $3$ roots in $\C\times \C^\times$.

\smallskip
On the other hand, Ku\v snirenko-Bern\v stein's genericity conditions amount to the non solvability of all of the initial systems corresponding to the faces of the pentagon in figure~\ref{Figure 2}. In the present example, these conditions fail because the initial system $f_{0,0} (s-1)=f_{1,0}(s-1)=0$ corresponding to the face $\Conv((0,0),(0,1))$ admits the root $(s,t)=(1,1)\in\T^2$. This explains why the system~(\ref{a system}) is generic with respect to our estimate but not with respect to Ku\v snirenko-Bern\v stein's one.  

\medskip
\subsection{Other examples}

\begin{example}\label{exmpl 2}{\rm
Consider a plane curve $S\subset \P^2$ of degree $D$, not necessarily smooth. For $j=0,1,2$ let $\ell_j \in \K[x_0,x_1,x_2]$ be a linear form that defines a  line $H_j\subset \P^2$ intersecting $S$ transversely in $D$ points ({\it a fortiori} smooth) and such that the obtained intersections are pairwise disjoint, namely $S\cap H_j\cap H_m= \emptyset$ for $j\ne m$. Let $k\ge 1$, we identify $\ell_j$ with the corresponding section of the universal line bundle $O(1)$ and we consider 
polynomials
$$
f_i= f_{i,0} \ell_1^{2k}+f_{i,1}\ell_0^k\ell_1^k t+f_{i,2} \ell_0^{2k-1}\ell_2 t^{2} \in \Gamma(O(2k))[t] \quad , \quad \mbox{ for } i=0,1 
$$
for some $f_{i,j}\in \K$. We compute the estimate in theorem~\ref{gengenthm} in this situation. This is an unmixed system and so the $v$-adic polytope of $f_0$ and $f_1$ coincide for every place $v\in V_S$. For $v\in H_j\cap S$ we have $\ord_v(\ell_j)= 1$ because of the transversality assumption. We explicit in figure~\ref{Figure 4} below the corresponding family of $v$-adic polytopes and roof functions, for $k=2$.  

\begin{figure}[htbp] 
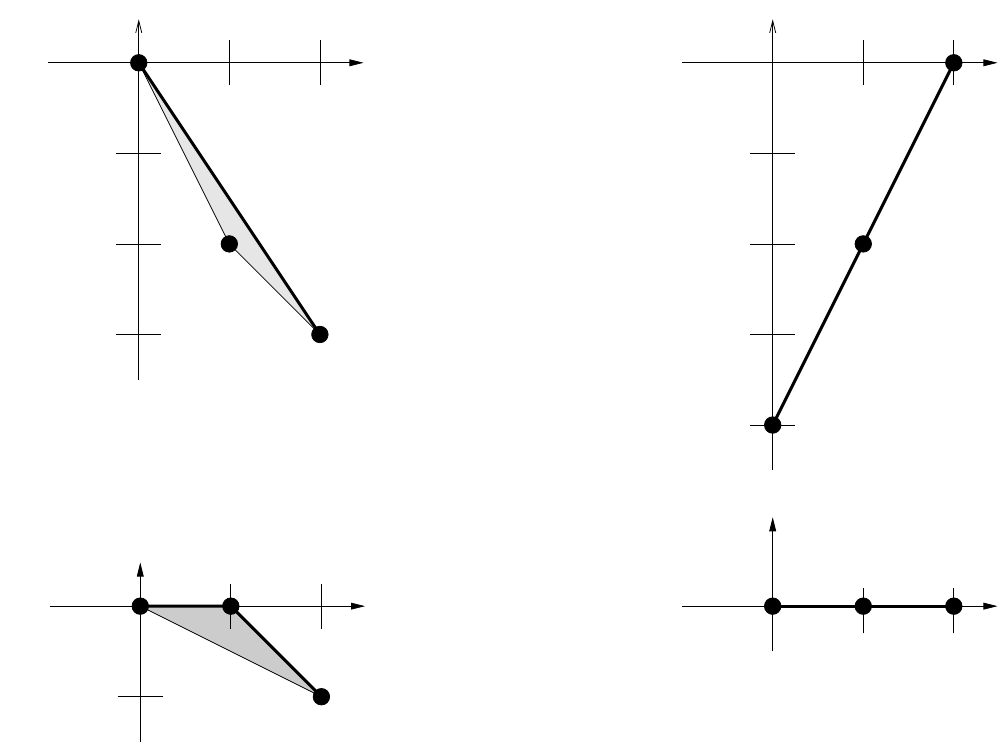
\vspace{-5mm}\caption{}\label{Figure 4}
\end{figure}

The polynomials $f_0$ and $f_1$ have no base points and so $\ov{\vartheta}_{i,v}=0$ for all $i,v$. The resulting estimate  for the number of roots  in $S\times \T^1$ 
of the system $f_0=f_1=0$ is
$$
2!\, \Big(
\deg(O(2k)|_S)\Vol_1([0,2])+ 
\sum_{v\in H_0\cap S} \int_0^{2} \vartheta_v\, \d u + \sum_{v\in H_1\cap S} \int_0^{2} \vartheta_v\, \d u 
+ \sum_{v\in H_2\cap S} \int_0^{2} \vartheta_v\, \d u \Big)
$$
which gives $8kD +(-4k+2)D -4kD -D= D$. This can be verified by solving explicitely the system of equations $f_0=f_1=0$, which for generic $f_{i,j}$'s is equivalent to a system $ \nu_0\ell_0+\nu_2\ell_2=0, t=\nu_1(\ell_1/\ell_0)^k$ over $S\times \T^1$, 
for some $\nu_j\in\K$.

\smallskip
The example in the introduction is the case $S=\{\bfx \in\P^2 : x_2=0\}$ 
with linear forms  $\ell_0=x_0, \ell_1=x_1-x_0, \ell_2=x_1$ in the above construction, so that $H_0\cap S= \{\infty\}$, $H_1\cap S= \{1\}$ and $H_0\cap S= \{0\}$. 
}\end{example}

\smallskip

\begin{example}\label{exempl 4}{\rm 
Consider a semi-abelian surface $G$, extension $0\to \T\to G\to E\to 0$ of an elliptic curve $E$ by a $1$-dimensional 
torus. In what follows we assume the reader is familiar with the material in~\cite[pp.~66-67 and~pp.191-193]{Wal87}. 

This extension corresponds to a point $u_0\in E$, and the open subset of $G$ over $E\setminus\{0,u_0\}$ 
identifies  with $(E\setminus\{0,u_0\})\times\T$. 
The algebra of regular functions of this open subset gives an embedding of $G$ into $\P^2\times\P^2$ 
in the following way: 
let  $\omega_1$, $\omega_2$ and $\eta_1$, $\eta_2$ denote the periods and quasi-periods of the Weierstrass 
function $\wp$ associated to $E$ and set 
$\Lambda:=\Z(\omega_1,\eta_1u_0)+\Z(\omega_2,\eta_2u_0)+\Z(0,2{\rm i}\pi)\subset\C^2$, so that $G\simeq \C^2/\Lambda$. 
Set $F(z,y)=\frac{\sigma(z-u_0)}{\sigma(z)}{\rm e}^y$ where  $\sigma$ denotes the Weierstrass sigma function. 
The embedding is 
$$\begin{array}{rcl}
\C^2/\Lambda &\longrightarrow &\P^2\times\P^2\\[2mm]
(z,y) &\longmapsto &\left(\left(1:\wp(z):\wp'(z)\right), \left(\frac{\wp'(z)+\wp'(u_0)}{\wp(z)-\wp(u_0)}:F(z,y):F(z,y)^{-1}\right)\right)
\enspace.\end{array}$$
The polynomials occurring in this situation have the form 
$$f(z,y) = \sum_{j=0}^N A_j\left(\wp(z),\wp'(z), \frac{\wp'(z)+\wp'(u_0)}{\wp(z)-\wp(u_0)}\right) \, F^{a_j}$$
for some $a_j\in \Z$ and $A_j$ three-variate polynomials. 
Considering the $A_j$'s as sections of a line bundle  $O(d_0 [0]+d_1 [u_0])$ over $E$ for some $d_0,d_1\ge 0$, we 
can apply theorem~\ref{genthm} for bounding the number of common zeros in $(E\setminus\{0,u_0\})\times\T$
of two such polynomials. 
The only positive contributions to the estimate come from the places $0$ and $u_0$ of $E$, 
whereas the contribution at all other places is negative or zero. 

In particular, for polynomials of the simpler form $f(\wp,\wp',F)$ which correspond to sections of the line bundle
$O(d_0[0])$,  the upper bound does not depend on the extension itself (namely on $u_0$). 
Even more particularly, consider integers $d$, $D$, $L$ and two polynomials $f_1$ and $f_2$ of the special form
$$f_i = \sum_{j=0}^{d} A_{i,j}(\wp)\,(\wp-\wp(u_1))^{(d-j)D}F^{jL}
$$
for some 
$u_1\not=0,\pm u_0$, and polynomials 
$A_{i,j}\in\C[s]$ such that $\deg(A_{i,j})\leq d$ and $(\wp-\wp(u_1))\nmid A_{i,d}(\wp)$ for $i=1,2$ and $0\le j\le d$. 
By considering the corresponding roof functions, we verify that the number of common isolated roots in $(E\setminus\{0,u_0\})\times\T$ 
of $f_1$ and $f_2$ 
is bounded above by $2d^2L$, 
which is most interesting when $d$ is significantly smaller than $D$ and $L$. 
}\end{example}

\smallskip

\begin{example}{\rm 
The examples presented so far
can all be easily 
handled by hand. 
However, this is not so for a typical system of equations as the required calculations can be too bulky; in those cases 
our results can help 
to determine the number of solutions. 
For instance, theorem~\ref{mainthm} together with proposition~\ref{equality in mainthm}
shows that 
the system
$$\displaylines{
(s+1)^2 + (s^2-1)t + (s^2-1)t^2 + (s-1)^2t^3 + (s-1)(s+2)t^4 + (s-1)(s+2)t^5\cr
2(s+1)^2 + (s^2-1)t - (s^2-1)t^2 + 3(s-1)^2t^3 - 4(s-1)(s+2)t^4 - 2(s-1)(s+2)t^5
}$$
has exactly $8$ solutions in $\A^1\times\T^1$, 
that can be calculated with a 
computer algebra software. 
In comparison, both the bihomogeneous B\'ezout theorem
and Ku\v snirenko-Bern\v stein's one allow up to $20$ solutions.
}\end{example}

\medskip\goodbreak

\subsection{Practical considerations}\label{Practical considerations}

\subsubsection{Computing the bound}\label{Computing the bound}
Let $f_i=\sum_{j=0}^{N_i} \alpha_{i,j}(s) \bft^{a_{i,j}}\in \K[s][\bft^{\pm1}]$ for 
$0\le i\le n$ be a family of primitive Laurent polynomial, as in the statement of theorem~\ref{mainthm}. The estimate~(\ref{mainthm_display}) only depends on the configuration of the roots of the coefficients of the $f_i$'s and {\it not} on their actual value, this configuration can be computed from factorizations 
\begin{equation}\label{facto_alphas}
\alpha_{i,j}(s)= \lambda_{i,j} \prod_{p\in P} p(s)^{e_p(i,j)} \quad \mbox{ for } 0\le i\le n \mbox{ and } 0\le j\le N_i
\end{equation}
for some finite set $P\subset \K[s]$ of pairwise coprime polynomials, $e_p(i,j)\in \N$ and $\lambda_{i,j}\in \K^\times$, 
in the following way. 

\begin{prop} \label{mainthm_calculable}
With notation as above, for $p\in P$ let $\rho_{i,p}:\NP(f_i)\to \R$ be the parameterization of the upper envelope of the polytope
$$
\Conv\big((a_{i,0},-e_p(i,0)), \dots, (a_{i,N_i},-e_p(i,N_i)\big) \subset \R^{n+1} \enspace.
$$
For $v\in \P^1$ let $\vartheta_{i,v}$ denote the roof function of $f_i$ at $v$, then  
if $v\ne\infty$ is a zero of some $p\in P$ with multiplicity $m_v\ge 1$ we have 
$\vartheta_{i,v} = m_v\rho_{i,p}$, otherwise $\vartheta_{i,v}=0$. 
Furthermore
\begin{equation} \label{mainthm_calculable_display}
\sum_{v\in \P^1} \MI_n (\vartheta_{0,v},\dots,\vartheta_{n,v})= 
\MI_n (\vartheta_{0,\infty},\dots,\vartheta_{n,\infty})
+ \sum_{p\in P} \deg(p) \MI_n (\rho_{0,p},\dots,\rho_{n,p})\enspace.
\end{equation}
\end{prop}

\begin{proof}
Let $p\in P$, then for $v\in Z(p)$ a zero of multiplicity $m_v$  we have 
$\ord_v(\alpha_{i,j})= m_v e_p(i,j)$  and so $\vartheta_{i,v} = m_v\rho_{i,p}$.
By proposition~\ref{transfolineaire}, for each such $v$  we have
$\MI_n(\bfvartheta_{v})=m_v \MI_n(\bfrho_{p})$ and
so 
$$
\sum_{v\in Z(p)} \MI_n(\bfvartheta_{v})= \sum_{v\in Z(p)} m_v\MI_n(\bfrho_{p})
= \deg(p) \MI_n(\bfrho_{p}) \enspace. 
$$
The identity~(\ref{mainthm_calculable_display}) follows by summing up over $p\in P$. 
\end{proof}

The  factorizations~(\ref{facto_alphas}) can be computed through the algorithm in the proof of lemma~\ref{factorization} below, with no need for extracting the roots of the $\alpha_{i,j}$'s. On the other hand, the roof functions $\vartheta_{i,\infty}$ can be computed from knowledge of the degree of the $\alpha_{i,j}$'s. Hence the estimate in theorem~\ref{mainthm} can be computed with operations in the field of definition of the $f_i$'s.

\begin{lem}\label{factorization}
Let  $G \subset \K[s]\setminus \{0\}$ be a finite set, we can compute a finite set $P\subset \K[s]$ of pairwise coprime polynomials and non negative integers $(e_p(g): g\in G, p\in P) $ such that $g= \gamma(g)\prod_{p\in P} p^{e_p(g)}$ for all $g\in G$ and some $\lambda(g)\in \K^\times$, with the operations \enspace $\cdot$\enspace,\enspace $/$ \enspace and \enspace $\gcd$\enspace over $\K[s]$ only. 
\end{lem}

\begin{proof}
Set $P=\emptyset$ and let $G=\{g_1,\dots,g_r\}$ be the given family of polynomials; we can suppose that not all of these polynomials are constant, otherwise we are done. Set $p_1$ for the {\it last} $g_j$ which is not constant, then compute the biggest power  $c_{1}\ge 0$ such that $p_1^{c_{1}}|g_1$. This can be done with division and $\gcd$ computations. If
$$
\gcd\Big(p_1 , \frac{g_1}{p_1^{c_{1}}}\Big) =1
$$
we do similarly with $g_2$ and so on. On the contrary, if this gcd is not constant we set $p_2:= \gcd(p_1 , {g_1}/{p_1^{c_{1}}})$ and we start over from $g_1$. Eventually, this procedure ends because $\deg(p_{k-1})>\deg(p_{k}) \ge 1$ and when it does, we have obtained partial factorizations
$$
g_j= p_k^{c_j} h_j \quad \mbox{ for } 1\le j\le r
$$ 
for some $h_j\in \K[s]$ such that $ \gcd(p_k, h_j)=1 $ for all $j$. We add the non constant polynomial $p:=p_k$ to $P$ and we set $e_p(g_j):=c_j$ for $1\le j\le r$, then we reapply the algorithm to the family $\{h_1,\dots, h_r\}$ instead of $G$, repeating the procedure until the factorization is completed. 
\end{proof}

\subsubsection{Comparing theorem~\ref{mainthm} and theorem~\ref{genthm}}
\label{Comparing}
Let $f_i \in \K[s][\bft^{\pm1}]$  for $0\le i\le n$ be a family of primitive Laurent polynomials, and  suppose that $f_0$ is reduced and irreducible and depends only in the two variables $s,t_1$. The number of isolated roots of such a system in $\A^1\times\T^n$ can be estimated through the following two options
\begin{enumerate}
\item directly with theorem~\ref{mainthm}; 
\item by applying theorem~\ref{genthm} (or  theorem~\ref{gengenthm} in the case of a singular curve) to the system $f_1,\dots, f_n$ with respect to the curve $S=\ov{Z(f_0)}\subset\P^2$.
\end{enumerate}
Which one is best? For $0\le i\le n$ write
$$
f_i = \sum_{k=0}^{M_i} \beta_{i,k} (s,t_1){(t_2,\dots,t_n)}^{b_{i,k}}
$$
for some $b_{i,k}\in \Z^{n-1}$ and $\beta_{i,k}\in \K[s, t_1^{\pm1}]$. Option (2) corresponds to this expression for the $f_i$'s, after interpreting the coefficients $\beta_{i,k}$ as sections $\sigma_{i,k}$ of a line bundle $O(\delta_i)$. 

Next expand each $\beta_{i,k}$ as $\beta_{i,k}=\sum_{j=0}^{N_{i,k}} \alpha_{i,k,j}(s) t_1^{a_{i,k,j}}$ for some $a_{i,k,j}\in \Z$ and $\alpha_{i,k,j} \in \K[s]$, so that 
$$
f_i=\sum_{k=0}^{M_i}\sum_{j=0}^{N_{i,k}} \alpha_{i,k,j}(s) \bft^{(a_{i,k,j},b_{i,k})}\enspace. 
$$
Option (1) corresponds to this expansion for the $f_i$'s; note that this places the system in a more generic situation than the first option. 

\medskip
Whenever the family of sections $(\sigma_{i,k}:0\le k\le M_i)$ has no base point for $0\le i\le n$, the estimate obtained from theorem~\ref{genthm} is generically 
attained in $Z(f_0) \times \T^{n-1}\subset \A^1\times\T^n$. In that case, option~(2) is  preferable to option~(1): the obtained estimate will not be worse since the system is put in a less generic situation, furthermore it is easier to compute since it involves mixed integrals of lower dimension. 

\medskip 
The following variant of example~\ref{exmpl 2} illustrates the above discussion. 

\begin{example} \label{exmpl 3}{\rm 
Consider the polynomials in $\K[s_1,s_2,t]$ 
\begin{align}\label{another system}
f_0&= s_2^{2k}-(s_1-1)^{2k} -1 \enspace, \nonumber \\
f_1&= (s_2^{2k} -1) + (s_1-1)^kt-s_1t^2 \enspace,
 \\ 
f_2&= (3-3s_2^{2k}) + (s_1-1)^kt-s_1t^2 \enspace.\nonumber
\end{align}

Firstly, we estimate the number $R$ of isolated roots of this system in $\A^1\times \T^2$ 
by applying theorem~\ref{genthm} to $f_1=f_2=0$ as a system 
over the 
smooth complete curve $S:=\ov{Z(f_0)}\subset\P^2$.
We have 
\begin{align*}
f_1 &\equiv f_1-f_0= (s_1-1)^{2k}  + (s_1-1)^kt-s_1t^2   &&\pmod{f_0}\\
f_2 &\equiv f_2+3f_0 =-3 (s_1-1)^{2k}  + (s_1-1)^kt-s_1t^2 &&\pmod{f_0}
\end{align*}
and so the considered system reduces to the one in example~\ref{exmpl 2} for the curve $S$ 
and linear forms $\ell_0=x_0, \ell_1=x_1- x_0, \ell_2=x_2$. These linear forms satisfy the required conditions of transversality and disjointness for the intersections, and the calculations in example~\ref{exmpl 2} show that the number of isolated roots of the system $f_1-f_0=f_2+3f_0=0$ 
over $S$ is bounded above by $\deg(S)=2k$.  
We can verify that this estimate $R\le 2k$ is exact by explicitely  solving the system, 
which turns to be equivalent to $s_1=2,s_2^{2k}=2, t=1$. 

\medskip
On the other hand, theorem~\ref{mainthm} applied to~(\ref{another system}) 
as a system in $\K[s_1][s_2^{\pm1},t^{\pm1}]$ bounds $R$ from above by the number of roots of the associated {\it generic} system 
\begin{align*}
F_0 &= F_{0,0} s_2^{2k} + F_{0,1}((s_1-1)^{2k} -1) \enspace, \\
F_i &= F_{i,0}s_2^{2k} +F_{i,1} + F_{i,2}(s_1-1)^kt+F_{i,3}s_1t^2  \quad (i=1,2)
\enspace.
\end{align*}
This system is equivalent to the generic system
\begin{align*}
F_0 &= F_{0,0} s_2^{2k} + F_{0,1}((s_1-1)^{2k} -1) \enspace, \\
G_i &= G_{i,0}(s_1-1)^{2k} +G_{i,1} + G_{i,2}(s_1-1)^kt+G_{i,3}s_1t^2 \in\K[s_1,t]  \quad (i=1,2)
\enspace.
\end{align*}
The number of roots of $G_1=G_2=0$ in $\A^1\times\T$ can be computed
by applying corollary~\ref{cor_mainthm}. 
The only non-zero contributions to the adelic formula~(\ref{cor_mainthm_display}) 
come from the places $0$ and $\infty$. The corresponding roof functions are those in the left hand side of figure~\ref{Figure 4}, therefore the 
number of roots of this system is $4k+1$. 
For each such root $(s_1,t)$ we obtain $2k$ values of $s_2$ by solving $F_0=0$, which shows that the number of roots of $F_0=F_1=F_2=0$ equals $2k(4k+1)=8k^2+2k$. 

This gives the estimate $R\leq 8k^2+2k$, which is much worse than the exact estimate $R=2k$ 
obtained from theorem~\ref{genthm}. 
}\end{example}

\bigskip

\section{Basic properties of the mixed integral}\label{bk_multiint}

In~\cite[\S~IV]{PS03} we introduced the mixed integral of a family of concave functions. 
In what follows we summarize its basic properties and pursue its study, 
in particular by establishing 
a decomposition formula (proposition~\ref{marius} below) expressing the mixed integral in terms of lower dimensional mixed integrals and volumes. 

\medskip 
By definition, a {\it convex body} of $\R^n$ is a non-empty, convex and compact subset. 
The {\it mixed volume} of a family of convex bodies $Q_1, \dots, Q_n$ of $\R^n$ is defined as
\begin{equation} \label{defmultivolume}
\MV_n(Q_1, \dots, Q_n) := \sum_{j=1}^n (-1)^{n - j} \sum_{1 \le i_1 < \cdots < i_j \le n} \Vol_n(Q_{i_1} + \cdots + Q_{i_j}) 
\end{equation} 
where $\Vol_n$ denotes the $n$-dimensional Hausdorff (or Lebesgue) measure of $\R^n$. 
This generalizes the volume of a convex body, since $\MV_n(Q, \dots, Q) = n! \, \Vol_n(Q)$. 
The mixed volume is symmetric and linear in each variable $Q_i$ with respect to the Min\-kows\-ki sum, 
and monotone with respect to inclusion~\cite[chap.~IV]{Ewa96}, \cite[chap.~5]{Schn93}.

\smallskip
In what follows, all concave functions are supposed to be defined on convex bodies. 
The mixed integral (definition~\ref{defmultiintegrale}) is the natural extension to concave functions of the mixed volume of convex bodies. 
It is symmetric and linear in each variable $\rho_i$ with respect to the sup-convolution~$\boxplus$ and for a function $\rho:Q\to\R$ we have
$ \MI_n(\rho, \dots, \rho) = (n+1)! \, \int_Q \rho(u) \, \d \Vol_n(u) $. 

It is possible to express the mixed integral in terms of mixed volumes: for a concave function $\rho: Q \to \R$ and a constant $\gamma \le \min(\rho,0) $ consider the 
polytope $Q_{\rho,\gamma} \subset \R^{n+1}$,  defined as the convex hull
$$Q_{\rho,\gamma} := \Conv\big( \Graph(\rho), Q \times \{\gamma\}\big) = \Conv\big( (u, \rho(u)), (u, \gamma) \, : \ u \in Q \big) \enspace.$$
Note that $\int_Q\rho(u)\, \d \Vol_n(u) = \Vol_{n+1}(Q_{\rho,\gamma}) + \gamma\Vol_n(Q)$. 
Then for $\gamma_i \le \min(\rho_i,0)$ we have~\cite[prop.~IV.5(d)]{PS03}
\begin{align} \label{masure_gorbeau}
\MI_n(\rho_0, \dots, \rho_n) =& \MV_{n+1} (Q_{\rho_0, \gamma_0}, \dots, Q_{\rho_n, \gamma_n}) \\[0mm]
\nonumber &+ \sum_{i=0}^n \gamma_i \, \MV_n(Q_0, \dots, Q_{i-1}, Q_{i+1}, \dots, Q_n) \enspace.
\end{align}
This identity
together with the monotonicity of the mixed volume readily  
implies that the mixed integral is monotone too: 

\begin{prop} \label{fantine}
For $0\le i\le n$ let $\rho_i$ and $\sigma_i$ be concave functions defined over the same convex body $Q_i$
and such that $\rho_i \ge \sigma_i $, then $\MI_n(\rho_0,\dots, \rho_n)\ge \MI_n(\sigma_0,\dots, \sigma_n)$. 
\end{prop}

In particular $\MI_n(\rho_0,\dots, \rho_n)\ge 0$ whenever the $\rho_i$'s are non-negative. The mixed integral behaves well with respect to 
linear changes of variables: 

\begin{prop}\label{transfolineaire}
Let $\ell:\R^n\rightarrow\R^n$ be an invertible linear transformation and 
for $0\le i\le n$ let $\rho_i$ be a concave function 
defined over a convex body  of $\R^n$, 
then
$$
\MI_n(\rho_0\circ\ell,\dots,\rho_n\circ\ell) = |\det(\ell)|^{-1}\MI_n(\rho_0,\dots,\rho_n)
\enspace.$$
\end{prop}
\begin{proof}
By the very definition of the mixed integral the formula reduces to the same one for integrals, where it is just the formula for a linear change of variables.
\end{proof}

\begin{prop} \label{bagne}
For $0\le i\le n$ let $R_i\subset \R^{n+1}$ be a convex body sitting above $Q_i \subset \R^n$ {\it via} the projection 
$\pi:\R^{n+1}\to \R^n$ which forgets the last coordinate, and set $u(R_i), \ell(R_i) :Q_i\to \R$ for the parameterization of the upper and lower envelope of $R_i$, respectively. Then 
$$
\MI_n\big(u(R_0), \dots, u(R_n)\big) + \MI_n\big(-\ell(R_0), \dots,-\ell(R_n) \big) = \MV_{n+1}(R_0, \dots,R_n) \enspace.$$
\end{prop}

\begin{proof}
For convex bodies $R,S\subset \R^{n+1}$  
we have that $u(R+S)= u(R)\boxplus u(S)$ and $-\ell(R+S) = (-\ell(R))\boxplus(-\ell(S))$. 
This remark together with the definitions
of the mixed integral and volume
allows to deduce the equality from the (trivial)  
unmixed case $R_0=\cdots=R_n$. 
\end{proof}

Let $\delta_0, \dots,\delta_n \in \R$, as a further consequence of the identity~(\ref{masure_gorbeau}) 
applied separately to the $\rho_i$'s and to the 
$\rho_i+\delta_i$'s we obtain a useful relationship 
between their mixed integrals:
\begin{align}\label{multiintegralefonctionstranslates}
\MI_n(\rho_0+\delta_0,\dots,\rho_n+\delta_n) =&
 \MI_n(\rho_0,\dots,\rho_n) \\[0mm]
&+ \sum_{i=0}^n\delta_i\MV_n(Q_0,\dots,Q_{i-1},Q_{i+1},\dots,Q_n) \enspace. \nonumber
\end{align}

\begin{exmpl} The mixed volume of a parallelepiped is equal to 
the permanent of the matrix of the lengths of the edges of the given parallelepiped times the volume of the similar parallelepiped with edges of unit length. The mixed integral of constants functions on such parallelepipeds can be expressed by an analogous formula:

Let $\ell_1,\dots, \ell_n$ be linear forms of $\R^n$ and for $0\le i \le n$ let
$\bfc_i= (c_{i,1},\dots,c_{i,n})\in\R^n$. For each $i$  consider the parallelepiped $Q(\bfc_i):=\{x\in\R^n:|\ell_j(x)|\leq c_{i,j} \mbox{ for } j=1,\dots,n\}$ and a constant function $\rho_i:Q(\bfc_i)\rightarrow\R$, then 
$$\MI_n(\rho_0,\dots,\rho_n) = \Vol_n(Q(\bfun))\, \Perm
\begin{pmatrix}
\bfc_0 &\cdots &\bfc_n\\
\rho_0 &\cdots &\rho_n
\end{pmatrix}
\enspace.$$
\end{exmpl}

\medskip
The mixed volume of a family of polytopes $Q_1,\dots,Q_n$ can be decomposed in terms 
of the lower dimensional mixed volumes of their faces. 
For a convex body $Q\subset \R^n$ consider its {\it support function} 
$$
h_{Q}:\R^n\to \R \quad , \quad  u\mapsto \max\{\langle u,w\rangle\, : \ w\in Q\}
\enspace,$$
and for $u\in \R^n$ set $Q^u:=\{w\in Q : \langle u,w\rangle= h_Q(u) \}$ for its {\it face in the $u$-direction}. 
Let $\Es^{n-1}$ denote the unit sphere of $\R^n$, then~\cite[chap.~IV, thm.~4.10, p.~126]{Ewa96} 
or~\cite[formula~5.1.22 in~p.~276]{Schn93}
\begin{equation} \label{mme_thenardier}
\MV_n(Q_1,\dots,Q_n)=\sum_{u\in\Es^{n-1}} h_{Q_1}(u)\, \MV_{n-1}(Q_2^u,\dots, Q_n^{u}) \enspace.\end{equation}
This decomposition formula can be extended to general {convex bodies}, 
turning the sum into an integral and replacing the
mixed volume of the faces by the {mixed area measure}. For $w\in \R^n \setminus Q$ set $u(Q,w)\in\Es^{n-1}$ 
for the unit vector pointing from the nearest point in $Q$ towards $w$, and for $\varepsilon>0$ and  $U\subset \Es^{n-1}$ set
$$
B_\varepsilon(Q,U):= \{ w\in\R^n: 0<\dist(Q,w) \le \varepsilon \mbox{ and } u(Q,w)\in U\}\enspace.
$$
For a given convex body $Q\subset \R^n$,  the {\it area measure} $S_{n-1}(Q;\cdot)$ of $\Es^{n-1}$
is defined as 
the limit~\cite[formula~4.2.9 in p.~203]{Schn93}
$$
S_{n-1}(Q; U) =\lim_{\varepsilon\to0} \varepsilon^{-1} \Vol_n(B_\varepsilon(Q,U)) 
\quad \mbox{ for a measurable }U\subset \Es^{n-1} \enspace. 
$$
Then the  {\it mixed area measure} of a family of convex bodies $Q_2,\dots,Q_n$ of $\R^n$ is the measure of $\Es^{n-1}$ defined as~\cite[formula~5.1.20 in p.~276]{Schn93}
$$
S(Q_2,\dots,Q_n;\cdot) := \sum_{k=1}^{n-1}(-1)^{n+k-1}\kern-10pt
\sum_{2\le i_1<\dots<i_k\le n}S_{n-1}(Q_{i_1}+\dots+Q_{i_k};\cdot)
$$
In case the  $Q_i$'s are polytopes, this measure can be expressed as the finite sum of mixed volumes of faces~\cite[formula~5.1.21 in~p.~276]{Schn93}
$$S(Q_2,\dots,Q_n;U) = \sum_{u\in U} \MV_{n-1}(Q_2^u,\dots, Q_n^{u})\enspace.$$
With this notation, the extension of the decomposition~(\ref{mme_thenardier}) to 
general convex bodies $Q_1,\dots,Q_n$ is~\cite[thm.~5.1.6, p.~275]{Schn93}
\begin{equation}\label{l'egout de paris}
\MV_n(Q_1,\dots,Q_n) = \int_{\mathbb{S}^{n-1}}h_{Q_1}(u) \, \d S_{n-1}(Q_2,\dots,Q_n;u)
\enspace.
\end{equation}

\medskip
Let $\rho:Q\to \R$ be a given concave function and consider a continuous extension (not necessarily concave) 
to a neighborhood
of $Q$. It is always  possible to do this, since 
a concave function defined on a convex body is 
continuous. 
In analogy with the area measure, we define the (signed) measure $I_{n-1}(\rho; \cdot)$ on $\Es^{n-1}$
as the limit
$$
I_{n-1}(\rho; U) :=\lim_{\varepsilon\to0} \varepsilon^{-1} \int_{B_\varepsilon(Q,U)} \rho(u) \, \d \Vol_n(u) 
\quad \mbox{ for a measurable }U\subset \Es^{n-1} \enspace. 
$$
For concave functions $\rho_1,\dots,\rho_n$ on convex bodies of $\R^n$, we consider the signed measure 
on $\mathbb{S}^{n-1}$ defined by
$$I_{n-1}(\rho_1,\dots,\rho_n, \cdot):=
\sum_{k=1}^n(-1)^{n+k}\kern-10pt\sum_{1\le i_1<\dots<i_k\le n} I_{n-1}(\rho_{i_1}\boxplus\dots\boxplus\rho_{i_k}; 
\cdot)
\enspace.$$
For piecewise affine $\rho_i$'s defined on polytopes, this measure can be expressed as the finite sum 
$$I_{n-1}(\rho_1,\dots,\rho_n, U) = \sum_{u\in U} \MI_{n-1}(\rho_1|_{Q_1^u},\dots,\rho_n|_{Q_n^u})\enspace.$$
We denote by $\mathbb{S}^{n}_+\subset \R^{n+1}$ the subset of $\mathbb{S}^{n}$ of vectors the last coordinate of which is positive. The following is the analog of~(\ref{l'egout de paris}) for mixed integrals: 

\begin{prop} \label{marius}
Let $\rho_0 : Q_0 \to \R, \dots,  \rho_n: Q_n\to \R$ be a family of concave functions defined on convex bodies. 
Set $Q_{\rho_i}:= \Conv( \Graph(\rho_i)) \subset \R^{n+1}$ for the convex hull of the graph of $\rho_i$ over $Q_i$, then 
\begin{align*}
\MI_n(\rho_0,\dots,\rho_n) = &\int_{\mathbb{S}^{n-1}}
h_{Q_0}(u) \,\d I_{n-1}(\rho_1,\dots,\rho_n;u) 
\\ 
&+ \int_{\mathbb{S}^{n}_+} h_{Q_{\rho_0}}(r) \, \d S_{n-1}(Q_{\rho_1},\dots,Q_{\rho_n};r)
\enspace.
\end{align*}
\end{prop}

For piecewise affine functions this formula takes the finite form:
\begin{align}\label{MI_decomposition}
\MI_n(\rho_0,\dots,\rho_n) =& \sum_{u\in\mathbb{S}^{n-1}} h_{Q_0}(u) \, \MI_{n-1}(\rho_1|_{Q_1^u},\dots,\rho_n|_{Q_n^u}) \\[0mm]
&+\sum_{r\in\mathbb{S}^{n}_+} h_{Q_{\rho_0}}(r) \, \MV_n(Q_{\rho_1}^{r},\dots,Q_{\rho_n}^{r})
\enspace. \nonumber
\end{align}

\begin{proof}
We first prove the proposition for piecewise affine functions defined on polytopes;  
the proof relies on a reduction to mixed volumes.
Take $\gamma_i:=\min(\rho_i,0)$ in the identity~(\ref{masure_gorbeau}), applying the decomposition formula~(\ref{mme_thenardier}) to the resulting mixed volumes we obtain
\begin{equation}\label{dspreuvepropmarius}
\MI_n(\bfrho) = \Phi+\gamma_0\MV_{n}(Q_1,\dots,Q_n)+
\sum_{i=1}^n\gamma_i\Phi_i
\end{equation}
with
$$\Phi = \MV_{n+1}(Q_{\rho_0, \gamma_0}, \dots, Q_{\rho_n, \gamma_n}) = \sum_{r\in\Es^{n}} h_{Q_{\rho_0,\gamma_0}}(r)\, \MV_{n}(Q_{\rho_1,\gamma_1}^r,\dots, Q_{\rho_n,\gamma_n}^{r})$$ 
and, for $1\le i\le n$, 
\begin{align*}
\Phi_i &= \MV_{n}(Q_{0}, \dots,Q_{i-1},Q_{i+1},\dots, Q_{n})\\[1mm]
&= \sum_{u\in\Es^{n-1}} h_{Q_{0}}(u)\, \MV_{n-1}(Q_{1}^u,\dots,Q_{i-1}^u,Q_{i+1}^u,\dots, Q_{n}^{u})
\enspace.\end{align*}
Writing the index variable $r=(r_1,\dots,r_{n+1})\in\Es^n$ we split $\Phi$ into three sums
$\Phi=\Sigma_++\Sigma_0+\Sigma_-$  according to whether $r_{n+1}$ is positive, zero or negative:

\smallskip
{\leftskip=\parindent\parindent=-\parindent
\un{Case $r_{n+1}>0$}~:   we have $Q_{\rho_i,\gamma_i}^r= Q_{\rho_i}^{r}$ for $0\le i\le n$, since $Q_{\rho_i,\gamma_i}$ and $Q_{\rho_i}$ have the same upper envelope, and so 
\begin{equation}\label{mlle guillenormand}
\Sigma_+= \sum_{r\in\Es^{n}_+} h_{Q_{\rho_0}}(r)\, \MV_{n}(Q_{\rho_1}^r,\dots, Q_{\rho_n}^{r})\enspace. 
\end{equation}

\un{Case  $r_{n+1}=0$}~: write  $r=(u,0)$ for some $u\in\Es^{n-1}$.
Then $h_{Q_{\rho_0,\gamma_0}}(u,0) = h_{Q_0}(u)$ and the identity~(\ref{masure_gorbeau}) 
implies that $\MV_{n}(Q_{\rho_1,\gamma_1}^{(u,0)},\dots, Q_{\rho_n,\gamma_n}^{(u,0)})$ is equal to
$$\MI_{n-1}(\rho_1|_{Q_1^u},\dots,\rho_n|_{Q_n^u}) - \sum_{i=1}^n\gamma_i\MV_{n-1}(Q_1^u,\dots,Q_{i-1}^u,Q_{i+1}^u,\dots,Q_n^u)
$$
thus
\begin{equation}\label{soeur_crucifixion}
\Sigma_0 = \sum_{u\in\Es^{n-1}}h_{Q_0}(u)\MI_{n-1}(\rho_1|_{Q_1^u},\dots,\rho_n|_{Q_n^u}) - \sum_{i=1}^n\gamma_i\Phi_i
\enspace.\end{equation}

\un{Case $r_{n+1}<0$}~: we have $\MV_{n}(Q_{\rho_1,\gamma_1}^r,\dots,$ $Q_{\rho_n,\gamma_n}^{r})=0$ for $r\not=({\bf 0},-1)$, 
because each $Q_{\rho_i,\gamma_i}^r$ lies in a translate of the linear 
space $r^\bot \cap (\R^n\times\{0\})$
which for $r\not=({\bf 0},-1)$ has codimension $2$. 
On the other hand, for $r=({\bf 0},-1)$ we have $Q_{\rho_i,\gamma_i}^{r} = Q_i\times\{\gamma_i\}$ and $h_{Q_{\rho_0,\gamma_0}}({\bf 0},-1) =-\gamma_0$, from where follows
\begin{equation}\label{catherina}
\Sigma_- = -\gamma_0 \MV_{n}(Q_{1},\dots, Q_{n})
\enspace.\end{equation}
}

Identities~(\ref{mlle guillenormand}), (\ref{soeur_crucifixion}) and~(\ref{catherina}) together with~(\ref{dspreuvepropmarius}) establish the proposition for the piecewise affine case ({\it i.e.}~(\ref{MI_decomposition})). The general case follows by approximating the $\rho_i$'s  by piecewise affine concave functions and applying the continuity of the mixed integral and of the support functions together with the weak continuity of the mixed area and the $I_{n-1}$ measures. 
\end{proof}

For a single piecewise affine and non-negative function $\rho$, the formula~(\ref{MI_decomposition}) corresponds to the decomposition of the integral into the sum of volumes of pyramids with apex at the point ${\bf 0}_{n+1}$ and base either a wall (for the terms in the first sum) or a face of the roof (for the terms in the second sum) of $Q_{\rho,0}$, as shown in the figure below: 
\vspace{0mm}
\begin{center}
\begin{picture}(0,0)%
\includegraphics{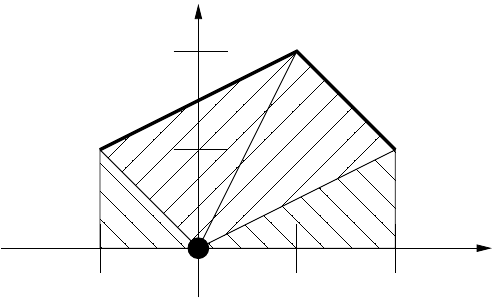}%
\end{picture}%
\setlength{\unitlength}{2072sp}%
\begingroup\makeatletter\ifx\SetFigFont\undefined%
\gdef\SetFigFont#1#2#3#4#5{%
  \reset@font\fontsize{#1}{#2pt}%
  \fontfamily{#3}\fontseries{#4}\fontshape{#5}%
  \selectfont}%
\fi\endgroup%
\begin{picture}(4524,2724)(3139,-5473)
\put(6751,-3661){\makebox(0,0)[lb]{\smash{{\SetFigFont{11}{10.8}{\rmdefault}{\mddefault}{\updefault}{\color[rgb]{0,0,0}$\mbox{Gr}(\rho)$}%
}}}}
\end{picture}%

\end{center}
\vspace{-2mm}

\begin{rem} In the unmixed case $\rho_0=\cdots=\rho_n=\rho$ for some piecewise affine $\rho$
defined by integral conditions, 
the decomposition~(\ref{MI_decomposition}) can be interpreted 
in geometric terms 
as the B\'ezout theorem for Chow weights applied to 
the intersection of a projective toric variety with a monomial divisor, {\it see}~\cite[\S~IV]{PS04} for the details. 
It is possible that the general (integral) 
case of this  decomposition
might be interpreted {\it via} an extension of this result 
to the multiprojective setting. 
\end{rem}

As illustration, consider the functions $\rho:[0,3]\to \R$,  $\sigma:[0,2]\to\R$ 
with  graph given by the figure below

\vspace{1mm}
\begin{center}
\begin{picture}(0,0)%
\includegraphics{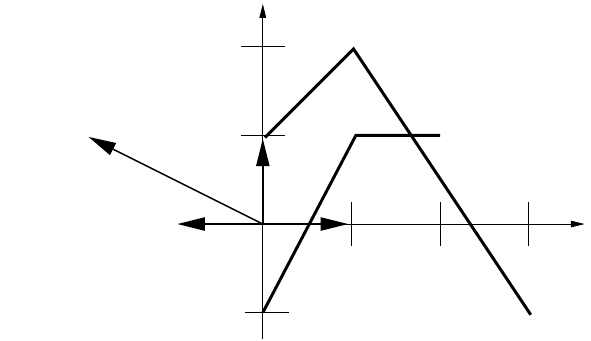}%
\end{picture}%
\setlength{\unitlength}{1865sp}%
\begingroup\makeatletter\ifx\SetFigFont\undefined%
\gdef\SetFigFont#1#2#3#4#5{%
  \reset@font\fontsize{#1}{#2pt}%
  \fontfamily{#3}\fontseries{#4}\fontshape{#5}%
  \selectfont}%
\fi\endgroup%
\begin{picture}(5967,3444)(5881,-7543)
\put(5896,-5866){\makebox(0,0)[lb]{\smash{{\SetFigFont{11}{8.4}{\rmdefault}{\mddefault}{\updefault}{\color[rgb]{0,0,0}$(-2,1)$}%
}}}}
\put(7156,-6721){\makebox(0,0)[lb]{\smash{{\SetFigFont{11}{8.4}{\rmdefault}{\mddefault}{\updefault}{\color[rgb]{0,0,0}$-1$}%
}}}}
\put(7561,-5326){\makebox(0,0)[lb]{\smash{{\SetFigFont{11}{8.4}{\rmdefault}{\mddefault}{\updefault}{\color[rgb]{0,0,0}$(0,1)$}%
}}}}
\put(9541,-6226){\makebox(0,0)[lb]{\smash{{\SetFigFont{11}{8.4}{\rmdefault}{\mddefault}{\updefault}{\color[rgb]{0,0,0}$1$}%
}}}}
\put(9046,-6946){\makebox(0,0)[lb]{\smash{{\SetFigFont{11}{8.4}{\rmdefault}{\mddefault}{\updefault}{\color[rgb]{0,0,0}$\sigma$}%
}}}}
\put(9811,-4696){\makebox(0,0)[lb]{\smash{{\SetFigFont{11}{8.4}{\rmdefault}{\mddefault}{\updefault}{\color[rgb]{0,0,0}$\rho$}%
}}}}
\end{picture}%

\end{center}
\vspace{-2mm}

\noindent Proposition~\ref{marius} reads in this case
\begin{align*}
\MI_1(\rho,\sigma)&= h_{[0,3]}(-1) \, \sigma(0) + h_{[0,3]}(1) \, \sigma(2) 
\\[1mm]
&\kern10pt+h_{Q_{\rho}}(-2,1) \, \Vol([0,1])+h_{Q_{\rho}}(0,1) \, \Vol([1,2])
\\[2mm]
&=  0+3+1+2=6
\enspace.
\end{align*}
The decomposition formula~(\ref{MI_decomposition}) can be a convenient
alternative 
for computing mixed integrals, since it avoids the costly calculation
of sup-convolutions. 

\bigskip


\end{document}